\documentclass{amsart}
\usepackage{amssymb, amscd}
\usepackage{pst-all}
\makeindex
\numberwithin{equation}{section}

\def\CC{{\mathbb C}}

\def\EE{{\mathbb E}}
\def\FF{{\mathbb F}}
\def\GG{{\mathbb G}}

\def\LL{{\mathbb L}}
\def\PP{{\mathbb P}}
\def\QQ{{\mathbb Q}}
\def\RR{{\mathbb R}}
\def\TT{{\mathbb T}}
\def\VV{{\mathbb V}}
\def\WW{{\mathbb W}}
\def\ZZ{{\mathbb Z}}

\def\Acal{{\mathcal A}}

\def\Ecal{{\mathcal E}}
\def\Fcal{{\mathcal F}}

\def\Hcal{{\mathcal H}}

\def\Lcal{{\mathcal L}}
\def\Mcal{{\mathcal M}}

\def\Xcal{{\mathcal X}}

\def\Hg{{\Hcal}_g}
\def\Hgmin{{\Hcal}_{g-1}}
\def\H#1{{\Hcal}_{#1}}
\def\SLZ{{\rm SL}(2,\ZZ)}
\def\SLR{{\rm SL}(2,\RR)}
\def\SpgZ{{\rm Sp}(2g,\ZZ)}

\def\SpgR{{\rm Sp}(2g,\RR)}

\def\symp{\langle \, , \, \rangle}
\def\ABCD{\left( \begin{matrix}A & B \cr C & D \cr \end{matrix} \right)}

\newcommand\Definition[1]{\emph{#1}}

\newcommand\proofsquare{\nobreak\hfill \hbox{%
\vrule height 5pt
\kern-.4pt
 \vbox{%
\hrule width 5pt depth0pt height.4pt
 \kern4.6pt \hrule  }
\kern-3.75pt
\vrule height 5pt}\kern1pt
\par}

\newtheorem{theorem}{Theorem}[section]
\newtheorem{lemma}[theorem]{Lemma}
\newtheorem{proposition}[theorem]{Proposition}
\newtheorem{corollary}[theorem]{Corollary}
\newtheorem{conjecture}[theorem]{Conjecture}
\newtheorem{definition-lemma}[theorem]{Definition-Lemma}

\theoremstyle{definition}
\newtheorem{definition}[theorem]{Definition}
\newtheorem{example}[theorem]{Example}

\theoremstyle{remark}
\newtheorem{remark}[theorem]{Remark}
\newtheorem{exercise}[theorem]{Exercise}

\newtheorem{question}[theorem]{Question}

\begin{document}

\title[Siegel Modular Forms]{
Siegel Modular Forms
\footnote{\tt siegel.tex\today}}
\author{Gerard van der Geer}

\address{Faculteit Wiskunde en Informatica, University of
Amsterdam, Plantage Muidergracht 24, 1018 TV Amsterdam, The Netherlands}

\email{geer@science.uva.nl}


\begin{abstract}
These are the lecture notes of the lectures on Siegel modular
forms at the Nordfjordeid Summer School on Modular Forms and their
Applications. We give a survey of Siegel modular forms and
explain the joint work with Carel Faber on vector-valued Siegel
modular forms of genus $2$ and present evidence for a conjecture
of Harder on congruences between Siegel modular forms of genus 
$1$ and $2$.
\end{abstract}

\maketitle

\begin{section}{Introduction}
\label{sec: intro}
Siegel modular forms generalize the usual modular forms on 
${\rm SL}(2,\ZZ)$ in that the group ${\rm SL}(2,\ZZ)$ is
replaced by the automorphism group ${\rm Sp}(2g,\ZZ)$
of a unimodular symplectic form on $\ZZ^{2g}$
and the upper half plane is replaced by the Siegel upper half
plane $\Hg$. The integer $g \geq 1$ is called
the degree or genus.

Siegel \index{Siegel} pioneered the generalization of the theory of
elliptic modular forms to the modular forms in more variables
now named after him. 
He was motivated by his work on the Minkowski-Hasse principle for
quadratic forms over the rationals, cf., \cite{Sie0}.
He investigated the geometry of the
Siegel upper half plane, determined a fundamental domain
and its volume and proved a central result equating an
Eisenstein series with a weighted sum of theta functions.

No doubt, Siegel modular forms are of fundamental importance
in number theory and algebraic geometry, but unfortunately, their
reputation does not match their importance. And although 
vector-valued rather than scalar-valued Siegel modular forms 
are the natural generalization of elliptic modular forms, 
their reputation amounts to even less. 
A tradition of ill-chosen notations may have contributed to this, 
but the lack of attractive examples that can be handled decently seems
to be the main responsible. Part of the beauty of elliptic
modular forms is derived from the ubiquity of
easily accessible examples. The accessible examples that we have of Siegel
modular forms are scalar-valued Siegel modular forms given by
Fourier series and for $g>1$ it is difficult to extract the
arithmetic information (e.g., eigenvalues of Hecke operators) from
the Fourier coefficients.

The general theory of automorphic representations provides a generalization
of the theory of elliptic modular forms. 
But despite the obvious merits of this approach 
some of the attractive explicit 
features of the $g=1$ theory are lost
in the generalization.

The elementary theory of elliptic modular forms ($g=1$) requires
little more than basic function theory, while a good grasp of the
elementary theory of Siegel modular forms requires a better understanding of
the geometry involved, in particular of  the compactifications of
the quotient space ${\rm Sp}(2g,\ZZ)\backslash \Hg$. A singular
compactification was provided by Satake and Baily-Borel 
and a smooth compactification
by Igusa in special cases and by Mumford c.s.\ by an intricate machinery
in the general case.

The fact that  ${\rm Sp}(2g,\ZZ)\backslash \Hg$ is the moduli space of
principally polarized abelian varieties plays an important role
in the arithmetic theory of modular forms. Even for $g=1$ one needs the 
understanding of the geometry of moduli space as a scheme (stack)
over the integers and its cohomology as Deligne's proof
of the estimate $|a(p)|\leq 2p^{(k-1)/2}$ (the Ramanujan conjecture)
for the Fourier coefficients
of a Hecke eigenform of weight $k$ showed. 
For quite some time the lack of a well-developed theory of moduli
spaces of principally polarized abelian varieties over the integers
formed a serious hurdle for the development of the
arithmetic theory. Fortunately, Faltings' work
on the moduli spaces of abelian varieties has provided us with the
first necessary ingredients of the arithmetic theory, both the smooth
compactification over $\ZZ$ as well as the Satake compactification 
over $\ZZ$. It also gives the analogue of the Eichler-Shimura
theorem which expresses Siegel modular forms in terms of the cohomology
of local systems on ${\rm Sp}(2g,\ZZ)\backslash {\mathcal H}_g$.
The fact that the vector-valued Siegel modular forms are the
natural generalization of the classical elliptic modular forms
becomes apparent if one studies the cohomology of the universal
abelian variety.

Examples of modular forms for ${\rm SL}(2,\ZZ)$ are easily constructed
using Eisenstein series or theta series.
These methods are much less effective when dealing with the case 
$g\geq 2$,
especially if one is interested in vector-valued Siegel modular 
forms. Some examples can be constructed using theta series, but it is 
not always easy to
calculate the Fourier coefficients and more difficult to extract
the eigenvalues of the Hecke operators.

We show that there is an alternative approach that uses the
analogue of the classical Eichler-Shimura theorem. Since cohomology
of a variety over a finite field can be calculated by determining
the number of rational points over extension fields one can
count curves over finite fields to calculate traces
of Hecke operators on spaces of vector-valued cusp forms for $g=2$.
This is joint work with Carel Faber. 
It has the pleasant additional feature that our forms all live in level $1$, i.e.\ on the full Siegel modular group.

We illustrate this by providing convincing evidence for a 
conjecture of Harder on congruences between the
eigenvalues of Siegel modular forms of genus $2$ and elliptic modular forms.

In these lectures we concentrate on modular forms for the full Siegel
modular group ${\rm Sp}(2g,\ZZ)$ and leave modular forms on
congruence subgroups aside. We start with the elementary theory
and try to give an overview of the various interesting aspects
of Siegel modular forms. An obvious omission are the Galois 
representations associated to Siegel modular forms.

A good introduction to the Siegel modular group and Siegel modular
forms is Freitag's book \cite{Fr1}. The reader may also consult the
introductory book by Klingen~\cite{Kl}.
Two other references to the literature are the two books
\cite{Shim3, Shim4} by Shimura. Vector-valued Siegel modular forms
are also discussed in a paper by Harris, \cite{Harris}.

{\sl Acknowledgements.}
I would like to thank Carel Faber,
Alex Ghitza, Christian Grundh, Robin de Jong, Winfried Kohnen,
Sam Grushevsky, Martin Weissman and Don Zagier for reading the manuscript and/or
 providing helpful comments. Finally I would like to thank Kristian Ranestad
for inviting me to lecture in Nordfjordeid in 2004.
\end{section}
\eject

\tableofcontents
\vfill
\eject

\begin{section}{The Siegel Modular Group}\label{SiegelModularGroup}
\bigskip
The ingredients of the definition of `elliptic modular form'
are the group $\SLZ$, the upper half plane ${\Hcal}$, the action of
$\SLZ$ on ${\Hcal}$, the concept of a holomorphic function
and the factor of automorphy $(cz+d)^k$. So if we want to
generalize the concept `modular form' we need to generalize
these notions.
But the upper half plane can be expressed in terms of the group
as $\SLR /{\rm SO(2)}$, where ${\rm SO}(2)={\rm U}(1)$, a maximal
compact subgroup, is the stabilizer
of the point $i=\sqrt{-1}$. Therefore, the group is the central
object and we start by generalizing
the group. The group $\SLZ$ is the
automorphism group of the lattice $\ZZ^2$ with the standard
alternating form $\langle \, , \, \rangle$ with
$$
\langle (a,b), (c,d)\rangle = ad-bc.
$$
This admits an obvious generalization by taking for $g \in \ZZ_{\geq 1}$
the lattice $\ZZ^{2g}$ of rank $2g$ with basis $e_1,\ldots, e_g,
f_1,\ldots,f_g$ provided with the symplectic form $\symp$
with
$$
\langle e_i,e_j \rangle=0, \, \langle f_i,f_j \rangle=0
\quad {\rm and}\quad \langle e_i,f_j \rangle =\delta_{ij},
$$
where $\delta_{ij}$ is Kronecker's delta. The {\sl symplectic group}
\index{symplectic group}
$\SpgZ$ is by definition the automorphism group of this symplectic lattice
$$
\SpgZ := {\rm Aut}(\ZZ^{2g}, \symp).
$$
By using the basis of the $e$'s and the $f$'s we can write
the elements of this group as matrices
$$
\left( \begin{matrix} A & B \cr C & D \cr
\end{matrix} \right),
$$
where $A$, $B$, $C$ and $D$ are $g\times g$ integral matrices
satisfying $AB^t=BA^t$, $CD^t=DC^t$ and $AD^t-BC^t=1_g$. Here we write
$1_g$ for the $g\times g$ identity matrix.
For $g=1$ we get back the group $\SLZ$. The group $\SpgZ$ is called the
\Definition{Siegel modular group} 
\index{Siegel modular group}
(of degree $g$) and often denoted $\Gamma_g$.

\begin{exercise}
Show that the conditions on $A,B,C$ and $D$ are equivalent to
$C^t\cdot A-A^t\cdot C=0$, $D^t\cdot B-B^t\cdot D=0$ and $D^t\cdot A -B^t\cdot C=1_g$.
\end{exercise}

The upper half plane ${\Hcal}$ can be given as a coset space
${\rm SL}(2,\RR) /K$ with $K={\rm U}(1)$ a maximal compact subgroup,
and this admits a generalization,
but the desired generalization also
admits a description as a half plane and with this we start:
the {\sl Siegel upper half plane} $\Hg$ 
\index{Siegel upper half plane} 
is defined as
$$
\Hg= \{ \tau \in {\rm Mat}(g\times g, \CC) 
\colon \, \tau^t=\tau, \, {\rm Im}(\tau) >0 \},
$$
consisting of $g\times g$ complex symmetric matrices which have
positive definite imaginary part (obtained by taking the imaginary
part of every matrix entry). Clearly, we have ${\Hcal}_1={\Hcal}$.

An element $\gamma=\ABCD$ of the
group $\SpgZ$, sometimes denoted by $(A,B;C,D)$,
 acts on the Siegel upper half
plane by
$$
\tau \mapsto \gamma(\tau)=(A\tau +B)(C\tau+D)^{-1}. \eqno(1)
$$
Of course, we must check that this is well-defined, in particular that
$C \tau +D$ is invertible. 
For this we use the identity
$$
(C\bar{\tau}+D)^t(A\tau+B)-(A\bar{\tau}+B)^t(C\tau +D)=\tau-\bar{\tau}=2iy,
\eqno(2)
$$
where we write $\tau=x+iy$
with $x$ and $y$ symmetric real $g\times g$ matrices. 
We claim that $\det(C\tau+D)\neq 0$. Indeed,
if the equation $(C\tau+D)\xi=0$ has a solution $\xi\in \CC^g$ then
equation (2) implies $\bar{\xi}^t y \xi =0$ and by the assumed
positive definiteness of $y$ that $\xi=0$. 

One can also check directly the identity
\begin{align}
(C\tau +D)^t(\gamma(\tau)-\gamma(\tau)^t)(C\tau +D)&=
(C\tau+D)^t(A\tau+B)-(A\tau+B)^t(C\tau+D)\cr
&=\tau -\tau^t=0\cr  \notag
\end{align}
that shows that $\gamma(\tau)$ is symmetric. Moreover, again by (2)
and this last identity
 we find the relation between  $y'={\rm Im}(\gamma(\tau))$ 
and $y$
$$
(C\bar{\tau}+D)^t y' (C\tau+D)=
\frac{1}{2i} (C\bar{\tau}+D)^t(\gamma(\tau)-(\overline{\gamma(\tau)})^t)
(C\tau+D)=y
$$
and this shows that $y'={\rm Im}(\gamma(\tau))$ is positive definite.
Using these details one easily checks that (1) defines indeed an action of
$\SpgZ$, and even of $\SpgR$ on $\Hg$. 

The group $\SpgR/\{\pm 1\}$ acts effectively on $\Hg$ and it is the 
biholomorphic automorphism group of $\Hg$. The action is transitive
and the stabilizer of $i\, 1_g$ is
$$
{\rm U}(g):=\{\left( \begin{matrix} A & B \cr -B & A\cr \end{matrix} \right)
 \in \SpgR \colon 
A\cdot A^t + B \cdot B^t = 1_g \},
$$
the unitary group.
We may thus view  $\Hg$ as the coset space
$\SpgR /{\rm U}(g)$  of a simple Lie group
by a maximal compact subgroup (which is unique up to conjugation).

The disguise of $\H1$ as the unit disc $\{ z \in \CC \colon |z| <1 \}$
also has an analogue for $\Hg$.
The space $\Hg$ is analytically equivalent to a bounded symmetric domain
$$
D_g := \{ Z \in {\rm Mat}(g \times g, \CC) \colon Z^t=Z, \, Z^t \cdot Z < 1_g \}
$$
and the generalized Cayley transform \index{Cayley transform}
$$
\tau \mapsto z = (\tau - i 1_g)(\tau + i 1_g)^{-1}, \qquad
z \mapsto \tau=i \cdot (1_g+z)(1_g-z)^{-1}
$$
makes the correspondence explicit. The `symmetric' in the name refers to the
existence of an involution on ${\mathcal H}_g$ (or $D_g$)
$$
\tau \mapsto -\tau^{-1} \qquad (z \mapsto -z)
$$
having exactly one isolated fixed point. Note that we can write $\Hg$ also as
$S_g + i S_g^+$ with $S_g$ (resp.\ $S_g^+$) the $\RR$-vector space (resp.\
cone) of real symmetric (resp.\ real positive definite symmetric) 
matrices of size $g \times g$.

The group 
$\SpgZ$ is a discrete subgroup of $\SpgR$ and acts properly discontinuously
on $\Hg$, i.e., for every $\tau \in \Hg$ there is an open  neighborhood $U$
of $\tau$
such that $\{ \gamma \in \SpgZ \colon \gamma(U) \cap U \neq \emptyset \}$
is finite. In fact, this follows immediately from the properness of the map
 $\SpgR \to \SpgR /U(g)$.

For $g=1$ usually one proceeds after these introductory remarks on the action
to the construction of a fundamental domain for the action of $\SLZ$ and all
the texts display the following archetypical figure.

$$
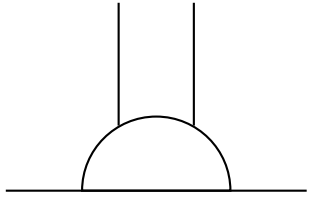
\setlength{\unitlength}{1mm}
\begin{picture}(50,30)
\psline{-}(2.5,0.866)(2.5,2.5)
\psline{-}(1.5,0.866)(1.5,2.5)
\psline{-}(0,0)(4,0)
\pswedge(2,0){1}{0}{180}
\end{picture}
$$

Siegel (see \cite{Sie1}) constructed also a fundamental domain for $g\geq 2$, 
\index{fundamental domain}
namely the set of $\tau=x+iy \in \Hg$ satisfying the following three 
conditions:
\begin{itemize}
\item $|\det (C\tau +D)| \geq 1$ for all $(A,B;C,D) \in \Gamma_g$;
\item $y$ is reduced in the sense of Minkowski;
\item the entries $x_{ij}$ of $x$ satisfy $|x_{ij}|\leq 1/2$. 
\end{itemize}
\noindent
Here Minkowski reduced means that $y$ satisfies the two properties 1)
$h^tyh\geq y_{kk}$ ($k=1,\ldots,g$) for all primitive 
vectors $h$ in $\ZZ^g$ and 2) $y_{k,k+1} \geq 0$ for $0\leq k \leq g-1$.
Already for $g=2$ the boundary of this fundamental domain is complicated;
Gottschling found that it posesses 28 boundary pieces, cf., \cite{Go},
and the whole thing does not help much to understand the nature of the 
quotient space $\SpgZ \backslash \Hg$.

The group $\SpgZ$ does not act freely on $\Hg$, but the subgroup 
$$
\Gamma_g(n):= 
\{ \gamma \in \SpgZ \colon \gamma \equiv 1_{2g} \, (\bmod \, n) \}
$$
acts freely if 
$n\geq 3$ as is easy to check, cf.\  \cite{Serre}.
 The quotient space (orbit space)
$$
Y_g(n):= \Gamma_g(n)\backslash \Hg
$$
is then for $n\geq 3$ a complex manifold of dimension $g(g+1)/2$.
Note that the finite group
${\rm Sp}(2g,\ZZ / n \ZZ)$ acts on $Y_g(n)$ as a group of biholomorphic
automorphisms and we can thus view
$$
\SpgZ \backslash \Hg
$$
as an orbifold (quotient of a manifold by a finite group).

The Poincar\'e metric on the upper half plane also generalizes
to the Siegel upper half plane. The corresponding volume form
is given by
$$
(\det y)^{-(g+1)} \prod_{i\leq j} dx_{ij} \, dy_{ij}
$$
which is $\partial \overline{\partial} \log \det {\rm Im}(\tau)^g$.
The volume of the fundamental domain was calculated by Siegel,
\cite{Sie2}. If we
normalize the volume such that it gives the orbifold Euler
characteristic 
\index{Euler characteristic}
the result is (cf. Harder \cite{Harder1})
$$
{\rm vol}(\SpgZ \backslash \Hg)= \zeta(-1)\zeta(-3) \cdots \zeta(1-2g).
$$
with $\zeta(s)$ the Riemann zeta function. In particular, for $n\geq 3$
the Euler number of the manifold $\Gamma_g(n)\backslash {\Hcal}_g$ equals
$[\Gamma_g(1):\Gamma_g(n)] \zeta(-1) \cdots \zeta(1-2g)$.
We first present two exercises for the solution of which we
refer to \cite{Fr1}.
\begin{exercise}
Show that the Siegel modular group $\Gamma_g$ is generated by the elements
$\left( \begin{matrix} 1_g & s \\ 0 & 1_g \\ \end{matrix}\right)$
with $s=s^t$ symmetric and the element $\left( \begin{matrix}
0 & 1_g \\ -1_g & 0 \\ \end{matrix} \right)$.
\end{exercise}
\begin{exercise} 
Show that ${\rm Sp}(2g,\ZZ)$ is contained in ${\rm SL}(2g,\ZZ)$.
\end{exercise}

We close with another model of the domain ${\mathcal H}_g$ that
can be obtained
as follows. Extend scalars of our symplectic lattice 
$(\ZZ^{2g},\langle \, , \, \rangle)$ to $\CC$ and let 
$Y_g$ be the Lagrangian Grassmann variety 
\index{Lagrangian Grassmann variety}
parametrizing totally isotropic subspaces of dimension $g$:
$$
Y_g:=\{ L \subset \CC^{2g} \colon \dim(L)=g, \, \langle x,y\rangle =0
\, \hbox{\rm for all $x,y \in L$}\}.
$$
Since the group ${\rm Sp}(2g,\CC)$ acts transitively on the set of
totally isotropic subspaces we may identify $Y_g$ with the compact
manifold ${\rm Sp}(2g,\CC)/Q$, where $Q$ is the parabolic subgroup
that fixes the first summand $\CC^g$. Consider now in $Y_g$ the open set
$Y_g^{+}$ of Lagrangian subspaces $L$ such that 
$-i\langle x, \bar{x}\rangle >0$
for all non-zero $x$ in $L$. Then $Y_g^{+}$ is stable under the action of
${\rm Sp}(2g,\RR)$ and the stabilizer of a point is isomorphic to the unitary group
$U(g)$. A basis of such an $L$ is given by the columns of a unique
$2g\times g$ matrix 
$\left( \begin{matrix} -1_g \\ \tau\\ \end{matrix} \right)$ 
with $\tau \in {\mathcal H}_g$ and this embeds
${\mathcal H}_g$ in $Y_g$ as the open subset $Y_g^{+}$; 
for $g=1$ we get the upper
half plane in $\PP^1$. The manifold $Y_g$ is called the compact dual
\index{compact dual}
of ${\mathcal H}_g$.
\begin{remark} Just as for $g=1$ we could consider congruence subgroups
of ${\rm Sp}(2g,\ZZ)$, like for example $\Gamma_g(n)$, the kernel
of the natural homomorphism
${\rm Sp}(2g,\ZZ) \to {\rm Sp}(2g,\ZZ/n\ZZ)$ for natural numbers $n$.
We shall stick to the full symplectic group ${\rm Sp}(2g,\ZZ)$ here.
\end{remark}
\end{section}
\begin{section}{Modular Forms}\label{modularforms}
\bigskip
To generalize the notion of modular form 
as we know it for $g=1$ we still have to
generalize the `automorphy factor' $(cz+d)^k$. To do this we consider
a representation
$$
\rho \colon {\rm GL}(g,\CC) \to {\rm GL}(V)
$$
with $V$ a finite-dimensional $\CC$-vector space.

For reasons that become clear later, 
it is useful to provide $V$ with a hermitian metric
$(\, ,\, )$ such that $(\rho(g) v_1,v_2)=(v_1,\overline{\rho({g}^t)}v_2)$
and we shall put $\| v\|= (v,v)^{1/2}$. Such a hermitian metric can always be
found and is unique up to a scalar for irreducible representations.

\noindent
\begin{definition} A holomorphic map $f \colon \Hg \to V$ is called
a {\sl Siegel modular form of weight}
\index{Siegel modular form} \index{weight} $\rho$ if
$$
f(\gamma(\tau))= \rho(C\tau+D) f(\tau)
$$
for all $\gamma = \left(\begin{matrix} A & B \cr C & D \cr \end{matrix}\right)
\in \SpgZ$ and all $\tau \in \Hg$, plus for $g=1$ the requirement
that $f$ is holomorphic at $\infty$.
\end{definition}
Before we proceed, a word about notations. The subject has been plagued with
unfortunate choices of notations, and the tradition of using capital
letters for the matrix blocks of elements of the symplectic group is one
of them. I propose to use lower case letters, so I will write
$f(\gamma(\tau))=\rho(c\tau+d)f(\tau)$ for all $\gamma =
(a,b;c,d) \in \Gamma_g$ for our condition.
 
The modular forms we consider here are vector-valued modular forms.
As it turns out, the holomorphicity condition is not necessary for 
$g>1$, see the Koecher principle hereafter. 

Modular forms of weight $\rho$ form a $\CC$-vector space
$M_{\rho}=M_{\rho}(\Gamma_g)$ and we shall see later that 
all the $M_{\rho}$ are finite-dimensional. 
If $\rho$  is a direct sum of two representations
$\rho=\rho_1\oplus \rho_2$ then $M_{\rho}$
is isomorphic to the direct sum $M_{\rho_1}\oplus M_{\rho_2}$
and this allows us to restrict ourselves 
to studying $M_{\rho}$ for the irreducible 
representations of ${\rm GL}(g,\CC)$.  

As is well-known (see \cite{F-H}, but see also the later Section \ref{roots}),
the irreducible finite-dimensional
representations of ${\rm GL}(g,\CC)$ correspond bijectively to
the $g$-tuples $(\lambda_1,\ldots,\lambda_g)$
of integers with $\lambda_1 \geq \lambda_2 \geq \cdots \geq \lambda_g$,
the highest weight of the representation $\rho$.
That is, for each irreducible $V$
there exists a unique $1$-dimensional subspace 
$\langle v_{\rho} \rangle$ of $V$ such that $\rho({\rm diag}(a_1,\ldots,a_g))$
acts on $v_{\rho}$ by multiplication by $\prod_{i=1}^g a_i^{\lambda_i}$.
For example, the $g$-tuple $(1,0,\ldots,0)$ corresponds to the
tautological representation $\rho(x)=x$ for $x\in {\rm GL}(g,\CC)$,
while the determinant representation corresponds to 
$\lambda_1=\ldots=\lambda_g=1$. Tensoring a given irreducible representation
with the $k$-th power of the determinant changes the $\lambda_i$ to
$\lambda_i+k$. We thus can arrange that $\lambda_g=0$
or that $\lambda_g \geq 0$ (i.e.\ that the representation is 
`polynomial').
Let $R$ be the set of isomorphism classes of representations of 
${\rm GL}(g,\CC)$. This set forms a ring with $\oplus$ as addition 
and $\otimes$ as multiplication. It is called the representation
ring of ${\rm GL}(g,\CC)$.

For $g=1$ one usually forms a graded ring of modular forms by taking
 $M_{*}(\Gamma_1)=\oplus M_k(\Gamma_1)$.
We can try do something similar for $g>1$ and try to make the direct sum
$\oplus_{\rho \in R} M_{\rho}(\Gamma_g)$ into a graded ring. But of course, this
is a huge ring, even for $g=1$ much larger than $M_*(\Gamma_1)$ since it 
involves also the reducible representations and it is not really what we want.

The classes of the irreducible representations of  ${\rm GL}(g,\CC)$
form a subset of all classes of representations. For $g=1$ and $g=2$ 
the fact is that the tensor product of two irreducible representations 
is a direct sum of irreducible representations with multiplicity $1$. 
In fact, for $g=1$ the tensor product of the irreducible
representations $\rho_{k_1}$
and $\rho_{k_2}$ of degree $k_1+1$ and $k_2+1$
is the irreducible representation
$\rho_{k_1+k_2}$. For $g=2$, a case that will play a prominent role in
these lecture notes,
we let $\rho_{j,k}$ denote the irreducible
representation of ${\rm GL}(2,\CC)$ that is ${\rm Sym}^j(W) \otimes \det(W)^k
$ with $W$ the standard $2$-dimensional representation; it corresponds
to highest weight $(\lambda_1,\lambda_2)=(j+k,k)$. Then there is
the formula
$$
\rho_{j_1,k_1} \otimes \rho_{j_2,k_2} \cong \sum_{r=0}^{\min (j_1,j_2)}
 \rho_{j_1+j_2-2r,k_1+k_2+r}.
$$
So we can decompose $M_{\rho_{j_1,k_1}}$ as a direct sum 
$\sum_{r=0}^{\min (j_1,j_2)} M_{\rho_{j_1+j_2-2r,k_1+k_2+r}}$,
but this is not canonical as it requires 
a choice of isomorphism in the above formula. Nevertheless,
this decomposition is useful to construct modular forms in new weights
by multiplying modular forms.

To make $\oplus_{\rho \in {\rm Irr}}M_{\rho}(\Gamma_2)$ into a ring
requires a consistent 
choice for all these identifications. We can avoid this by remarking
that multiplication of polynomials defines a canonical map
${\rm Sym}^{j_1}(W) \otimes {\rm Sym}^{j_2}(W) \to
{\rm Sym}^{j_1+j_2}(W)$. Using this and the obvious map
$\det(W)^{k_1}\otimes \det(W)^{k_2} \to \det(W)^{k_1+k_2}$ the direct
sum $\oplus_{\rho \in {\rm Irr}}M_{\rho}(\Gamma_2)$
becomes a ring; \index{ring of modular forms}
we just `forgot' the terms in the above sum with $r>0$.
For $g\geq 3$ the tensor products come in general with multiplicities,
given by Littlewood-Richardson numbers. Nevertheless, one
can define a ring structure on $\oplus_{\rho \in {\rm Irr}}M_{\rho}(\Gamma_g)$
that extends the multiplication of modular forms for $g=1$ and
the one given here for $g=2$ as Weissman shows. 
We refer to his interesting paper, \cite{weissman}.

For every $g$ one obtains a subring of the representation ring by 
taking the powers of the determinant
$\det : {\rm GL}(g,\CC) \to \CC^*$. This leads to a ring of 
`classical' modular
forms.
\smallskip
\noindent
\begin{definition} A {\sl classical Siegel modular form of weight}
 $k$ (and degree $g$) \index{classical Siegel modular form}
is a holomorphic function $f \colon \Hg \to \CC$ such that
$$f(\gamma(\tau))=\det (c\tau +d)^k f(\tau)
$$
for all $\gamma=(a,b;c,d) \in \SpgZ$ 
(with for $g=1$ the usual holomorphicity requirement at $\infty$).
\end{definition}

Classical Siegel modular forms are also known as scalar-valued Siegel 
modular forms.

Let $M_k=M_k(\Gamma_g)$ be the vector space of classical Siegel modular
forms of weight $k$. Together these spaces 
form a graded ring $M^{\rm cl}:=\oplus M_k$ of $M$
of classical Siegel modular forms. Of course, for $g=1$ the notion
of classical modular form 
reduces to the usual notion of modular form on ${\rm SL}(2,\ZZ)$.

\end{section}
\begin{section}{The Fourier Expansion of a Modular Form}\label{Fourier}

The classical Fourier expansion of a modular form on $\SLZ$
has an analogue. To define it we need the following definition.

\begin{definition} A symmetric $g\times g$-matrix $n\in {\rm GL}(g,\QQ)$ 
is called  {\sl half-integral} \index{half-integral}
if $2n$ is an integral matrix the diagonal entries of which are even.
\end{definition}

Every half-integral $g \times g$-matrix $n$ 
defines a linear form with integral coefficients
in the coordinates
$\tau_{ij}$ with $ 1\leq i \leq j \leq g$
of $\Hg$, namely
$$
{\rm Tr}(n \tau)=\sum_{i=1}^g  n_{ii}\tau_{ii} + 
2\sum_{1 \leq i < j \leq g} n_{ij} \tau_{ij}
$$
and every linear integral combination of the coordinates is of this form.

Let us write $\tau = x+iy$ with $x$ and $y$ symmetric $g\times g$ matrices.
A function $f: \Hg \to \CC$ that is periodic in the sense that
 $f(\tau+s)=f(\tau)$ for all integral
symmetric $g\times g$-matrices $s$ admits a Fourier expansion
\index{Fourier expansion}
$$
f(\tau)=\sum_{n \, {\rm half-integral}} a(n) e^{2\pi i {\rm Tr}(n\tau)}
$$
with $a(n) \in \CC$
given by
$$
a(n)= \int_{x \, \bmod \, 1} f(\tau) e^{-2\pi i {\rm Tr}(n\tau)} dx
$$
with $dx$ the Euclidean volume of the space of $x$-coordinates
and the integral runs over $-1/2 \leq x_{ij} \leq 1/2$. 
This is a series which is uniformly convergent on compact subsets.
If $f$ is a vector-valued modular form  in $M_{\rho}$ we have a 
similar \Definition{Fourier series}
$$
f(\tau)=\sum_{n \, {\rm half-integral}} a(n) e^{2\pi i {\rm Tr}(n\tau)}
$$
with $a(n) \in V$. One could also use the suggestive notation
$$
f(\tau)=\sum_{n \, {\rm half-integral}} a(n) \, q^n,
$$
where we write $q^n$ for $e^{2\pi i {\rm Tr}(n\tau)}$.
Moreover, we have the property
$$
a(u^tnu)=\rho(u^t) \, a(n) \qquad \hbox{\rm for all $u \in {\rm GL}(g,\ZZ)$}.
\eqno(4)
$$
Indeed, we have
\begin{align}
a(u^tnu)&=\int_{x \, \bmod 1} f(\tau)e^{-2\pi i {\rm Tr}(u^tnu \tau)} dx\cr
&=\rho(u^t) \int_{x \, \bmod 1}
 f(u\tau u^t)e^{-2 \pi i {\rm Tr}(n \, u\tau u^t)} dx\cr
&= \rho(u^t) a(n).\cr \notag \end{align}
A direct corollary of formula (4) 
(proof left to the reader) restricts the
weight of non-zero forms.

\begin{corollary} A classical Siegel modular form of weight $k$ with
$kg \equiv 1 (\bmod \, 2)$ vanishes.
\end{corollary}

A basic result is the following theorem.
\begin{theorem}
Let $f \in M_{\rho}(\Gamma_g)$. Then $f$ is bounded
on any subset of $\Hg$ of the form
 $\{ \tau \in \Hg \colon {\rm Im}(\tau)> c\cdot 1_g\}$ with $c>0$.
\begin{proof} For $g=1$ the boundedness comes from the requirement in the
definition that the Fourier expansion $f= \sum_n a(n) q^n$ has no negative
terms.
So suppose that $g\geq 2$ and let $f =\sum_n a(n) e^{2\pi i{\rm Tr} n \tau} 
\in M_{\rho}(\Gamma_g)$. 
Since $f$ converges absolutely on $\Hg$ we see by substitution of $\tau
= i \cdot 1_g$ that there exists a constant $c>0$ such that for all 
half-integral matrices we have $|a(n)| \leq c e^{2\pi {\rm Tr} n \tau}$.
We first will show that $a(n)$ vanishes for $n$ that are not positive
semi-definite. 

Suppose that $n$ is not positive semi-definite. Then there exists a primitive
integral 
(column) vector $\xi$ such that $\xi^t n \xi <0$. We can complete $\xi$
to a unimodular matrix $u$. Using the relation $a(u^tnu)=\rho(u^t)a(n)$
and replacing $n$ by $u^t n u$ we may assume that entry $n_{11}$ of $n$
is negative. Consider now for $m \in \ZZ$ the matrix
$$
v=\left( 
\begin{matrix} 
1 & m & \cr 0 & 1 & \cr
&& 1_{g-2} & \cr 
\end{matrix} 
\right) 
\in {\rm GL}(g,\ZZ),
$$
where the omitted entries are zero.
We have 
$$
|a(n)|=|\rho(v^t)^{-1}| \, |a(v^tnv)| \leq c e^{2\pi {\rm Tr} v^tnv}.
$$
But ${\rm Tr}(v^tnv)={\rm Tr}(v)+ n_{11}m^2+2n_{12}m$ and if $m \to \infty$
then this expression goes to $-\infty$, so $|a(n)|=0$.

We conclude that $f=\sum_{n\geq 0} a(n)e^{2\pi i{\rm Tr} n \tau}$.
We can now majorize by the value at $c \, i \cdot 1_g$ of $f$, viz.\
$\sum_{n\geq 0} |a(n)|e^{-2\pi{\rm Tr}nc}$, uniformly in $\tau$ on
$\{\tau \in \Hg \colon {\rm Im}(\tau) > c \cdot 1_g \}$. 
\end{proof}
\end{theorem}
The proof of this theorem shows the validity of the so-called Koecher principle
announced above.

\begin{theorem}{\rm (Koecher Principle)} \index{Koecher principle}
Let $f=\sum_n a(n) q^n \in M_{\rho}(\Gamma_g)$ with $q^n=
e^{2\pi i{\rm Tr}(n\tau)}$ be a modular form of weight $\rho$.
Then $a(n)=0$ if the half-integral matrix $n$ is not positive semi-definite.
\end{theorem}
The Koecher principle was first observed in 1928 
by G\"otzky for Hilbert modular forms and in general by Koecher in
1954, see \cite{Koe}. 
\begin{corollary}
A classical Siegel modular form of negative weight vanishes.
\begin{proof}
Let $f\in M_k(\Gamma_g)$ with $k<0$. 
Then the function $h=\det(y)^{k/2}|f(\tau)|$ is invariant under 
$\Gamma_g$ since  ${\rm Im}(\gamma(\tau))=
(c\tau+d)^{-t}({\rm Im}(\tau))\overline{(c\tau+d)}^{-1}$. 
It is not difficult to see that a fundamental domain is contained in
$\{ \tau \in \Hg \colon
{\rm Tr}(x^2)< 1/c, \, y > c \cdot 1_g \}$ for some suitable $c$.
This implies that for negative $k$ the expression $\det(y)^{k/2}$
is bounded on a fundamental domain,  and by the Koecher principle
$f$ is bounded on
$\{ \tau \in \Hg \colon \det y \geq c \}$.
It follows that $h$ is bounded on $\Hg$, say $h \leq c'$ and with
$$
a(n) e^{-2\pi {\rm Tr} ny}=\int_{x \, \bmod \, 1} f(\tau)e^{-2\pi {\rm Tr} nx}
dx
$$
we get
$$
|a(n)|e^{-2\pi {\rm Tr}ny}=\sup_{x \, \bmod \, 1} |f(x+iy)| \leq c'
 \, \det y^{-k/2}.
$$
If we let $y \to 0$ then for $k<0$ we see $|a(n)|=0$ for all $n\geq 0$.
\end{proof}
\end{corollary}
This corollary admits a generalization for vector-valued Siegel modular
forms, cf., \cite{Fr3}:
\begin{proposition} Let $\rho$ be a non-trivial irreducible representation of
${\rm GL}(g,\CC)$ with highest weight $\lambda_1\geq \ldots \geq \lambda_g$.
If $M_{\rho}\neq \{ 0 \}$ then we have $\lambda_g \geq 1$.
\end{proposition}
One proves this by taking a totally real field $K$ of degree $g$ over $\QQ$
and by identifying the symplectic space $O_K\oplus O_K^{\vee}$
(with $O_K^{\vee}$ the dual of $O_K$ with respect to the trace)
with our standard symplectic space $(\ZZ^{2g},\langle \, , \, \rangle )$.
This induces an embedding ${\rm SL}(2,O_K) \to {\rm Sp}(2g,\ZZ)$
and a map ${\rm SL}(2,O_K)\backslash {\mathcal H}_1^g \to
{\rm Sp}(2g,\ZZ)\backslash {\mathcal H}_g$. Pulling back Siegel
modular forms yields Hilbert modular forms on ${\rm SL}(2,O_K)$.
Now use that a Hilbert modular form of weight $(k_1,\ldots,k_g)$
vanishes if one of the $k_i\leq 0$, cf., \cite{vdG2}. By varying $K$
one sees that if $\lambda_g \leq 0$
then a non-constant $f$ vanishes on a dense subset of ${\mathcal H}_g$.

\end{section}
\begin{section}{The Siegel Operator and Eisenstein Series}\label{siegeloperator}
Since modular forms $f \in M_{\rho}(\Gamma_g)$ are bounded in the sets
of the form 
$\{ \tau \in {\Hcal}_g  \colon {\rm Im}(\tau)> c\cdot 1_g\}$ 
we can take the limit.

\begin{definition} We define an operator $\Phi$ on $M_{\rho}i(\Gamma_g)$ by 
$$
\Phi f=
\lim_{t \to \infty} f(\begin{matrix}\tau^{\prime} & 0 \cr 0 & it\cr 
\end{matrix})
\qquad \hbox{\rm with $\tau^{\prime} \in \Hgmin, t \in \RR$}.
$$
\end{definition}
In view of the convergence
we can also apply this limit to
all terms in the Fourier series and get
$$
(\Phi f)(\tau^{\prime})= \sum_{n^{\prime} \geq 0} 
a(\begin{matrix} n^{\prime} & 0 \cr 0 & 0 \cr \end{matrix} )
e^{2 \pi i{\rm Tr}(n^{\prime} \tau^{\prime})}.
$$
The values of $\Phi f$ generate a subspace $V' \subseteq V$ that is 
invariant under the action of 
the subgroup of matrices $\{(a,0;0,1) 
\colon a \in {\rm GL}(g-1,\CC)\}$ and that defines a representation 
$\rho^{\prime}$ of ${\rm GL}(g-1,\CC)$.
The operator $\Phi$ defined on Siegel modular forms of degree $g$
is called the \Definition{Siegel operator} 
\index{Siegel operator}
and defines a linear
map $M_{\rho}(\Gamma_g) \to M_{\rho'}(\Gamma_{g-1})$.
If $\rho$ is the irreducible representation with highest weight
$(\lambda_1,\ldots,\lambda_g)$ then $\Phi$ maps $M_{\rho}(\Gamma_g)$
to $M_{\rho'}(\Gamma_{g-1})$ with $\rho'$ the irreducible representation
of ${\rm GL}(g-1,\CC)$ with highest weight $(\lambda_1,\ldots,\lambda_{g-1})$.
\smallskip
\noindent
\begin{definition} A modular form $f \in M_{\rho}$ is called a
{\sl cusp form} if $\Phi f=0$. The subspace of $M_{\rho}$ of cusp
forms is denoted by $S_{\rho}=S_{\rho}(\Gamma_g)$. \index{cusp form}
\end{definition}
\begin{exercise} Show that  a modular $f=\sum a(n) e^{2\pi i{\rm Tr}(a \tau)}
 \in M_{\rho}$  is a cusp form if and only if $a(n)=0$ for all semi-definite
$n$ that are not definite.
\end{exercise}

We can apply the Siegel operator repeatedly (say $r\leq g$ times)
to a Siegel modular form on $\Gamma_g$ and one thus obtains a Siegel 
modular form on $\Gamma_{g-r}$. If $\rho$ is irreducible with highest weight
$(\lambda_1,\ldots,\lambda_g)$ and $\Phi F=f\neq 0$ for some 
$F \in M_{\rho}(\Gamma_g)$ then necessarily $\lambda_g\equiv 0 (\bmod \, 2)$
because with $\gamma$ also $-\gamma$ lies in $\Gamma_g$.

Let now $f_1$ and $f_2$ be modular forms of weight $\rho$, one of them
a cusp form. Then we define the {\sl Petersson product} 
\index{Petersson product} of $f_1$ and $f_2$ by
$$
\langle f_1 , f_2 \rangle = \int_F (\rho({\rm Im}(\tau))f_1(\tau), 
f_2(\tau)) d\tau,
$$
where $d\tau= \det(y)^{-(g+1)}\prod_{i\leq j} dx_{ij}dy_{ij}$ is an 
invariant measure on $\Hg$, $F$ is a fundamental domain for the action
of $\Gamma_g$ on $\Hg$ and the brackets  $(\, , \, )$
refer to the Hermitian product defined in 
Section \ref{modularforms}. One checks that it converges exactly because 
at least one of the two forms is a cusp form.
Furthermore, we define
$$
N_{\rho}= S_{\rho}^{\bot},
$$
for the orthogonal complement of $S_{\rho}$ and then have an orthogonal
decomposition $M_{\rho}=S_{\rho}\oplus N_{\rho}$.

Just as in the case $g=1$ one can construct modular forms
explicitly using Eisenstein series. \index{Eisenstein series}
We first deal with the case of classical Siegel modular forms.
Let $g\geq 1$ be the degree and let $r$ be a natural number 
with $0\leq r \leq g$.
Suppose that $f\in S_k(\Gamma_r)$ is a (classical Siegel modular)
cusp form of {\sl even} weight $k$.
For a matrix $ \left( \begin{matrix} \tau_1 & z \\ z & \tau_2 \\ \end{matrix} 
\right)$ 
with $\tau_1 \in {\mathcal H}_r$ and $\tau_2 \in {\mathcal H}_{g-r}$ we write
$ \tau^{*}=\tau_1 \in {\Hcal}_r$. (For $r=0$ we let
$\tau^*$ be the unique point of $
{\Hcal}_0$.) If $k$ is positive and even we define the
\Definition{Klingen Eisenstein series}, 
\index{Klingen Eisenstein series} a formal series,
$$
E_{g,r,k}(f):= \sum_{A=(a,b;c,d) \in P_r\backslash \Gamma_g}
f((a\tau +b)(c\tau +d)^{-1})^*) \det(c\tau+d)^{-k},
$$
where $P_r$ is the subgroup
$$
P_r :=\left\{  \left(  
\begin{matrix} a' & 0 & b' & * \\
* & u & * & * \\
c' & 0 & d' & * \\
0 & 0 & 0 & u^{-t} \\
\end{matrix} \right) \in \Gamma_g : 
\left( \begin{matrix} a' & b ' \\ c' & d' \\ \end{matrix}\right)
\in \Gamma_r, \,
u \in {\rm GL}(g-r, \ZZ) \right\}.
$$
For an interpretation of this subgroup we refer to Section 
\ref{compactifications}. 
In case $r=0$, $f$ constant, say $f=1$, we get the old Eisenstein series 
$$
E_{g,0,k}= \sum_{(a,b;c,d)} \det(c\tau +d)^{-k},
$$
where the summation is over a full set of
representatives for the cosets
${\rm GL}(g,\ZZ) \backslash \Gamma_g$.

\begin{theorem} Let $g\geq 1$ and $0 \leq r \leq g$ and $k> g+r+1$ be integers
with $k$ even. Then for every cusp form $f \in S_k(\Gamma_r)$ the series
$E_{g,r,k}(f)$ converges to a classical Siegel modular form of weight $k$
in $M_k(\Gamma_g)$ and $\Phi^{g-r}E_{g,r,k}(f)=f$.
\end{theorem}

This theorem was proved by Hel Braun\footnote{Hel Braun was a student of 
Carl Ludwig Siegel (1896--1981)
\index{Braun}\index{Siegel},
the mathematician after whom our modular forms are named. She sketches
 an interesting portrait
of Siegel in \cite{Br}} in 1938 for $r=0$ and $k>g+1$.

The Fourier coefficients of these Eisenstein series were determined by
Maass, see \cite{Maa1}.
Often we shall restrict the summation
over co-prime $(c,d)$ in order to avoid an unnecessary factor.

\begin{corollary} The Siegel operator 
$\Phi: M_k(\Gamma_g) \to M_k(\Gamma_{g-1})$ is surjective for even $k>2g$.
\end{corollary}

Weissauer improved the above result and proved
that $\Phi$ is surjective if $k>(g+r+3)/2$, see \cite{Weis2}.
He also treated the case of vector-valued modular forms and showed that
the image  $\Phi(M_{\rho}(\Gamma_g))$ contains the space of cusp forms
$S_{\rho'}(\Gamma_{g-1})$ if $k=\lambda_g\geq g+2$.

If $k$ is odd we have no good Eisenstein series; for example look at the Siegel
operator $M_k(\Gamma_g) \to M_k(\Gamma_{g-1})$ for $k\equiv g \equiv 1 \, 
(\bmod \, 2)$. Then $M_k(\Gamma_g)=(0)$ while the target space 
$M_k(\Gamma_{g-1})$ is non-zero for sufficiently large $k$ 
(e.g.\ $M_{35}(\Gamma_2) \neq (0)$ as we shall see later).

Just as for $g=1$ one can construct Poincar\'e series \index{Poincar\'e series}
and use these to generate the spaces of cusp forms
if the weight is sufficiently high. These Poincar\'e series
behave well with respect to the Petersson product.
We refer to \cite{Kl}, Ch.\ 6,
or \cite{B-B} for the general setting.
\end{section}

\begin{section}{Singular Forms}
\label{singularforms} \index{singular modular form}

A particularity of $g>1$ are the so-called singular modular forms.

\begin{definition} A modular form $f=\sum_n a(n) e^{2 \pi i{\rm Tr} n\tau}
 \in M_k(\Gamma_g) $ is called \Definition{singular} 
if $a(n)\neq 0$ implies that $n$ is a singular matrix ($\det(n)=0$).
\end{definition}

Modular forms of small weight are singular as the following theorem shows,
\cite{Weis1}.

\begin{theorem} {\rm (Freitag, Salda{\~n}a, Weissauer )} 
Let $\rho$ be 
irreducible with highest weight $(\lambda_1,\ldots,\lambda_g)$.
A non-zero modular form 
$f \in M_{\rho}(\Gamma_g)$ is singular if and only if 
$2\lambda_g < g$.
\end{theorem} 
In particular, there are no cusp forms of weight $2\lambda_g < g$.
One defines the \Definition{co-rank} of an irreducible representation
as $\# \{ 1\leq i \leq g \colon \lambda_i=\lambda_g \}$. For a modular
form $f=\sum_n a(n) \exp(2 \pi i {\rm Tr} n \tau) \in M_{\rho}(\Gamma_g)$,
Weissauer introduced the \Definition{rank} and \Definition{co-rank}
of $f$ by
$$
{\rm rank}(f)= \max \{ {\rm rank}(n) \colon a(n)\neq 0\} 
$$
and
$$
{\rm co-rank}(f)=g -\min \{ {\rm rank}(n) \colon a(n)\neq 0\}.
$$
In particular, modular forms of rank $<g$ are singular
while cusp forms have co-rank~$0$ \index{rank of a modular form}
\index{co-rank of a modular form} and Siegel-Eisenstein
forms $E_{g,0,k}$ have co-rank~$g$; $\Phi$ applied $k+1$ 
times to forms of co-rank $k$ should be zero. 
Weissauer proved  (see \cite{Weis1})
for irreducible $\rho$ that
$\hbox{\rm co-rank}(f) \leq \hbox{\rm co-rank}(\rho)$
and also that $M_{\rho}(\Gamma_g)=(0)$ if $\lambda_g 
\leq g/2 - \hbox{\rm co-rank}(\rho)$.
More precisely, he proved
\begin{theorem}\label{Verschwindungssatz}
Let $\rho=(\lambda_1,\ldots,\lambda_g)$ be an irreducible
representation of co-rank $< g-\lambda_g$. If
$\#\{i \colon 1\leq i \leq g, \lambda_i=\lambda_g+1 \}
< 2(g-\lambda_g-\hbox{\rm co-rank}(\rho))$ then $M_{\rho}=(0)$.
\end{theorem}

Finally, Duke and Imamo{\v g}lu prove in \cite{D-I2} that there are no 
cusp forms of small weights; \index{small weights}
for example, $S_{6}(\Gamma_g)=(0)$
for all $g$.
\end{section}
\begin{section}{Theta Series}\label{thetaseries}\index{theta series}

Besides Eisenstein series one can construct Siegel modular forms 
using theta series. We begin with the so-called 
\Definition{theta-constants}. \index{theta constant}
Let $\epsilon=\left(
\begin{matrix}
\epsilon' \\ 
\epsilon^{\prime \prime}\\
\end{matrix}
\right)$ with 
$\epsilon', \epsilon^{\prime \prime} \in \{0,1\}^g$ 
and consider the rapidly converging
series
$$
\theta[\epsilon]=\sum_{m \in \ZZ^g} \exp 2 \pi i \{
(m+\frac{1}{2}\epsilon')^t\tau(m+\frac{1}{2}\epsilon') + \frac{1}{2}
(m+\frac{1}{2}\epsilon')^t (\epsilon^{\prime \prime})\}.
$$
This vanishes identically if 
$\epsilon$ is odd, that is, if
$\epsilon' (\epsilon^{\prime \prime})^t$ is odd. The other
$2^{g-1}(2^g+1)$ cases (the `even' ones) yield the so-called even 
theta characteristics. These are modular forms on a level 2 congruence subgroup
of ${\rm Sp}(2g,\ZZ)$ of weight $1/2$, cf.\ \cite{Ig5}. These can be used
to construct classical Siegel modular forms on ${\rm Sp}(2g,\ZZ)$.
For example, for $g=1$ one has
$$
(\theta[\begin{matrix} 0 \\ 0 \\ \end{matrix}]\theta [
\begin{matrix} 0 \\ 1 \\ \end{matrix}]\theta [
\begin{matrix} 1 \\ 0 \\ \end{matrix}])^8 = 2^8 \Delta \, \in S_{12}(\Gamma_1).
$$
For $g=2$ the product $-2^{-14}\prod \theta[\epsilon]^2$ of the squares of the
ten even theta characteristics gives a cusp form $\chi_{10}$ of weight $10$
on ${\rm Sp}(4,\ZZ)$, cf.\ \cite{Ig1,Ig2,Ig3}. Similarly, an expression
$$
(\prod \theta[\epsilon]) \sum \pm (\theta[\epsilon_1] \theta[\epsilon_2]
\theta[\epsilon_3])^{20},
$$
where the product is over the even theta characteristics and the sum 
is over so-called azygous triples of theta characteristics (i.e., triples
such that $\epsilon_1+\epsilon_2+\epsilon_3$ is odd) defines (up to a
normalization $-2^{-39}5^{-3} i$ ?) a cusp form $\chi_{35}$ of weight $35$
on ${\rm Sp}(4,\ZZ)$. Similarly, for $g=3$ the product of the $36$ even
theta characteristics defines a cusp form of weight $18$ on ${\rm Sp}(6,\ZZ)$.
The reason why one needs such a complicated expression is that the
theta characteristics are modular forms on a subgroup $\Gamma_g(4,8)$ of
${\rm Sp}(2g,\ZZ)$ and the quotient group ${\rm Sp}(2g,\ZZ)/\Gamma_g(4,8)$
permutes them and creates signs in addition so that we need a sort of
symmetrization to get something invariant.

Another source of Siegel modular forms are theta series associated to 
even unimodular lattices. Let $B$ be a positive definite symmetric 
even unimodular matrix of size $r \equiv \, 0 (\bmod \,  8)$. We denote
by $H_k(r,g)$ the space of harmonic polynomials $P: \CC^{r \times g} \to
\CC$ satisfying for $M \in {\rm GL}(g,\CC)$ the identity $P(zM)=
\det(M)^k P(z)$. Recall that harmonic means that $\sum_{i,j} \partial^2 /
\partial z_{ij}^2 \, P(z)=0$ if $z_{ij}$ are the coordinates 
on $\CC^{r \times g}$.
For a pair $(B,P)$ with $P\in H_k(r,g)$ we set
$$
\theta_{B,P}(\tau)= \sum_{A \in \ZZ^{r \times g}} 
P(\sqrt{B} A) e^{\pi i {\rm Tr}(A^t B A \tau )},
$$
where $\sqrt{B}$ is a positive matrix with square $B$.
Then $\theta_{B,P}$ is a classical 
Siegel modular form in $M_{k+r/2}(\Gamma_g)$,
see \cite{Fr1}. Such theta series for $P \in H_{k-r/2,g}$ and $B$
as above span  a subspace of $M_k(\Gamma_g)$ that is invariant
under the Hecke-operators that will be introduced later, 
cf.\ Section \ref{Heckeoperators}.
There are analogues of these that give vector-valued Siegel modular
forms if we require that $P$ is a vector-valued polynomial satisfying
the relation $P(zM)=\rho(M)P(z)$. See also Section \ref{vvgenus2} and
\cite{Ibu1, Ibu2}
for an example.

Finally, we would like to make a reference to
 Siegel's Hauptsatz \cite{Sie0} \index{Siegel's Hauptsatz}
(or \cite{Fr1}, p.\ 285)
on representations of quadratic forms by quadratic forms which can be
viewed as an identity between an Eisenstein series and a weighted
sum of theta series,  and to its far-reaching generalizations, cf.\
\cite{KR}.
\end{section}
\begin{section}{The Fourier-Jacobi Development 
of a Siegel Modular Form}\label{Fourier-Jacobi}
As we saw above, just as for $g=1$ we have a Fourier expansion of
a Siegel modular form $f= \sum_{n\geq 0} a(n) e^{2 \pi i {\rm Tr}(n\tau)}$.
But for $g>1$ there are other developments that provide more information,
like the so-called Fourier-Jacobi development, a concept due to
Piatetski-Shapiro.
 
We consider classical Siegel modular forms of weight $k$ on $\Gamma_g$.
We write $\tau \in \Hg$ as 
$$
\tau =\left( \begin{matrix} \tau^{\prime} & z \\ z^t & \tau^{\prime \prime}
\end{matrix} \right) \hbox{\rm with $\tau^{\prime} \in \H1$, $z \in \CC^{g-1}$
and $\tau^{\prime \prime} \in {\Hcal}_{g-1}$}. \eqno(5)
$$
From the definition of modular form it is clear that $f$ is invariant
under $\tau^{\prime} \mapsto \tau^{\prime}+b$ for $b \in \ZZ$ (given by 
an element of ${\rm Sp}(2g,\ZZ)$), hence we
have a Fourier series
$$
f= \sum_{m=0}^{\infty} \phi_m(\tau^{\prime \prime},z) e^{2 \pi i m \tau^{\prime}}.
$$
Here the function $\phi_m$ is a holomorphic function on ${\Hcal}_{g-1} \times
\CC^{g-1}$ satisfying certain transformation rules. More generally, if
we split $\tau$ as in (5) but with $\tau^{\prime} \in {\Hcal}_r$, $z \in 
\CC^{r (g-r)}$ and $\tau^{\prime \prime} \in {\Hcal}_{g-r}$ we find
a development
$$
\sum_m \phi_m(\tau^{\prime \prime},z)e^{2\pi i {\rm Tr}(m\tau^{\prime})},
$$
where the sum is over positive semi-definite half-integral matrices
$r \times r$ matrices $m$  and the functions
$\phi_m$ are holomorphic on ${\Hcal}_r \times \CC^{r(g-r)}$. For
$r=g$ we get back the Fourier expansion  and for $r=1$ we get
what is called the Fourier-Jacobi \index{Fourier-Jacobi development}
development.

For ease of explanation and to simplify matters we start with $g=2$. 
Then the function $\phi_m(\tau^{\prime},z)$ turns out to be a Jacobi form
\index{Jacobi form} 
of weight $k$ and index $m$, i.e., $\phi_m \in J_{k,m}$ which amounts to 
saying that it satisfies
\begin{enumerate}
\item $\phi_m((a\tau^{\prime}+b)/(\tau^{\prime}+d),z/(c\tau^{\prime}+d))=
(c\tau^{\prime}+d)^k e^{2\pi i m c z^2/(c\tau^{\prime}+d)} \phi_m(\tau^{\prime},z)$,
\item $\phi_m(\tau^{\prime},z+\lambda \tau^{\prime}+\mu)= e^{-2 \pi i m (\lambda^2 \tau^{\prime}+2\lambda z)} \phi_m(\tau^{\prime},z)$,
\item $\phi_m$ has a Fourier expansion of the form
$$
\phi_m = \sum_{n=0}^{\infty} \sum_{r \in \ZZ,\, r^2 \leq 4mn}
c(n,r) e^{2 \pi (n\tau^{\prime}+rz)}.
$$
\end{enumerate}
This gives a relation between Siegel modular forms for genus $2$
and Jacobi forms (see \cite{E-Z}) that we shall exploit later.
In the general case, if we split $\tau$ as
$$
\tau =\left( \begin{matrix} \tau^{\prime} & z \\ z^t & \tau^{\prime \prime}
\end{matrix} \right) \hbox{\rm with $\tau^{\prime} \in {\Hcal}_r$, 
$z \in \CC^{r(g-r)}$
and $\tau^{\prime \prime} \in {\Hcal}_{g-r}$}
$$
and a symmetric matrix $n$ as 
$(\begin{matrix} n' & \nu \\ \nu^t & n^{\prime \prime} \end{matrix} )$
and if we use the fact that ${\rm Tr}(n\tau)={\rm Tr}(n' \tau')+
2 {\rm Tr}(\nu z) + {\rm Tr}(n^{\prime \prime}\tau^{\prime \prime})$
then we can decompose the Fourier series of $f \in M_{\rho}(\Gamma_g)$
as 
$$
\sum_{n^{\prime \prime} \geq 0} \phi_{n^{\prime \prime}}(\tau', z)
e^{2 \pi i {\rm Tr}(n^{\prime \prime}\tau^{\prime \prime})}
$$
with $V$-valued holomorphic functions $\phi_{n^{\prime \prime}}(\tau', z)$
that satisfy the rules
\begin{enumerate}
\item For $\lambda, \mu \in \ZZ^g$ we have
$$
\phi_{n^{\prime \prime}}(\tau', z+\tau' \lambda +\mu)=
\rho(( \begin{matrix}
1_r & -\lambda \\ 0 & 1_{g-r}\\ \end{matrix} 
)) e^{-2\pi i {\rm Tr}(2 \lambda^t z + \lambda^t \tau' \lambda)} \phi_{
n^{\prime \prime}}(\tau', z).
$$
\item For $\gamma'=(a',b;c',d') \in \Gamma_{g-1}$ we have
$$
\phi_{n^{\prime \prime}}(\gamma'(\tau')),(c'\tau'+d')^{-t}z)=
e^{2\pi i {\rm Tr}(n^{\prime \prime} z^t (c' \tau'+d')^{-1} c' z)}
\rho(\left( 
\begin{matrix} c'\tau'+d' & c'z \\ 0 & 1_{g-r}\\ \end{matrix} 
\right))
\phi_{n^{\prime \prime}}(\tau',z).
$$
\item
$\phi_{n^{\prime \prime}}(\tau',z)$ is regular at infinity.
\end{enumerate}
The last condition means  that
$\phi_{n^{\prime \prime}}(\tau',z)$ has a Fourier expansion
$\phi_{n^{\prime \prime}}(\tau',z)=\sum c(m,r) \exp (2 \pi i {\rm Tr}
(m\tau'+ 2r^t z))$ for which $c(m,r)\neq 0$ implies that 
$(\begin{matrix} m & r \\ r^t & n^{\prime \prime} \\ \end{matrix} ) $ 
is positive semi-definite.
A holomorphic $V$-valued function $\phi(\tau',z)$ satisfying 1),2) and 3)
is called a Jacobi form of weight $(\rho',n^{\prime \prime})$.
The sceptical reader may frown upon this unattractive set of
transformation formulas, but there is a
natural geometric explanation 
for this transformation behavior that we shall see in section 
\ref{compactifications}. 

\end{section}
\begin{section}{The Ring of Classical Siegel Modular Forms for Genus Two}
\index{ring of classical modular forms}
So far we have not met any striking examples of Siegel modular forms. To
convince the reader that the subject is worthy of his attention we turn
to the first non-trivial case: classical Siegel modular forms of genus $2$.

For $g=1$ we know the structure of the graded ring $M_*(\Gamma_1)=\oplus_k M_k(
\Gamma_1)$. It is a 
polynomial ring generated by the Eisenstein series $e_4=E_4^{(1)}$ and $e_6
=E^{(1)}_6$ and
the ideal of cusp forms is generated by the famous cusp form $\Delta=
(e_4^3-e_6^2)/1728$ of weight $12$.

In comparison to this our knowledge of the graded ring 
$\oplus_{\rho \in \rm Irr} M_{\rho}$
of Siegel modular forms for $g=2$ 
is rather restricted and most of what we know concerns classical Siegel 
modular forms. A first
basic result was the determination by Igusa \cite{Ig1} of the ring of 
classical Siegel modular forms for $g=2$. We now know also the
structure of the ring of 
classical Siegel modular forms for $g=3$, a result of Tsuyumine,
\cite{Tsu}.

Recall that we have the Eisenstein Series $E_k^{(g)} \in M_k(\Gamma_g)$
for $k>g+1$. In particular, for $g=2$ we have $E_4=E_4^{(2)} \in M_4(\Gamma_2)$
and $E_6=E_6^{(2)} \in M_6(\Gamma_2)$. Let us normalize them here so that
$$
E_k=\sum_{(c,d)} \det(c\tau+d)^{-k},
$$
where the sum is over non-associated pairs of {\sl co-prime} 
symmetric integral matrices (non-associated w.r.t.\ to the multiplication on
the left by ${\rm GL}(g,\ZZ)$). The Fourier expansion of these
modular forms is known. If we write $\tau=\left( \begin{matrix}
\tau_1 & z \cr z & \tau_2\cr \end{matrix} \right)$
then
$$
E_k= \sum_{N} a(N) e^{2 \pi i {\rm Tr}(N\tau)},
$$
with constant term $1$ and 
for non-zero 
$N=(\begin{matrix} n & r/2 \\ r/2 & m\\ \end{matrix} )$ 
the coefficient $a(N)$ given as
$$
a(N)=\sum_{d|(n,r,m)} d^{k-1} H(k-1, \frac{4mn-r^2}{d^2})
$$
with $H(k-1,D)$ Cohen's function, i.e.,
$H(k-1,D)=L_{-D}(2-k)$, where $L_D(s)=L(s, \left( \frac{D}{} \right))$
is the Dirichlet $L$-series associated to $D$
if $D$ is $1$ or a discriminant of a real quadratic field, 
cf., \cite{E-Z}, p.\ 21. (this is essentially a class number.)
Explicitly we have with $q_j=e^{2\pi i \tau_j}$ and $\zeta=e^{2\pi i z}$
the developments (cf., \cite{E-Z})
\begin{align} \notag
E_4=1 +  240(q_1+q_2) +& 2160(q_1^2+q_2^2)+ \\ \notag
(240\, & \zeta^{-2}+13440\, \zeta^{-1}+30240+
13440 \, \zeta + 240 \, \zeta^2)q_1q_2
+ \ldots \\ \notag
\end{align}
and
\begin{align} \notag
E_6=1-504(q_1+q_2)&
-16632(q_1^2+q_2^2)+ \\ \notag 
+(-504\, \zeta^{-2}& +44352 \, \zeta^{-1}+166320+ 
44352 \, \zeta  -504\, \zeta^2)q_1q_2+ \ldots . \\ \notag
\end{align}
Under Siegel's operator $\Phi \colon
M_k(\Gamma_2) \to M_k(\Gamma_1)$ the Eisenstein series $E_k$
on $\Gamma_2$ maps to the Eisenstein series $e_k$ 
on $\Gamma_1$ for $k\geq 4$. In particular,
the modular form $E_{10} -E_4E_6$ maps to 
$e_{10} -e_4e_6$, and this is zero since $\dim M_{10}(
\Gamma_{1})=1$ and the $e_k$ are normalized so that their Fourier expansions
have constant term $1$. We thus find a cusp form. Similarly, $E_{12}
-E_6^2$ defines a cusp form of weight $12$ on $\Gamma_2$. 
To see that these are not zero
we restrict to the `diagonal' locus as follows.

Consider the map $\delta \colon 
\H1 \times \H1 \to \H2$ given by $(\tau_1,\tau_2) 
\mapsto \left( \begin{matrix} \tau_1 & 0 \cr 0 & \tau_2 \cr 
\end{matrix} \right)$. 
There is a corresponding map ${\rm SL}(2,\ZZ)\times {\rm SL}(2,\ZZ)
\to {\rm Sp}(4,\ZZ)$ by sending $(\begin{matrix} a & b\cr
c & d \cr \end{matrix}), (\begin{matrix} a' & b' \cr c' & d' \end{matrix})$
to $(A,B;C,D)$ (difficult to avoid capital letters here)
with $A=\left( \begin{matrix}
a &0 \cr 0 & a' \cr\end{matrix}\right)$, etc.\ that induces $\delta$
(on $({\rm SL}(2,\RR)/U(1))^2\to {\rm Sp}(4,\RR)/U(2)$).
If we use the coordinates 
$$
\tau=(\begin{matrix}\tau_1 & z \cr z & \tau_2\cr
\end{matrix}) \in {\Hcal}_2
$$
then the image of the map $\delta$ is given by $z=0$ and it is 
the fixed point locus of the involution on $\H2$ given by 
$(\tau_1,z,\tau_2) \mapsto (\tau_1,-z,\tau_2)$ induced by the
element $(A,B;C,D)$ with $A=(1,0;0,-1)=D$ and $B=C=0$.

An element $F \in M_k(\Gamma_2)$ 
can be developed around this locus $z=0$
$$
F= f(\tau_1,\tau_2)z^n + O(z^{n+1}) \qquad 
\hbox{\rm for some $n \in \ZZ_{\geq 0}$}. \eqno(1)
$$
It is now easy to check that
\begin{enumerate}
\item$f(\tau_1,\tau_2)\in M_{k+n}(\Gamma_1) \otimes M_{k+n}(\Gamma_1)$;
\item $f(\tau_2,\tau_1)= (-1)^k f(\tau_1,\tau_2)$;
\item$f(\tau_1,-z,\tau_2)=(-1)^kf(\tau_1,z,\tau_2)$.
\end{enumerate}
the first by looking at the action of ${\rm SL}(2,\ZZ)\times {\rm SL}(2,\ZZ)$
and the second by applying the involution 
$(A,B;C,D)$ with $A=D=(0,1;1,0)$ and
$B=C=0$ which interchanges $\tau_1$ and $\tau_2$ and the last by using the
involution $z \mapsto -z$.
The idea of developing along the diagonal locus was first used by
Witt, \cite{Wi}.

Developing $E_{10} -E_4E_6$ along $z=0$ and
writing $q_j=e^{2\pi i \tau_j}$ one finds
$c q_1q_2 z^2 +O(z^3)$, with $c \neq 0$,
so we normalize to get a cusp form $\chi_{10}=
E_{10}^{(1)}(\tau_1)\otimes E_{10}^{(1)}(\tau_2) z^2+ O(z^3)$. 
Similarly, the form $E_{12}^{(2)}
-(E_6^{(2)})^2$ gives after normalization a non-zero cusp form $\chi_{12}
= \Delta(\tau_1)\otimes \Delta(\tau_2) z^2+ O(z^3)$.

As we saw above in Section  \ref{thetaseries}
we also know the existence of a  cusp form $\chi_{35}$ of
odd weight $35$. 

   We now describe the structure of the ring of classical Siegel modular forms
for $g=2$. The theorem is due to Igusa and various proofs have been
recorded in the literature, cf.\ \cite{Ig1,Ig2,Ig3,Fr4,Aoki,Ham}. 
Here is another variant.

\begin{theorem}{\rm (Igusa)} The graded ring $M= \oplus_k  M_k(\Gamma_2)$
of classical Siegel modular forms of genus $2$ is generated by
$E_4,\, E_6, \, \chi_{10}, \, \chi_{12}$ and $\chi_{35}$ and
$$
M \cong \CC[ E_4,E_6,\chi_{10},\chi_{12},\chi_{35}]/(\chi_{35}^{2}=
R),
$$
where $R$ is an explicit (isobaric) polynomial 
in $E_4, E_6, \chi_{10}$ and $\chi_{12}$ (given on \cite{Ig2}, p.\ 849).
\begin{proof} (Isobaric means that every monomial has weight $70$
if $E_4$, $E_6$, $\chi_{10}$ and $\chi_{12}$ are given weights $4$, $6$,
$10$ and $12$.)
We start by introducing the vector spaces of modular forms:
$$
M_k^{\geq n}(\Gamma_1) = \{ f \in M_k(\Gamma_1) : f=O(q^n) 
\hbox{\rm \, at $\infty$} \} = \Delta^n \, M_{k-12n}(\Gamma_1)
$$
and
$$
M_k^{\geq n}(\Gamma_2)=\{ F \in M_k(\Gamma_2) : F=O(z^n)
\hbox{\rm \, near $\delta(\H1 \times \H1)$ } \} 
$$
We distinguish two cases depending on the parity of $k$.

{\bf $k$ even}. As we saw above  (use properties (1),(2), (3))
any element $F \in M_{k}^{\geq 2n}(\Gamma_2)$ can be
written as $F(\tau_1,z,\tau_2)=f(\tau_1,\tau_2)z^{2n} + O(z^{2n+2})$
with $f\in M_{k+2n}(\Gamma_1) \otimes M_{k+2n}(\Gamma_1)$ symmetric
(i.e.\ $f(\tau_1,\tau_2)=f(\tau_2,\tau_1)$) and $f=O(q_1^n,q_2^n)$.
This last fact follows from the observation that each Fourier-Jacobi 
coefficient $\phi_m(\tau_1,z)$ of $F$ is also $O(z^{2n})$, so is zero if
$2n > 2m$. We find an exact sequence
$$
0\to M_k^{\geq 2n+2}(\Gamma_2) \to M_k^{\geq 2n}(\Gamma_2) {\buildrel
r \over \longrightarrow} {\rm Sym}^2(M_{k+2n}^{\geq n}(\Gamma_1)) \to 0,
$$
where the surjectivity of $r$ is a consequence of the fact that
$$
{\rm Sym}^2(M_{k+2n}^{\geq n}(\Gamma_1)) = 
\CC[e_4\otimes e_4, e_6\otimes e_6, \Delta \otimes \Delta]
$$
and $\chi_{10}=\Delta(\tau_1) \Delta(\tau_2) z^2 + O(z^4)$
so that a modular form $\chi_{10}^n P(E_4,E_6,\chi_{12})$
with $P$ an isobaric polynomial maps to 
$P(e_4 \otimes e_4, e_6 \otimes e_6, \Delta \otimes \Delta)$. It follows that 
$$
\dim M_k(\Gamma_2)= \sum_{n=0}^{\infty} \dim {\rm Sym}^2(M_{k+2n}^{\geq n}
(\Gamma_1)) 
= \sum_{0 \leq n \leq k/10} \dim {\rm Sym}^2(M_{k-10n}(\Gamma_1)), 
$$
i.e., we get
\begin{align}
\sum_{k \, \text{\rm even}} \dim M_k(\Gamma_2) t^k&= \frac{1}{1-t^{10}}
\sum_{k \geq 0} \dim {\rm Sym }^2(M_k(\Gamma_1))t^k \cr
&= \frac{1}{1-t^{10}} \text{Hilbert series of 
$\CC[e_4\otimes e_4, e_6 \otimes e_6, \Delta \otimes \Delta]$} \cr
&= \frac{1}{(1-t^4)(1-t^6)(1-t^{10})(1-t^{12})}.\cr \notag
\end{align}

{\bf $k$ odd}. For $F \in M_k^{\geq 2n+1}(\Gamma_2)$ we find 
$f=O(q_1^{n+2},q_2^{n+2} )$. Since our Fourier-Jacobi coefficients 
$\phi_m(\tau_1,z)$ have 
a zero of order $2n+1$ 
 at $z=0$ and another three at the $2$-torsion points we see 
$2m\geq (2n+1)+3$ for non-zero $\phi_m$. 
Also we know that $f$ is anti-symmetric now, so 
$\dim M_k(\Gamma_2)  \leq
\sum_{n \geq 0} \dim \wedge^2(M_{k+2n+1}^{\geq n+2}(\Gamma_1))$
and this shows that for odd $k<35$ $\dim M_k(\Gamma_2)=0$. Since we
have a non-trivial form of weight $35$ we see that
$$
\sum_{k \, \rm odd} \dim M_k(\Gamma_2) t^k=
\frac{t^{35}}{(1-t^4)(1-t^6)(1-t^{10})(1-t^{12})}.
$$
The square $\chi_{35}^2$ is a modular form of even weight, hence can be
expressed as an polynomial $R$ in $E_4,E_6, \chi_{10}$ and $\chi_{12}$.
This was done by Igusa in \cite{Ig2}.
This completes the proof.
\end{proof}
\end{theorem}
\end{section}
\begin{section}{Moduli of Principally Polarized Complex Abelian Varieties}
For $g=1$ the quotient space $\Gamma_1\backslash {\Hcal}_1$ has an
interpretation as the moduli space of elliptic curves over the complex
numbers (complex tori of dimension $1$). To a point $\tau \in \H1$
we associate the complex torus $\CC / \ZZ +\ZZ \tau$. Then to a
point $(a\tau+b)/(c\tau+d)$ in the $\Gamma_1$-orbit of $\tau$
we associate the torus $\CC/\ZZ + \ZZ(a\tau+b)/(c\tau+d)$, and the
homothety $z \mapsto (c\tau+d) z$ defines an isomorphism of this torus
with $\CC/ \ZZ(c\tau+d) + \ZZ (a\tau+b) = \CC/\ZZ + \ZZ \tau$ since
$(c\tau+d, a \tau+b)$ is a basis of $\ZZ + \ZZ \tau$ as well. Conversely,
every $1$-dimensional complex torus can be represented as $\CC/
\ZZ + \ZZ \tau$. This can be generalized to $g>1$ as follows. A point
$\tau \in \Hg$ determines a complex torus $\CC^g/\ZZ^g + \ZZ^g \tau$,
but we do not get all complex $g$-dimensional tori. The following lemma,
usually ascribed to Lefschetz,
tells us what conditions this imposes. 

\begin{lemma} The following conditions on a complex torus $X=V/\Lambda$
are equivalent:
\begin{enumerate}
\item $X$ admits an embedding into a complex projective space;
\item $X$ is the complex manifold associated to an algebraic variety;
\item There is a positive definite Hermitian form $H$ on $V$
such that ${\rm Im}(H)$ takes integral values on $\Lambda \times \Lambda$.
\end{enumerate}
\end{lemma}
A complex torus satisfying these requirements is called a complex 
\Definition{abelian variety}. \index{abelian variety} 
For $g=1$ we could take $H(z,w)=z\bar{w}/{\rm Im}(\tau)$
on $\Lambda=\ZZ+\ZZ\tau$ and indeed, the map $\CC/\Lambda \to \PP^2$
given by $z \mapsto (\wp(z):\wp'(z):1)$ for $z\notin \Lambda$ with $\wp$ the
Weierstrass $\wp$-function defines the embedding. For $g>1$
we can take $H(z,w)=z^t ({\rm Im}(\tau))^{-1} \overline{w}$.
An $H$ as in the lemma is called a \Definition{polarization}. It is
called a \Definition{principal polarization} \index{polarization} if the map
${\rm Im}(H): \Lambda \times \Lambda \to \ZZ$ is unimodular. We shall
write $E={\rm Im}(H)$ for the alternating form that is the imaginary part
of $H$. Given a complex torus $X=V/\Lambda$ and a principal polarization 
on $\Lambda$ we can normalize things as follows. We choose an isomorphism
$V \cong \CC^g$ and choose a symplectic basis $e_1,\ldots,e_{2g}$ of the
lattice $\Lambda$ such that $E$ takes the standard form
$$
J = \left( \begin{matrix} 0_g & 1_g \\ -1_g & 0_g \\ \end{matrix} \right)
$$
with respect to this basis. These two bases yield us a period matrix 
\index{period matrix} $\Omega \in {\rm Mat}(g\times 2g,\CC)$ expressing 
the $e_i$ in terms of the chosen $\CC$-basis of $V$. A natural question
is which period matrices occur.  For this we note that $E$ is the imaginary
part of a Hermitian form $H(x,y)=E(ix,y)+\sqrt{-1} E(x,y)$
if and only if $E$ satisfies the condition $E(iz,iw)=E(z,w)$
for all $z,w \in V$ and this translates into (Exercise!)
$$
\Omega \, J^{-1} \Omega^t =0
$$
while the positive definiteness of $H$ translates into the condition
$$
2i (\bar{\Omega} J^{-1} \Omega^t)^{-1} \hbox{\rm is positive definite}.
$$
These conditions were found by Riemann \index{Riemann}
in his brilliant 1857 paper
\cite{Ri}.

If we now associate to $\Omega=(\Omega_1 \, \Omega_2)$ with $\Omega_i$
complex $g\times g$ matrices we see that the two conditions just found say that
if we put $\tau=\Omega_2^{-1} \Omega_1$ we have
$
\tau=\tau^t$, ${\rm Im}(\tau)>0$ i.e., $\tau$ lies in $\Hg$.
A change of basis of $\Lambda$ changes $(\tau \, 1_g)$ into
$(\tau a + c, \tau b +d)$ with $(a, b; c,d) \in {\rm Sp}(2g,\ZZ)$, 
but the corresponding torus is isomorphic to
$\CC^g/\ZZ^g(\tau b+d)^{-1}(\tau a + c)+ \ZZ^g$. In this way we see that
the isomorphism classes of complex tori with a principal polarization are
in 1-1 correspondence with the points of the orbit space 
$\Hg / {\rm Sp}(2g,\ZZ)$. If we transpose we can identify this orbit space
with the orbit space ${\rm Sp}(2g,\ZZ) \backslash \Hg$ for the usual
action $\tau \mapsto (a\tau +b)(c \tau +d)^{-1}$. 
\begin{proposition} There is a canonical bijection between the set of 
isomorphism classes of principally polarized abelian varieties of dimension 
$g$ and the orbit space $\Gamma_g \backslash \Hg$.
\end{proposition}
If we try to construct the whole family of abelian varieties we encounter
a difficulty. The action of the semi-direct product $\Gamma_g \ltimes \ZZ^{2g}$
on $\Hg \times \CC^g$ given by the usual action of $\Gamma_g$ on $\Hg$
and the action of $(\lambda, \mu) \in \ZZ^{2g}$ on a fibre $\{ \tau\} \times
\CC^{2g}$ by $z \mapsto z + \tau \lambda +\mu$ forces $-1_{2g}
\in \Gamma_g$ to act by
$-1$ on a fibre, so instead of finding the
complex torus $\CC^{g}/\ZZ^g+\tau \ZZ^g$ we get its quotient by the action
$z \mapsto -z$. However, if we replace $\Gamma_g$ by the congruence subgroup
$\Gamma_g(n)$ with $n\geq 3$ (see \cite{Serre}) then we get an honest family 
${\Xcal}_g(n)=\Gamma_g(n)\ltimes \ZZ^{2g} \backslash \Hg \times \CC^g
$ of abelian varieties
over $\Gamma_g(n)\backslash \Hg$. If we insist on using $\Gamma_g$ then
we have to work with orbifolds or stacks to have a universal family available;
the orbifold in question is the quotient of ${\Xcal}_g(n)$ under
the action of the finite group ${\rm Sp}(2g,\ZZ /n \ZZ)$.

The cotangent bundle of the family of abelian varieties over ${\Acal}_g(n)=
\Gamma_g(n) \backslash \Hg$ along the zero section defines a vector bundle
of rank $g$ on ${\Acal}_g(n)$. It can be constructed explicitly as a quotient
$\Gamma_g(n) \backslash \Hg \times \CC^g$ under the action of $\gamma
\in \Gamma_g(n)$ by $(\tau,z) \mapsto (\gamma(\tau),(c\tau +d)^{-t} z)$.
The bundle is called the \Definition{Hodge bundle} \index{Hodge bundle}
and denoted by $\EE=\EE_{g}$. The finite group
${\rm Sp}(2g,\ZZ / n \ZZ)$ acts on the bundle $\EE$ on ${\Acal}_g(n)$.
A section of $\det(\EE)^{\otimes k}$ 
that is ${\rm Sp}(2g,\ZZ/n\ZZ)$-invariant comes from a
holomorphic function on $\Hg$ that is a classical Siegel modular form of
weight $k$. Classical modular forms thus get a geometric interpretation.
In particular,
the determinant of the cotangent bundle of ${\Acal}_g(n)$, i.e., the
canonical bundle, is isomorphic to $\det(\EE)^{\otimes g+1}$; so to a
modular form $f$ of weight $g+1$ we can associate a top differential
form on $\Hg$ that is $\Gamma_g$-invariant via $f \mapsto
f(\tau) \prod_{i\leq j} d\tau_{ij}$. In a similar way 
one can construct
for each $\rho$ a vectorbundle over ${\Acal}_g(n)$ whose 
${\rm Sp}(2g,\ZZ /n \ZZ)$-invariant sections are the Siegel modular forms
of weight $\rho$ by applying the quotient $\Hg \times V$ by
$\Gamma_g$ under $(\tau,z) \mapsto (\gamma(\tau),\rho(c\tau +d) z)$,
see Section \ref{vectorbundles}.

The Hermitian form $H$ on the lattice $\Lambda \subset \CC^g$ can be 
viewed as the first Chern class 
(in $H^2(X,\ZZ)\cong \wedge^2 (H_1(X,\ZZ)^{\vee})\cong 
(\wedge^2 \Lambda)^{\vee}$)
of a line bundle $L$ on $X=\CC^g/\Lambda$
with $\dim_{\CC}H^0(X,L)=1$. A non-zero section determines an effective
divisor $\Theta$ on $X$. The line bundle $L$ and the corresponding
 divisor $\Theta$ are determined by $H$ up to translation over $X$.
If we require that $\Theta$ be invariant under $z \mapsto -z$ then
$\Theta$ is unique up to translation over a point of order $2$ on $X$
and then $2\Theta$ is unique.

If we pull a non-zero section $s$ of $L$ back to the universal cover $\CC^g$
then we obtain a holomorphic function with a certain transformation
behavior under translations by elements of $\Lambda$. An example of such
a function is provided by Riemann's theta function 
\index{Riemann's theta function}
$$
\theta(\tau,z)=\sum_{n \in \ZZ^g} e^{\pi i (n^t \tau n +2 n^tz)},
\qquad (\tau \in \Hg, z \in \CC^g)
$$
a series that converges very rapidly and defines a holomorphic function
that satisfies for all $\lambda, \mu \in \ZZ^g$
$$
\theta(\tau,z+\tau\lambda+\mu)=e^{-\pi i(\lambda^t \tau \lambda +2 \lambda^t z)}
\theta(\tau,z).
$$
Conversely, a holomorphic function $f$ on $\CC^g$ that satisfies
for all $\lambda, \mu \in \ZZ^g$
$$
f(z+\tau\lambda+\mu)=e^{-\pi i(\lambda^t \tau \lambda +2 \lambda^t z)}
f(z)
$$
is up to a multiplicative constant precisely 
$\theta(\tau,z)$ as one sees by
developing $f$ in a Fourier series $f=\sum_{n} c(n)\exp{ 2\pi i n^tz}$
and observing that addition of a column $\tau_k$ of $\tau$ to $z$ produces
$$
f(z+\tau_k)=\sum_n  c(n) \exp{(2\pi i n^t (z+\tau_k))}=
\sum_n c(n) \exp{(2\pi i n^t \tau_k)} \exp{(2\pi i n^t z)}
$$
from which one obtains $c(n+e_k)=
c(n) \exp{(2\pi i n^t \tau_k+ \pi i \tau_{kk})}$
and gets that $f$ is completely determined by $c(0)$.

If $S$ is a compact Riemann surface of genus $g$ it determines a
Jacobian variety ${\rm Jac}(S)$ which is a principally polarized 
complex abelian variety of dimension $g$. Sending $S$ to
${\rm Jac}(S)$ provides us with a
map ${\Mcal}_g(\CC) \to \Gamma_g \backslash \Hg$ from the moduli
space of compact Riemann surfaces of genus $g$ to the moduli
of complex principally polarized abelian varieties of dimension $g$
which is injective by
a theorem of Torelli\index{Torelli's theorem}. 
The geometric interpretation given for Siegel
modular forms thus pulls back to the moduli of compact Riemann surfaces.

\end{section}

\begin{section}{Compactifications}\label{compactifications}
\index{compactifications}
It is well known that $\Gamma_1 \backslash \H1$ is not compact, 
but can be compactified by adding the cusp, that is, the orbit
of $\Gamma_1$ acting on $\QQ \subset \bar{\Hcal}_1$.
Or if we use the equivalence of $\H1$ with the unit disc $D_1$
given by $\tau \mapsto (\tau -i)/(\tau +i)$ then we add to $D_1$ the
rational points of the boundary of the unit disc and take the orbit
space of this enlarged space. We can do something similar for $g>1$
by considering the bounded symmetric domain
$$
D_g =\{ z \in {\rm Mat}(g\times g, \CC) : z^t=z, z^t \cdot \bar{z} < 1_g \}
$$
which is analytically equivalent to $\Hg$. We now enlarge this space by
adding not the whole boundary but only part of it as follows. Let
$$
D_r =\{ \left( \begin{matrix} z' & 0 \\ 0 & 1_{g-r} \\ \end{matrix} \right)
: z' \in D_r \} \subset \bar{D}_g
$$
and define now $D_g^*$ to be the union of all $\Gamma_g$-orbits of these
$D_r$ for $0 \leq r \leq g $. 
Note that $\Gamma_g$ acts on $D_g$ and on its closure $\bar{D}_g$.
Then $\Gamma_g$ acts in a natural way on $D_g^*$
and the orbit space decomposes naturally as a disjoint union
$$
\Gamma_g \backslash D_g^* = \sqcup_{i=0}^g \Gamma_i \backslash D_i.
$$
Going back to the upper half plane model this means that we consider
$$
\Gamma_g \backslash {\Hcal}_g^* = \sqcup_{i=0}^g \Gamma_i \backslash {\Hcal}_i
$$
Satake has shown how to make this space into a normal analytic space,
the Satake compactification
\index{Satake compactification}. One first defines a topology on $\Hg^*$
and then a sheaf of holomorphic functions. The quotient 
$\Gamma_g \backslash \Hg^*$ then becomes a normal analytic space.
By using explicitly constructed modular forms one then shows that classical
modular forms of a suitably high weight separate points and tangent vectors and
thus define an embedding of $\Gamma_g \backslash \Hg^*$ into projective
space. By Chow's lemma it is then a projective variety. The following
theorem is a special case of a general theorem due to Baily and Borel,
\cite{B-B}. \index{Baily-Borel compactification}

\begin{theorem} Scalar Siegel modular forms of an appropriately high weight 
define an embedding of $\Gamma_g \backslash {\Hg}^*$ into projective 
space and the image of $\Gamma_g \backslash \Hg$ (resp.\
$\Gamma_g \backslash \Hg^*$) is a quasi-projective (resp.\
a projective) variety.
\end{theorem}

The resulting Satake or Baily-Borel compactification is for $g>1$ very singular.
As a first attempt at constructing a smooth compactification we
reconsider the case $g=1$. In ${\Hcal}_{1,c} = \{ \tau \in \H1 \colon
{\rm Im}(\tau)\geq c \}$ with $c>1$ the action of $\Gamma_1$ reduces to
the action of $\ZZ$ by translations $\tau \mapsto \tau +b$. So consider
the map ${\Hcal}_{1,c} \to \CC^*$, $\tau \mapsto q=\exp{2\pi i \tau}$.
It is clear how to compactify ${\Hcal}_{1,c}/\ZZ$: just add the origin $q=0$
to the image in $\CC^* \subset \CC$. 
In other words, glue 
$\Gamma_1\backslash {\Hcal}_1$ with $\ZZ \backslash {\Hcal}_{1,c}^*$ 
over $\ZZ \backslash {\Hcal}_{1,c}$.
To do something similar for $g>1$
we consider the subset (for a suitable real symmetric $g\times g$-matrix
$c\gg 0$ which is sufficiently positive definite)
$$
{\Hcal}_{g,c}=\{ \tau =\left( \begin{matrix} \tau_1 & z \\ z^t & \tau_2 \\
\end{matrix} \right)\in \Hg : {\rm Im}(\tau_2) - {\rm Im}(z^t){\rm Im}(\tau_1)^{-1}
{\rm Im}(z) \geq c \}
$$
The action of $\Gamma_g$ in ${\Hcal}_{g,c}$ reduces to the action of 
the subgroup  $P$ 
$$
\left\{ \left(
\begin{matrix}
a & 0 & b & * \\
* & \pm 1 & * & * \\
c & 0 & d & * \\
0 & 0 & 0 & \pm 1 \\
\end{matrix}
\right) \in \Gamma_g, \, 
\left( \begin{matrix} a & b \\ c & d \\ \end{matrix} \right) 
\in \Gamma_{g-1} 
 \right\},
$$
the normalizer of the `boundary component' ${\Hcal}_{g-1}$.
We now make a map 
$$
{\Hcal}_{g,c} \to {\Hcal}_{g-1} \times \CC^{g-1} \times \CC^*,
\quad \tau \mapsto (\tau_1, z, q_2=\exp(2 \pi i \tau_2)).
$$
The associated parabolic subgroup $P$ acts on 
${\Hcal}_{g-1} \times \CC^{g-1} \times \CC^{*}$ and this action can be extended
to an action on ${\Hcal}_{g-1} \times \CC^{g-1} \times \CC$, where
$$
\left( \begin{matrix}
1_{g-1} & 0 & 0 & 0 \\
0 & 1 & 0 & b \\
0 & 0 & 1_{g-1} & 0 \\
0 & 0 & 0 & 1 \\
\end{matrix} \right)
$$
acts now by $(\tau_1,z,q_2) \mapsto (\tau_1,z, e^{2 \pi i b} q_2)$
while the matrix
$$
\left(\begin{matrix}
1_{g-1} & 0 & 0 & l \\
m & 1 & l & 0 \\
0 & 0 & 1_{g-1} & -m \\
0 & 0 & 0 & 1 \\
\end{matrix} \right)
$$ 
acts by 
$(\tau_1,z,q_2) \mapsto (\tau_1, z+ \tau_1 m + l, e^{2 \pi i (m\tau_1 m
+ 2 m z + lm)} q_2)$, and the diagonal matrix with entries
$(1,\ldots,-1,1,\ldots,1,-1)$ acts by $(\tau_1,z,\zeta) \mapsto
(\tau_1,-z,\zeta)$
and finally 
$(a,b;c,d) \in \Gamma_{g-1}$ acts on ${\Hcal}_{g-1} \times \CC^{g-1} 
\times \CC$ by
$$
(\tau_1,z,q_2) \mapsto  (\gamma(\tau_1), (a- (\gamma(\tau_1)c)\, z,
, \tau_2-z^t (c\tau_1+d)^{-1}c z)
$$
and this action can be extended similarly.

We now have an embedding $\Gamma_g \backslash {\Hcal}_{g,c} \longrightarrow
P\backslash {\Hcal}_{g-1} \times \CC^{g-1}\times \CC$ 
and by taking the closure of
the image we obtain a `partial compactification'.
The quotient of ${\Hcal}_{g-1}\times \CC^{g-1} 
\times \{ 0\}$ by this action is the `dual universal abelian variety' 
$\hat{\Xcal}_{g-1}=
\Gamma_{g-1}\ltimes \ZZ^{2g-2} \backslash {\Hcal}_{g-1} \times \CC^{g-1}$
over $\Gamma_{g-1}\backslash{\Hcal}_{g-1}$. Note that a principally polarized
abelian variety is isomorphic to its dual,
so we can enlarge our orbifold $\Gamma_g \backslash \Hg$
by adding this orbifold quotient ${\Xcal}_{g-1}=\Gamma_{g-1} \times \ZZ^{2g-2} 
\backslash {\Hcal}_{g-1} \times \ZZ^{2g-2}$. The result is a partial
compactification
$$
{\Acal}_g^{(1)}={\Acal}_g \sqcup {\Xcal}^{\prime}_{g-1},
$$
where the prime refers to the fact that we are dealing
with orbifolds and have to divide by (at least) an extra involution
since a semi-abelian variety generically has $\ZZ/2 \times \ZZ/2$ as
its automorphism group, while a generic abelian variety has only $\ZZ/2$.

This space parametrizes principally polarized complex
abelian varieties of dimension
$g$ or degenerations of such (so-called semi-abelian varieties
\index{semi-abelian variety}
of torus rank $1$)
that are extensions 
$$
1 \to \GG_m \to \tilde{X} \to X\to 0
$$
of a $g-1$-dimensional principally polarized complex abelian variety 
by a rank $1$ torus $\GG_m=\CC^*$. Such extension classes
are classified by the dual abelian variety $\hat{X} \cong X$ 
(associate to a line bundle on $X$ the $\GG_m$-bundle obtained by
deleting the zero section) which explains
why we find the universal abelian variety of dimension $g-1$ in the `boundary'
of ${\Acal}_g$. (There is the subtlety whether one allows isomorphisms to be
$-1$ on $\GG_m$ or not.)
This partial compactification \index{partial compactification}
is canonical. If we wish to
construct a full smooth compactification one can use Mumford's
theory of toroidal compactifications, but unfortunately there is (for $g\geq 4$)
no unique such compactification. We refer e.g.\ to \cite{Nam}.

This partial compactification
enables us to reinterpret the Fourier-Jacobi series 
\index{Fourier-Jacobi series}
of a Siegel modular form. 
In particular, the formulas in Section \ref{Fourier-Jacobi} 
tell us that the pull back of $f$
to a fibre of ${\Xcal}_{g-1} \to {\Acal}_{g-1}$
is an abelian function and that $f$ restricted to the zero-section
of ${\Xcal}_{g-1} \to {\Acal}_{g-1}$ is a Siegel modular form
of weight $k-1$ on $\Gamma_{g-1}$. 

We can be more precise. We work with a group $\Gamma_g(n)$ with
$n\geq 3$ or interpret everything in the orbifold sense.
The normal bundle of ${\Xcal}_{g-1}$ is the
line bundle $O(-2\Theta)$, as one can deduce from the action given above.

We can also extend the Hodge bundle $\EE=\EE_g$ to a vector bundle on
${\Acal}_g^{(1)}$. On the boundary divisor ${\mathcal X}_{g-1}^{\prime}$ it
is the extension of the pull back $\pi^*\EE_{g-1}$
from ${\Acal}_{g-1}$ to ${\Xcal}_{g-1}$ by a line bundle.

So if we are given a classical Siegel modular form of weight $k$
we can interpret it as a section of $\det(\EE)^{\otimes k}$ 
and develop (the pull back of) 
$f$ along the boundary ${\Xcal}_{g-1}$ where the $m$-th term in the
development is a section of 
$$
(\det (\EE)_{|{\Xcal}_{g-1}})^{\otimes k} \otimes O(-2m\Theta)
$$
on ${\Xcal}_{g-1}$. This gives us a geometric interpretation of the
Fourier-Jacobi development. \index{Fourier-Jacobi development}

Of course, it is useful to have not only a partial 
compactification, but a smooth compactification. 
The theory of toroidal compactifications
developed by Mumford and his co-workers Ash, Rapoport and Tai
provides such compactifications $\tilde{\mathcal A}_g$.
They depend on the choice of a certain cone decomposition
of the cone of positive definite bilinear forms in $g$
variables, cf.\ \cite{AMRT}. The `boundary' $\tilde{\mathcal A}_g
-{\mathcal A}_g$ is a divisor with normal crossings and one has
a universal semi-abelian variety over $\tilde{\mathcal A}_g$
in the orbifold sense.

\end{section}
\begin{section}{Intermezzo: Roots and Representations}\label{roots}
Here we record a few concepts and notations that we shall need
in the later sections. The reader may want to skip this on a first reading.

Recall that we started out in Section \ref{SiegelModularGroup}
with a symplectic lattice
$(\ZZ^{2g},\langle \, , \, \rangle)$ with a basis 
$e_1,\ldots,e_g,f_1,\ldots,f_g$ with $\langle e_i,f_j\rangle =\delta_{ij}$
and $\langle e_1,\ldots,e_g \rangle$ and $\langle f_1,\ldots,f_g\rangle$
isotropic subspaces. We let $G:={\rm GSp}(2g,\QQ)$ be the group
of rational symplectic similitudes (transformations that preserve the form up to a scalar), 
viz.,
$$
G:={\rm GSp}(2g,\QQ)=\{ \gamma \in {\rm GL}(\QQ^{2g}) :
\gamma^t J \gamma = \eta(\gamma) J\}
$$
and $G^{+}=\{ \gamma \in G \colon \eta(\gamma)>0\}$. 
Note that $\det(\gamma)=\eta(\gamma)^g$ for $\gamma \in G$
and that $G^0={\rm Sp}(2g,\ZZ)$ is the kernel of the map that sends 
$\gamma$ to $\eta(\gamma)$ on $G^{+}(\ZZ)$. For $\gamma \in G$
the element $\eta(\gamma)$ is called the \Definition{multiplier}. 
\index{multiplier}
Note that we view elements of $\ZZ^{2g}$ as column vectors and $G$
acts from the left.

There are several important subgroups that play a role in the sequel.
Given our choice of basis there is a natural Borel subgroup $B$
respecting the symplectic flag 
$\langle e_1 \rangle \subset\langle e_1,e_2\rangle \subset \ldots
\subset 
\langle e_1,e_2\rangle^{\bot} \subset \langle e_1 \rangle^{\bot}$. 
It consists of the matrices $(a,b;0,d)$ with $a$ upper triangular
and $d$ lower triangular.

Other natural subgroups are: 
the subgroup $M$ of elements respecting the
decomposition $\ZZ^g \oplus \ZZ^g$ of our symplectic space. It is
isomorphic to ${\rm GL}(g)\times \GG_m$ and consists of the matrices
$\gamma=(a,0;0,d)$ with $ad^t=\eta(\gamma) 1_g$. Furthermore, we
have the Siegel (maximal) parabolic subgroup $Q$ of elements that stabilize
the first summand $\ZZ^g =\langle e_1,\ldots,e_g\rangle$; it 
consists of the matrices $(a,b;0,d)$. It contains the subgroup $U$
(unipotent radical) of matrices of the form $(1_g,b;0,1_g)$ with $b$
symmetric that act as the identity of the first summand $\ZZ^g$.

Another important subgroup of $G$ is the diagonal torus $\TT$ 
isomorphic to $\GG_m^{g+1}$ of matrices 
$\gamma={\rm diag}(a_1,\ldots,a_g,d_1,\ldots,d_g)$ 
with $a_id_i=\eta(\gamma)$. 
Let $X$ be the character group of $\TT$;
it is generated by the characters $\epsilon_i: \gamma \mapsto 
a_i$ for $i=1,\ldots,g$ and $\epsilon_0(\gamma)=\eta(\gamma)$.
 Let $Y$ be the co-character group 
of $\TT_m$, i.e., $Y={\rm Hom}(\GG_m,\TT)$. 
This group is  isomorphic to the group $\ZZ^{g+1}$ of $g+1$-tuples
with $(\alpha_1,\ldots,\alpha_g,c)$ corresponding to
 the co-character
$t \mapsto {\rm diag}(t^{\alpha_1},\ldots,t^{\alpha_g},
t^{c-\alpha_1},\ldots,t^{c-\alpha_g})$.
We fix a basis of $Y$ by letting $\chi_i$ for
$i=1,\ldots,g$ correspond to $\alpha_j=\delta_{ij}$ and $c=0$
and $\chi_0$ to $\alpha_j=0$ and $c=1$. Then the
characters and co-characters pair via $\langle \epsilon_i ,
\chi_j\rangle =\delta_{ij}$.

The adjoint action of $\TT$ on the Lie algebras of $M$ and $G$ defines
root systems $\Phi_M$ and $\Phi_G$ in $X$. Concretely, we may take
as simple roots $\alpha_i=\epsilon_i-\epsilon_{i+1}$ for $i=1,\ldots,g-1$
and $\alpha_g=2\epsilon_g-\epsilon_0$ and coroots $\alpha_i^{\vee}=
\chi_i-\chi_{i+1}$ for $i=1,\ldots,g-1$ and $\alpha_g^{\vee}=\chi_g$.

The set $\Phi_G^{+}$ of positive roots (those occuring in the Lie algebra
of the nilpotent radical of $B$) consists of the 
so-called compact roots 
$\Phi_M^{+}=\{\epsilon_i-\epsilon_j \colon 1\leq i < j \leq g \}$
and the non-compact roots 
$\Phi^{+}_{\rm nc}=\{ \epsilon_i +\epsilon_j-\epsilon_0 
\colon 1 \leq i, j \leq g\}$. We let $2\varrho=2\varrho_G$ (resp.\
$2\varrho_M$) be the sum of the positive roots in $\Phi_G^{+}$
(res.\ $\Phi_M^{+}$). When viewed as characters  $2\varrho_M$
corresponds to $\gamma \mapsto \prod_{i=1}^g a_i^{g+1-2i}$
and $2\varrho_G$ to $\gamma \mapsto \eta(\gamma)^{-g(g+1)/2}
\prod_{i=1}^g a_i^{2g+2-2i}$.

There is a symmetry group acting on our situation, the Weyl group
\index{Weyl group}
$W_G=N(\TT)/\TT$, with $N(\TT)$ the normalizer of $\TT$ in $G$. This group
$W_G$ is isomorphic to $S_g \ltimes (\ZZ/2\ZZ)^g$, where the generator
of the $i$-th factor $\ZZ/2\ZZ$ acts on a matrix of the form
${\rm diag}(a_1,\ldots,a_g,d_1,\ldots,d_g)$ by
interchanging $a_i$ and $d_i$ and the symmetric group $S_g$
acts by permuting the $a$'s and $d$'s. The Weyl group
of $M$ (normalizer this time in $M$)
is isomorphic to the symmetric group $S_g$.
We have positive Weyl chambers 
$
P_G^{+}=\{ \chi \in Y \colon \langle \chi , \alpha \rangle \geq 0
\, \hbox{\rm for all $\alpha \in \Phi_G^{+}$ } \}
$
and similarly for $M$:
$
P_M^{+}=\{ \chi \in Y \colon \langle \chi , \alpha \rangle \geq 0
\, \hbox{\rm for all $\alpha \in \Phi_M^{+}$ } \}
$
giving the dominant weights.
\begin{lemma}
The irreducible complex representations of $G$ (resp.\ $M$) correspond
to integral weights in the chamber $P_G^{+}$ (res.\ $P_M^{+}$)
that come from characters of
$\TT$. 
\end{lemma}

Sometimes we just work with $G^0$ and $M^0=M \cap G^0$.
This means that we forget about the action of the multiplier
$\eta$. \index{multiplier}

We can give a set $W_0$ of $2^g$ canonical coset representatives of 
$W_M\backslash W_G$, the Kostant representatives\index{Kostant representatives}, which are 
characterized by the conditions
$$
W_0=\{ w \in W_G \colon \Phi_M^+ \subset w(\Phi_G^{+})\}
= \{ w \in W_G \colon w(\varrho)-\varrho \in P_M^{+}\}.
$$
With our normalizations we have $\varrho=(g,g-1,\ldots,1,0)$ and
$2\varrho_M=(g+1,\ldots,g+1,-g(g+1)/2)$. If we restrict to $G^0$ and $M^0$
then dominant weights for
$M^0\cong {\rm GL}(g)$ are given by $g$-tuples
$(\lambda_1,\ldots,\lambda_g)$ with $\lambda_i \geq \lambda_{i+1}$
for $i=1,\ldots,g-1$. A coset in $W_M \backslash W_G$ is given by
a vector $s$ (in $\{\pm 1\}^g$) of $g$ signs. 
The Kostant representative of $s$ is the  element $\sigma\, s$
such that 
$(s_{\sigma(1)}\lambda_{\sigma(1)},\ldots,
s_{\sigma(g)}\lambda_{\sigma(g)})$ is in $P^{+}_M$, i.e.,
$s_{\sigma(i)}\lambda_{\sigma(i)}\geq s_{\sigma(i+1)}
\lambda_{\sigma(i+1)}$ for $i=1,\ldots,g-1$
for all $(\lambda_1 \geq \ldots \geq \lambda_g)$. 
\end{section}
\begin{section}{Vector Bundles defined by Representations}\label{vectorbundles}

Let $\pi \colon {\mathcal X}_g \to {\mathcal A}_g$ be the universal
family of abelian varieties over ${\mathcal A}_g$.
The Hodge bundle $\EE =\pi_*\Omega_{{\mathcal X}_g/{\mathcal A}_g}$,
a holomorphic bundle of rank $g$, 
and the de Rham bundle $R^1\pi_* \CC$ on ${\mathcal A}_g$,
a locally constant sheaf of rank $2g$, are 
examples of vector bundles
associated to representations of ${\rm GL}(g)$ and ${\rm GSp}(2g)$. 
Their fibres at a point $[X] \in {\mathcal A}_g$ are 
$H^0(X,\Omega_X^1)$ and $H^1(X,\CC)$. 
The first is a holomorphic vector bundle, the second a local system.
Both are important for understanding Siegel modular forms.

To define these bundles
recall the description of ${\mathcal H}_g$ as an open part $Y_g^{+}$ 
of the symplectic Grassmann variety $Y_g$ given in Section 
\ref{SiegelModularGroup}. We can identify $Y_g$ with $G(\CC)/Q(\CC)$
with $Q$ the subgroup fixing the totally isotropic first summand $\CC^g$
of our complexified symplectic lattice $(\ZZ^g,\langle \, , \, \rangle)
\otimes \CC$.
If $\rho \colon Q^0 \to {\rm End}(V)$ is a complex representation
(with $Q^0=Q \cap G^0$)
then we can define a $G^0(\CC)$-equivariant vector bundle 
${\mathcal V}_{\rho}$
on $Y_g$ by ${\mathcal V}_{\rho}=G^0(\CC) \times^{Q^0(\CC)} V$ as the quotient
of $G^0(\CC) \times V$ under the equivalence relation 
$(g,v) \sim (g\, q, \rho(q)^{-1} v)$ for all $g \in G^0(\CC)$ and 
$q \in Q^0(\CC)$. Then $\Gamma_g$ (or any finite index subgroup $\Gamma'$) 
acts on ${\mathcal V}_{\rho}$ and the quotient is a vector bundle $V_{\rho}$
on ${\mathcal A}_g$ in the orbifold sense (or a true one if $\Gamma'$ acts
freely on ${\mathcal H}_g$). 

Recall that $M$ is the subgroup of ${\rm GSp}(2g,\QQ)$ respecting the 
decomposition $\QQ^g \oplus \QQ^g$ of our symplectic space and $M^0=
M \cap {\rm Sp}(2g,\QQ)$.
If we are given a complex representation of $M^0(\CC) \cong {\rm GL}(g)$
(or of $M \cong {\rm GL}(g)\times \GG_m$) we can obtain a vector bundle
by extending the representation to a representation on $Q^0(\CC)$
by letting it be trivial on the unipotent radical $U$ of $Q$.
(Note that $Q= M\cdot U$.) If we do this with the tautological 
representation of $M^0$ we get the Hodge bundle $\EE$. But there
is a subtle point here. If we work with $M$ instead of $M^0$ then
the Hodge bundle is given by the representation of $M$
that acts by $\eta(\gamma)^{-1} a $ on $\CC^g$ for $\gamma
=(a,0;0,d)$.

In any case we thus get a holomorphic vector bundle ${\mathcal W}(\lambda)$
associated to each dominant weight $(\lambda_1 \geq \ldots, \geq
\lambda_g)$ of ${\rm GL}(g)$. 
Another way of getting these vector bundles thus associated to the
irreducible representations of $M^0$ (or $M$) is by starting from
the Hodge bundle and applying Schur operators (idempotents)
to the symmetric
powers of $\EE$ analogously to the way one gets the 
corresponding representations from the standard one.
Since the Hodge bundle $\EE$ extends over a toroidal compactification
$\tilde{\mathcal A}_g$ this makes it clear that these vector bundles ${\mathcal W}(\lambda)$
can be extended over any toroidal compactification as constructed
by Mumford (or Faltings-Chai). The space of sections can be identified
with a space of modular forms $M_{\rho}$ and it thus
follows from general theorems in algebraic geometry that these
spaces of Siegel modular forms $M_{\rho}$ are finite dimensional.

Another important vector bundle is the bundle associated to the
first cohomology of the universal abelian variety ${\mathcal X}_g$
with fibre $H^1(X,\CC)$; more precisely, it is given by $\VV:= R^1\pi_*\CC$
with $\pi\colon {\mathcal X}_g \to {\mathcal A}_g$ the 
universal abelian variety.
It can be gotten from the construction just given by taking the
dual or contragredient of the
standard or tautological representation of ${\rm Sp}(2g,\CC)$
and restricting it to $Q^0(\CC)$. 
(Again, if one takes the multiplier into account--as one should--
then $R^1\pi_* \CC$ corresponds to $\eta^{-1}$ times the standard
representation.)
In this case we find a flat bundle:
all the bundles ${\mathcal V}_{\rho}$ on $Y_g$ come with a trivialization
given by $[(g,v)] \mapsto \rho(g)v$. So the quotient bundle 
carries a natural integrable connection. So $\VV$ is a local system
(locally constant sheaf). We thus find for each
dominant weight $\lambda=(\lambda_1 \geq \ldots \geq \lambda_g,c)$ 
of $G$ a local system ${\VV}_{\lambda}(c)$ on ${\mathcal A}_g$.
The multiplier representation defines a local system of rank $1$
denoted by $\CC(1)$ and we can twist $\VV_{\lambda}(c)$ by the $n$th power
of $\CC(1)$ to change $c$, cf.\ Section \ref{roots}.

\end{section}
\begin{section}{Holomorphic Differential Forms}
\index{holomorphic differential form}
Let $\Gamma' \subset \Gamma_g$ be a subgroup of finite index which
acts freely on ${\mathcal H}_g$, e.g., $\Gamma'=
\ker \{ {\rm Sp}(2g,\ZZ) \to {\rm Sp}(2g,\ZZ/n\ZZ\}$ for $n \geq 3$.
Let $\Omega^i$ be the sheaf of holomorphic $i$-forms on ${\mathcal H}_g$.
A section of $\Omega^1$ can be written as
$$
\omega ={\rm Tr}(f(\tau) d\tau),
$$
where $d\tau=(d\tau_{ij})$ and $f$ is a symmetric matrix of holomorphic
functions on ${\mathcal H}_g$. Then $\omega$ is invariant under the action of
$\Gamma'$ if and only if $f(\gamma(\tau))=(c\tau +d) f(\tau) (c\tau+d)^t$
for all $\gamma=(a,b;c,d) \in \Gamma'$. Note that if $r$ is the
standard representation of ${\rm GL}(g,\CC)$ on $V=\CC^g$ then the action
on symmetric bilinear forms ${\rm Sym}^2(V)$ is given by $b \mapsto
r(g)\, b \, r(g)^t$. So the space of holomorphic $1$-forms on
$\Gamma' \backslash {\mathcal H}_g$ can be identified with $M_{\rho}(\Gamma')$,
with $\rho$ the second symmetric power of the standard representation
and the space of holomorphic $i$-forms with $M_{\rho{\prime}}(\Gamma')$ 
with $\rho^{\prime}$
equal to the $i$th exterior power of ${\rm Sym}^2 V$. 
So we find an isomorphism $\Omega^1_{\Gamma^{\prime} \backslash 
{\mathcal H}_g} \cong {\rm Sym}^2{\EE}$ and this can be extended
over a toroidal compactification $\tilde{\mathcal A}$ to an isomorphism
$$
\Omega^1_{\tilde{\mathcal A}}(\log D) \cong {\rm Sym}^2(\EE)
$$
with $D$ the divisor at infinity.
(But again, one should be aware of the action of the multiplier:
if one looks at the action of ${\rm GSp}(2g,\RR)^+$
one has $d((a\tau +b)(c\tau+d)^{-1})= \eta(\gamma)(c\tau +d)^{-1} \,
d\tau \, (c\tau+d)^{-1}$.)

The question arises which
representations occur in $\wedge^i {\rm Sym}^2(V)$?

The answer is given in terms of roots. 
A theorem of Kostant \cite{Ko} tells that the irreducible
representations $\rho$ of ${\rm GL}(g,\CC)$ that occur in the exterior
algebra $\wedge^{\ast} {\rm Sym}^2(V)$ with $V$ the standard 
representation of ${\rm GL}(g,\CC)$ are those $\rho$ for which the
dual $\hat{\rho}$ is of the form $w\delta -\delta$ with $\delta=
(g,g-1,\ldots,1)$ the half-sum of the positive roots and $w$ in 
the set $W_0$ of Kostant representatives. Now if 
$\hat{\rho}=(\lambda_1\geq\lambda_1 \ldots \geq \lambda_g)$ 
occurs in this exterior algebra then $w\delta$ is of the form
$(g-\lambda_g,g-1-\lambda_{g-1},\ldots,1-\lambda_1)$. If $\alpha$
is the largest integer that occurs among the entries of $w\delta$
then either $\alpha=-1$ or $1 \leq \alpha \leq g$. In the latter case
$w\delta$ is of the form 
$(\alpha,*,\ldots,*,-\alpha -1,-\alpha -2,\ldots,-g)$
and it follows that $\lambda_{g-\alpha}=g+1$. This implies that
the number of $\lambda_j$ with $\lambda_j=\lambda_g$ (the co-rank
of $\hat{\rho}$, cf., Section \ref{siegeloperator}) 
plus the number of those with $\lambda_j=\lambda_g+1$
is at most $\alpha$. The vanishing theorem 
of Weissauer (Thm.\ \ref{Verschwindungssatz})
now implies that non-zero differentials
can only come from representations that are of the form
$$
\rho=(g+1,g+1,\ldots,g+1)
$$
which corresponds to top differentials ($\wedge^{g(g+1)/2} \Omega^1$)
and classical Siegel modular forms of weight $g+1$, or of the form
$$
\rho=(g+1,g+1,\ldots,g+1,g-\alpha,\ldots,g-\alpha),
$$
with $1 \leq \alpha \leq g$ and these occur in 
$\wedge^p \Omega^1$ with $p=g(g+1)/2-\alpha(\alpha+1)/2$.
For the following theorem of Weissauer we refer to \cite{Weis1}.

\begin{theorem}
Let $\tilde{\mathcal A}_g$ be a smooth compactification of 
${\mathcal A}_g$. If $p$ is an integer $0\leq p< g(g+1)/2$ then 
the space of holomorphic $p$-forms on
$\tilde{\mathcal A}_g$ is zero unless $p$ is of the form
$g(g+1)/2-\alpha (\alpha+1)/2$ with $1\leq \alpha \leq g$
and then 
$H^0(\tilde{\mathcal A}_g,\Omega^p_{\tilde{\mathcal A}_g})\cong M_{\rho}(\Gamma_g)$
with $\rho=(g+1,\ldots,g+1,g-\alpha,\ldots,g-\alpha)$ with
$g-\alpha$ occuring $\alpha$ times.
\end{theorem}

If $f$ is a classical Siegel modular form of weight $k=g+1$ on the group
$\Gamma_g$  then $f(\tau) \prod_{i\leq j} d\tau_{ij}$
is a top differential on the smooth part of quotient space
$\Gamma_g \backslash {\mathcal H}_g={\mathcal A}_g$.
It can be extended over the smooth part of the 
rank-$1$ compactification ${\mathcal A}_g^{(1)}$
if and only if $f$ is a cusp form. It is not difficult to see
that this form can be extended as a holomorphic form to the whole
smooth compactification $\tilde{\mathcal A}_g$.

\begin{proposition}\label{top}
The map that associates to a classical cusp form $f
\in S_{g+1}(\Gamma_g)$ of weight $g+1$ the top differential
$\omega= f(\tau) \prod_{i\leq j} d\tau_{ij}$ gives an 
isomorphism between $S_{g+1}(\Gamma_g)$ and the space
of holomorphic top differentials 
$H^0(\tilde{\mathcal A}_g, \Omega^{g(g+1)/2})$
on any smooth compactification
$\tilde{\mathcal A}_g$.
\end{proposition}

For this and an analysis of when the other forms extend over
the singularities in these cases we refer to \cite{Weis1,Fr1}.

Finally we refer to two papers of Salvati-Manni where he proves
the existence of differential forms of some weights, \cite{SM-1,SM-2}
and a paper of Igusa, \cite{Ig6}, where Igusa discusses the question
whether certain Nullwerte of jacobians of odd thetafunctions
can be expressed as polynomials or rational functions in theta Nullwerte.

\end{section}
\begin{section}{Cusp Forms and Geometry}
The very first cusp forms that one encounters often have a beautiful
geometric interpretation. We give some examples.

For $g=1$ the first cusp form is 
$\Delta=\sum \tau(n)q^n \in S_{12}(\Gamma_1)$. \index{Delta}
It is up to a normalization the discriminant $g_2(\tau)^3-27 g_3^2$ of the
equation $y^2=4x^3-g_2x -g_3$ for the Riemann surface $\CC / \ZZ \tau +\ZZ$
and does not vanish on $\H1$. Here $g_2=(4\pi^4/3)E_4(\tau)$
and $g_3=(8\pi^6/27)E_6(\tau)$ are the suitably normalized
Eisenstein series.

For $g=2$ there is a similar cusp form
$\chi_{10}$ of weight $10$ with development 
$$
\chi_{10}\left( \begin{matrix} \tau_{1} & z \\
z & \tau_{2} \\ \end{matrix} \right)=
((\exp 2\pi i \tau_{1} ) \exp(2 \pi i \tau_{2}) + \ldots) 
(\pi z)^2 + \ldots
$$
which vanishes (with multiplicity $2$) along the `diagonal' 
$z=0$.
So its zero divisor in ${\Acal}_2$ is
the divisor of abelian surfaces that are products of elliptic curves
with multiplicity $2$.
There is the Torelli map ${\Mcal}_2 \to {\Acal}_2$
that associates to a hyperelliptic complex
curve of genus $2$ given by $y^2=f(x)$ its Jacobian. Then the pull back of 
$\chi_{10}$ to ${\Mcal}_2$ is related to the discriminant of $f$,
cf.\ Igusa's paper \cite{Ig2} or \cite{RdJ}, Prop.\ 2.2.

For $g=3$ the ring of classical modular forms is generated by
$34$ elements, cf.\ \cite{Tsu}. 
As we saw above, there is a cusp form of weight $18$, namely the
product of the $36$ even theta constants $\theta [\epsilon]$
and its zero divisor is the closure of the hyperelliptic locus.
This expresses the fact that a genus $3$ Riemann surface with a vanishing
theta characteristic is hyperelliptic.

For $g=4$ there is the following beautiful example.
There is up to isometry only one isomorphism class of
even positive definite quadratic forms in $8$ variables, namely $E_8$. In $16$
variables there are exactly two such classes, 
$E_8\oplus E_8$ and $E_{16}$.
To each of these quadratic forms in $16$ variables we can associate
a Siegel modular form on $\Gamma_4$ by means of a theta series:
$\theta_{E_8\oplus E_8}$ and $\theta_{E_{16}}$. The difference
$\theta_{E_8\oplus E_8}-\theta_{E_{16}}$ is a cusp form of weight $8$.
Its zero divisor is the closure of the locus of Jacobians of Riemann surfaces
of genus $4$ in ${\Acal}_4$ as shown by Igusa, cf., \cite{Ig4}.
Here we also refer to \cite{Poor} for a proof.
We shall encounter this form again in Section \ref{lifting}.

Similarly, the theta series associated to the $24$ different Niemeier lattices
(even, positive definite) of rank $24$ produce in genus $12$ a linear
subspace of $M_{12}(\Gamma_{12})$ of dimension $12$. It intersects
the space of cusp forms in a $1$-dimensional subspace, as was
proved in \cite{B-F-W}. We thus find
a cusp form of weight $12$. As we shall see later, it is an Ikeda lift
of the cusp form $\Delta$ for $g=1$ (proven in \cite{B-K}). 
\begin{question} What is the
geometric meaning of this cusp form?
\end{question}

The paper \cite{N-V} contains explicit results on Siegel modular 
forms of weight $12$ obtained from lattices in dimension $24$. For 
example, it gives a non-zero cusp form of weight $12$ on $\Gamma_{11}$,
hence one has a top differential on  $\tilde{\mathcal A}_{11}$,
cf., \ref{top}, implying that this modular variety is not rational
or unirational (cf., \cite{Mum2} where it is proved that $\tilde{\mathcal A}_g$
is of general type for $g>6$).     
\end{section}

\begin{section}{The Classical Hecke Algebra}\label{Heckeoperators}
In the arithmetic theory of elliptic modular forms Hecke operators
play a pivotal role. They enable one to extract arithmetic
information from the Fourier coefficients of a modular form: if
$f=\sum_n a(n)q^n$ is a common eigenform of the Hecke operators
which is normalized ($a(1)=1$) then the eigenvalue $\lambda(p)$ of $f$
under the Hecke operator $T(p)$ equals the Fourier coefficient $a(p)$.

The classical theory of Hecke operators 
\index{Hecke operator} as for example exposed in 
Shimura's book (\cite{Shim2}) can be generalized to the setting of $g>1$
as Shimura showed in \cite{Shim1}, 
though the larger size of the matrices involved is a discouraging aspect of it.
It is worked out in the books \cite{An1,A-Z,Fr1}, 
of which the last, by Freitag, is certainly the most accessible. 
In this section we sketch this approach, in the next section
we give another approach. We refer to loc.\ cit.\ for details.
                                                                               
Recall the group
$
G:={\rm GSp}(2g,\QQ)=\{ \gamma \in {\rm GL}(\QQ^{2g}) :
\gamma^t J \gamma = \eta(\gamma) J, \, \eta(\gamma) \in \QQ^* \}
$
of symplectic similitudes of the symplectic vector space $(\QQ^{2g},
\langle \, , \, \rangle)$.
and $G^{+}= \{ \gamma \in G : \eta(\gamma)>0\}$.

We start by defining the abstract \Definition{Hecke algebra} \index{Hecke
algebra}
$H(\Gamma,G)$ 
for the pair $(\Gamma,G)$ with $\Gamma=\Gamma_g$ and $G={\rm GSp}(2g,\QQ)$.
Its elements are finite formal sums (with $\QQ$-coefficients)
of double cosets $\Gamma \gamma \Gamma$
with $\gamma \in G^{+}$. Each such double coset $\Gamma \gamma \Gamma$
can be written as a finite disjoint union of right cosets $L_i=\Gamma \gamma_i$
by virtue of the following lemma.

\begin{lemma}
Let $m$ be a natural number.
The set $O_g(m)=\{ \gamma \in {\rm Mat}(2g\times 2g,\ZZ) :
\gamma^t J \gamma = m J \}
$
can be written as a finite disjoint union of right cosets. Every right
coset has a representative of the form $(a, b; 0, d)$
with $a^t \, d = m \, 1_g$ and such that $a$ has zeros below the diagonal
\end{lemma}

So to each double coset $\Gamma \gamma \Gamma$ we can associate
a finite formal sum of right cosets.
Let ${\Lcal}$ be the $\QQ$-vector space of
finite formal expressions $\sum_i c_i L_i$ with $L_i=\Gamma \gamma_i$
a right coset and $c_i \in \QQ$.
The map $H(\Gamma,G) \to {\Lcal}$ is injective and induces an
isomorphism $H(\Gamma,G) \cong {\Lcal}^{\Gamma}$, where the action
of $\Gamma$ on ${\Lcal}$ is
$\Gamma \gamma_1 \mapsto \Gamma \gamma_1\gamma$.

We now make this into an algebra by specifying the product of
$\Gamma \gamma \Gamma=\sum_i \Gamma \gamma_i$ and $\Gamma \delta  \Gamma=
\sum_j \Gamma \delta_j$
by
$$
(\Gamma \gamma \Gamma)\cdot (\Gamma \delta  \Gamma)
=\sum_{i,j} \Gamma \gamma_i
\delta_j.
$$

To deal with these double cosets the following proposition is very helpful.

\begin{proposition} {\rm (Elementary divisors)} Let $\gamma \in
{\rm GSp}^+(2g,\QQ)$ be an element
with integral entries. Then double coset $\Gamma \gamma
\Gamma$ has a unique representative of the form
$$
\alpha={\rm diag}(a_1,\ldots,a_g,d_1,\ldots,d_g)
$$
with integers $a_j, d_j$ satisfying
$a_j>0$, $a_jd_j=\eta(\gamma)$ for all $j$, and furthermore
 $a_g|d_g$, $a_j|a_{j+1}$
for $j=1,\ldots,g-1$.
\end{proposition}
                                                                                
On $G$ we have the anti-involution
$$
\gamma= \left(\begin{matrix} a & b \\ c & d \end{matrix} \right)
\mapsto \gamma^{\vee}=\left(\begin{matrix} d^t & -b^t \\ -c^t & a^t \\
\end{matrix} \right)= \eta(\gamma)\, \gamma^{-1}.
$$
(Note that $\eta(\gamma^{\vee})=\eta(\gamma)$.) Another involution is
given by
$$
\gamma \mapsto J \, \gamma \, J^{-1}=\left(\begin{matrix} d & -c \\
-b & a\\ \end{matrix} \right)= \eta(\gamma) \gamma^{-t}.
$$
Because of the proposition we have $\Gamma \gamma \Gamma =
\Gamma \gamma^{\vee} \Gamma$ since we may choose $\gamma$ diagonal
and then $\gamma^{\vee}= J \gamma J^{-1}$ and $J \in \Gamma$.
This implies that for a sum of right cosets $\Gamma \gamma \Gamma =
\sum \Gamma \gamma_i$ we have $\Gamma \gamma \Gamma =
\sum \gamma_i^{\vee} \Gamma$. And it is easy to see that
$\gamma \mapsto \gamma^{\vee}$ defines an anti-involution of
$H(\Gamma, G)$ which acts trivially so that the Hecke algebra is commutative.

We can decompose these diagonal matrices as a product of matrices so that
in each of the factors only powers of one prime occur as non-zero entries.
This leads to a decomposition
$$
H(\Gamma,G)= \otimes_p H_p
$$
as a product of local Hecke algebras \index{local Hecke algebra}
$$
H_p= H(\Gamma, G\cap {\rm GL}(2g, \ZZ[1/p])),
$$
where we allow in the denominators only powers of $p$.
Now $H_p$ has a subring $H_p^0$ generated by
integral matrices. We have $H_p=H_p^0[1/T]$ with $T$ the element defined by
$T=\Gamma_g (p\, 1_{2g}) \Gamma_g$.
By induction one proves the following theorem, cf., \cite{A-Z,Fr1}.

\begin{theorem}\label{generators} 
The local Hecke algebra $H_p^0$ is 
generated by the element
$T(p)$ given by $\Gamma_g \left( \begin{matrix} 1_g & 0_g \\
0_g & p\, 1_g \\ \end{matrix}\right) \Gamma_g$ and the elements
$T_{i}(p^2)$ for $i=1,\ldots,g$ given by
$$
\Gamma_g \left( \begin{matrix} 1_{g-i} &&& \\ & p 1_i \\ && p^2 1_{g-i} \\
&&& p1_i \\ \end{matrix} \right) \Gamma_g
$$
\end{theorem}

For completeness sake
we also introduce the element $T_0(p^2)$ given by the double coset
$\Gamma_p  \left( 
\begin{matrix} 1_g & 0 \\ 0_g & p^2 1_g \\ \end{matrix}
\right) \Gamma_g$.
Note that $T_g(p^2)$ equals the $T=\Gamma_g(p\, 1_{2g})\Gamma_g$
given above.

\begin{definition} Let $T(m)$ be the element of $H(\Gamma,G)$ defined
by the  set $M=O_g(m)$ which is a finite disjoint union of double cosets.
\end{definition}

If $m=p$ is prime then $M=O_g(m)$ is
one double coset and $T(m)$ coincides with $T(p)$, introduced above.
For $m=p^2$ the set $O_g(p^2)$ is a union
of $g+1$ double cosets and the element
$T(p^2)$ is a sum $\sum_{i=0}^{g} T_{i}(p^2)$.

The Hecke algebra can be made to act on the space of Siegel modular
forms $M_{\rho}(\Gamma_g)$. We first define the `slash operator'.
\index{slash operator}

\begin{definition} Let $\rho: {\rm GL}(g,\CC) \to {\rm End}(V)$
be a finite-dimensional irreducible complex representation
corresponding to $(\lambda_1\geq \ldots \geq \lambda_g)$.
For a function $f\colon \Hg \to V$  and an element $\gamma \in 
{\rm GSp}^+(2g,\QQ)$
we set
$$
f|_{\gamma, \rho} (\tau)= \rho(c\tau+d)^{-1} f(\gamma(\tau)) \quad
\gamma=\left( \begin{matrix} * & * \\ c& d\\ \end{matrix} \right)
$$
(or possibly with 
$\eta(\gamma)^{\sum \lambda_i -g(g+1)/2}$, incorporate $\rho$ !).
\end{definition}
Note that
$f|_{\gamma_1, \rho} |_{\gamma_2,\rho} = f|_{\gamma_1\gamma_2,\rho}$.
So if $g>1$ then $f$ is a modular form of weight $\rho$ 
if and only if $f$ is holomorphic
and $f|_{\gamma,\rho} =f$ for every $\gamma\in \Gamma_g$.

Let now $M \subset {\rm GSp}(2g,\QQ)$ be a subset satisfying the two properties
\begin{enumerate}
\item $M=\sqcup_{i=1}^h \Gamma_g \, \gamma_i$ is a finite disjoint
union of right cosets $\Gamma_g \gamma_i$;
\item $M \, \Gamma_g\subset M$.
\end{enumerate}
The first condition implies that if for $f \in M_{\rho}$ we set
$$
T_M f := \sum_{i=1}^h f|_{\gamma_i,\rho}
$$
then this is independent of the choice of the
representatives $\gamma_i$, while the second condition implies that
$(T_M f)|_{\gamma}=T_M f$ for all $\gamma \in \Gamma_g$. Together
these conditions imply that $T_M$ is a linear operator on the space $M_{\rho}$.
                                                                                
Double cosets $\Gamma\gamma \Gamma$ satisfy condition 2) if $\Gamma =\Gamma_g$
and $\gamma \in {\rm Sp}(2g,\QQ)$ and also condition 1) by what was said
above.

An important observation is 
that 
$\langle Tf,g\rangle = \langle f, 
T^{\vee}g\rangle$, where $\langle \, , \, \rangle $ 
gives the Petersson product  and
thus the Hecke operators define Hermitian operators on the space
of cusp forms $S_{\rho}$.

Just as in the classical case $g=1$ we can associate correspondences
(i.e.\ divisors on ${\Acal}_g \times {\Acal}_g$)
to Hecke operators. The correspondence associated to $T_p$
sends a principally polarized abelian variety $X$
to the sum $\sum X^{\prime}$ of principally polarized $X^{\prime}$ which admit
an isogeny $X \to X'$ with kernel an isotropic (for the Weil pairing)
subgroup $H \subset X[p]$
of order $p^g$. Similarly, the correspondence associated to $T_{i}(p^2)$
sends $X$ to the sum $\sum X'$ with the $X'$ quotients $X/H$, where $H
\subset X[p^2]$ is an isotropic subgroup of order $p^{2g}$ with
$H \cap X[p]$ of order $p^{g+i}$.
\end{section}

\begin{section}{The Satake Isomorphism}
We can identify the local Hecke algebra $H_p$ with the $\QQ$-algebra of 
$\QQ$-valued locally constant functions on ${\rm GSp}(2g,\QQ_p)$
with compact support and which are invariant under the (so-called
hyperspecial maximal compact) subgroup $K=
{\rm GSp}(2g,\ZZ_p)$
acting both from the left and right. The multiplication in this
algebra is convolution 
$f_1 \cdot f_2= \int_{{\rm GSp}(2g,\QQ_p))} f_1(g)f_2(g^{-1}h) dg$, where
$dg$ denotes the unique Haar measure normalized such that the volume 
of $K$ is $1$.
The correspondence is obtained by sending the double coset
$K \gamma K$ to the characteristic function of $K \gamma 
K$. A compactly supported function in $H_p$ is constant on
double cosets and its support is a finite linear combination
of characteristic functions of double cosets.

Note that Proposition \ref{generators} 
tells us that $H_p$ is generated by
the double cosets of diagonal matrices.
In order to describe this algebra conveniently we compare it with the
$p$-adic Hecke algebras of two subgroups, the diagonal torus and the
Levi subgroup of the standard parabolic subgroup.

To be precise, recall the diagonal torus $\TT$ of
${\rm GSp}(2g,\QQ)$ isomorphic to $\GG_m^{g+1}$ and the Levi subgroup
$$
M= \{ \left( \begin{matrix} a& 0 \\ 0 & d \\ \end{matrix} \right) 
\in {\rm GSp}(2g,\QQ) \}
$$
of the standard parabolic $Q=\{ (a,b;0,d) \in {\rm GSp}(2g,\ZZ)\}$
that stabilizes the first summand $\ZZ^g$ of $\ZZ^g\oplus \ZZ^g$.
In particular for an element $(a,0;0,d) \in M$ we have $ad^t=\eta$ 
and the group $M$ is isomorphic to 
${\rm GL}(g)\times \GG_m$. Let $Y\cong \ZZ^{g+1}$ be the co-character group 
of $\TT_m$, i.e., $Y={\rm Hom}(\GG_m,\TT)$, cf.\ Section \ref{roots}. 

We can construct a local Hecke algebra $H_p(\TT)=H_p(\TT,\TT_{\QQ})$ 
for the group $\TT$ too
as the $\QQ$-algebra of $\QQ$-valued, bi-$\TT(\ZZ_p)$-invariant, locally
constant functions with compact support on $\TT(\QQ_p)$. This local 
Hecke-algebra is easy to describe: $H_p(\TT)\cong \QQ[Y]$, 
the group algebra over $\QQ$ of $Y$ where $\lambda \in Y$
corresponds to the characteristic function of the double coset
$D_{\lambda}=K \lambda(p) K$. Concretely, 
$H_p(\TT)$ is isomorphic to the ring
$\QQ[(u_1/v_1)^{\pm},\ldots, (u_g/v_g)^{\pm},
(v_1\cdots v_g)^{\pm}]$ under a map that sends $(a_1,\ldots,a_g,c)$ to 
the element
$(u_1/v_1)^{a_1} \cdots (u_g/v_g)^{a_g}(v_1 \cdots v_g)^c$.

Similarly, we have a $p$-adic Hecke algebra $H_p(M)=H_p(M,M_{\QQ})$ for $M$.

Recall that the Weyl group 
$W_G=N(\TT)/\TT$, with $G={\rm GSp}(2g,\QQ)$
and $N(\TT)$ the normalizer of $\TT$ in $G$,
acts. This group
$W_G$ is isomorphic to $S_g \ltimes (\ZZ/2\ZZ)^g$, where the generator
of the $i$-th factor $\ZZ/2\ZZ$ acts on a matrix of the form
${\rm diag}(\alpha_1,\ldots,\alpha_g,\delta_1,\ldots,\delta_g)$ by
interchanging $\alpha_i$ and $\delta_i$ and the symmetric group $S_g$
acts by permuting the $\alpha$'s and $\delta$'s. The Weyl group
of $M$ (normalizer this time in $M$)
is isomorphic to the symmetric group $S_g$. The algebra of invariants
$H_p(\TT)^{W_G}$ is of the form $\QQ[y_0^{\pm},y_1,\ldots,y_g]$, cf. 
\cite{Fr1}. 

We now give Satake's so-called spherical map of the Hecke algebra 
$H_p(\Gamma,G)$ to the Hecke algebras
$H_p(M)$ and $H_p(\TT)$, cf., \cite{Sa,Ca,F-C,Gross}. 
The images will land in the $W_M$-invariant (resp.\ the $W_G$-invariant) 
part.

We first need the following characters. 
The Borel subgroup $B$ of matrices $(a,b;0,d)$ with $a$ upper triangular
and $d$ lower triangular determines a set $\Phi^{+}$ of positive roots
in the set of all roots $\Phi$
(= characters that occur in the adjoint
representation of $G$ on ${\rm Lie}(B)$). We let $2\rho=\sum_{\Phi^+}
\alpha$.

Define $ e^{2\rho_n}: M \to \GG_m $
by $\gamma=(a,0; 0,d) \mapsto \det(a)^{g+1} \eta(\gamma)^{-g(g+1)/2}$,
where the multiplier $\eta(\gamma)$ is defined by $a\cdot d^t=\eta(\gamma) 1_g$.
(This corresponds to the adjoint action of $\TT$ on the Lie algebra
of the unipotent radical of $P$.)
Secondly, we have the character $e^{2\rho_M}: \TT \to \GG_m$ given by
$$
{\rm diag}(\alpha_1,\ldots,\alpha_g,\delta_1,\ldots,\delta_g)
\mapsto \prod_{i=1}^g \alpha_i^{g+1-2i}=\prod_{i=1}^g \delta^{2i-(g+1)}.
$$
and $2\rho_M$ is the sum of the positive roots
in $\Phi_M^{+}=\{ a_i/a_j: 1\leq i < j \leq g\}$.
Together they give a character $e^{2\rho} : \TT \to \GG_m$ given by
$e^{2\rho(t)}=e^{2\rho_n(t)}e^{2\rho_M(t)}$ for $t \in \TT$; explicitly,
$$
{\rm diag}(\alpha_1,\ldots,\alpha_g,\delta_1,\ldots,\delta_g)
\mapsto \eta^{-g(g+1)/2} \prod_{i=1}^g \alpha_i^{2g+2-2i}.
$$

Satake's spherical map \index{Satake's spherical map}
\index{Satake isomorphism}
$
S_{G,M}: H_p(\Gamma,G) \to H_p(M)
$
is defined by integrating
$$
S_{G,M}(\phi)(m)= | e^{\rho_n(m)}| \int_{U(\QQ_p)} \phi(mu)du,
$$
where $|p|=1/p$.
Similarly, we have a map
$$
S_{M,T} : H_p(M) \to H_p(\TT)
$$
given by
$$
S_T(\phi)(t)=|e^{\rho_{M}(t)}| \int_{M\cap N} \phi(tn)dn.
$$
In \cite{F-C} the authors define a `twisted' version of these spherical
maps where they put $| e^{2\rho_n(m)}|$ and $|e^{2\rho_{M}(t)}|$ instead
of the multipliers above. In this way one avoids square roots of $p$.
If one uses this twisted version one 
should also twist the action of the Weyl group on the co-character group 
$Y$ of $\TT$ by $e^{\rho}$ too: in the usual action $S_g$ permutes the $a_i$
and $d_i$ and the $i$-th generator $\tau_i$ of $(\ZZ/2\ZZ)^g$ 
interchanges $a_i$ and $d_i$. Under the twisted action $\tau_i$
sends $(u_i,v_i)$ to $(p^{g+1-i}v_i,p^{i-g-1}u_i)$, while the permutation
$(i \, i+1) \in S_g$ sends $(u_i/v_i)$ to $p u_{i+1}/v_{i+1}$. The formula
is
$
w\cdot \phi (t)= |e^{\rho(w^{-1}t)-\rho(t)}| \phi(w^{-1}t)
$
for $w \in W$ and $t \in \TT$, cf., \cite{F-C}.

The basic result is the following theorem. 

\begin{theorem} 
Satake's spherical maps  $S_{G,M}$ and $S_{M,T}$ define
isomorphisms of $\QQ$-algebras $H_p(G){\buildrel \sim \over
\longrightarrow}H_p(\TT)^{W_G}$ and $H_p(M){\buildrel \sim \over
\longrightarrow}H_p(\TT)^{W_M}$.
\end{theorem}

For the untwisted version there is a similar result but one 
needs to tensor with $\QQ(\sqrt{p})$.
One can calculate these maps explicitly. 
A right coset $K\lambda(p)$ with $\lambda \in Y$ is mapped 
under $S_{GT}$ to $p^{\langle \lambda,\rho\rangle} \lambda$.
Concretely, if $\gamma={\rm diag}
(p^{\alpha_1},\ldots,p^{c-\alpha_g})$ then $S_{G,T}(K\gamma)$ equals
$$p^{cg(g+1)/4} (v_1\cdots v_g)^c \prod_{i=1}^g (u_i/p^iv_i)^{
\alpha_i}.$$
If we write a double coset $K \lambda(p) K$
as a finite sum of right cosets  $K\gamma$ then we may take
$\gamma=\lambda(p)$ as one of these coset representatives.
Then the image of the double coset $K \lambda(p) K$
is a sum $p^{\langle \lambda,\rho\rangle} \lambda +
\sum_{\mu} n_{\lambda,\mu} \mu$ where the $\mu$ satisfy $\mu <
\lambda$ (i.e.\ $\lambda-\mu$ is positive on $\Phi^+$)
and the $n_{\lambda, \mu}$ are non-negative integers, cf., 
\cite{Ca,Gross}.
\end{section}

\begin{section}{Relations in the Hecke Algebra}\label{relations}
We derive some relations in the Hecke algebras. 
We first define elements $\phi_i$ in the Hecke
algebra $H_p(M)$ by
$$
p^{i(i+1)/2} \phi_i= M(\ZZ_p) \left( 
\begin{matrix} 1_{g-i} && \\ & p1_g & \\ && 1_i \\
\end{matrix} \right) M(\ZZ_p) \qquad i=0,\ldots,g
$$
From \cite{A-Z}, p.\ 142--145 one can derive the following result.

\begin{proposition}\label{phi} We have
$
S_{G,M}(T_p)= \sum_{i=0}^g \phi_i
$ and  for $i=1,\ldots,g$
$$
S_{G,M}(T_{i}(p^2))= \sum_{j,k\geq 0, j+i\leq k}^g  m_{k-j}(i) 
p^{-{k-j+1\choose 2}}\phi_j \phi_k,
$$
where $m_h(i)=\# \{A \in {\rm Mat}(h \times h , \FF_p)
\colon A^t=A, {\rm corank }(A)=i \}$.
Moreover, for $i=0,\ldots,g$ we have
$$
S_{M,T}(\phi_i)= (v_1\cdots v_g) \sigma_i(u_1/v_1,\ldots,u_g/v_g),
$$
where $\sigma_i$ denotes the elementary symmetric function of degree $i$.
\end{proposition}
\begin{example}
$g=1$. We have $T(p)\mapsto \phi_0+\phi_1$, $T_{0}(p^2) \mapsto
\phi_0^2+((p-1)/p) \phi_0\phi_1 +\phi_1^2$ and
$T_{1}(p^2) \mapsto \phi_0 \phi_1/p$. We derive that $T(p^2)=
T_0(p^2)+T_1(p^2)$ satisfies the well-known relation
$T(p^2)=T(p)^2-p T_{1}(p^2)$.

$g=2$. 
We find $T(p) \mapsto  \phi_0+ \phi_1+ \phi_2$ and
$T_{1}(p^2)\mapsto  \frac{1}{p} \phi_0 \phi_1+ \frac{p^2-1}{p^3} 
\phi_0 \phi_2 + 
\frac{1}{p} \phi_1 \phi_2$ and similarly $
T_2(p^2)\mapsto  \frac{1}{p^3} \phi_0 \phi_2$.
\end{example}
We denote the element $\phi_0$ corresponding to $(1_g,0;0,p 1_g)$ by
${\rm Frob}$\index{Frobenius}. This element of $H_p(M)$ generates the fraction
field of $H_p(M)$ over the fraction field of $H_p(\Gamma,G)$
as we can see from the calculation above. Indeed, we have that
$S_T(\phi_0)=v_1\cdots v_g$ and this element of $H_p(\TT)$
is fixed by $S_g$, but not by any other element of $W_G$. In particular,
it is a root of the polynomial
$$
\prod_{w \in (\ZZ/2\ZZ)^g} (X-w(\phi_0))= \prod_{I \subset \{ 1,\ldots,g\}}
(X- \prod_{i \in I} u_i \prod_{i \notin I} v_i).
$$
For example, for $g=1$ we find by elimination
 that $\phi_0$ is a root of 
$$
X^2-T_p X+p T_{p,1},
$$
while for $g=2$ we have that $\phi_0$ is a root of
$$
X^4-T(p) X^3 +(p\, T_{1}(p^2) +(p^3+p)T_2(p^2))X^2-p^3\,
T(p)T_{2}(p^2)X+p^6\, T_{2}(p^2)^2.
$$
Using the relation
$$
T(p)^2=T_0(p^2)+(p+1)T_1(p^2)+(p^3+p^2+p+1) T_2(p^2)
$$
this can be rewritten as a polynomial $F(X)$
given by
$$
X^4-T(p) X^3 +(T(p)^2-T(p^2)-p^2 T_{2}(p^2))\, X^2-p^3\,
T(p)T_{2}(p^2)\, X+p^6\, T_{2}(p^2)^2.
$$
Moreover, in the power series ring over the 
Hecke ring of ${\rm Sp}(4,\QQ)$ one has the formal relation (cf., 
\cite{Shim1}, \cite{A-Z}, p.\ 152)
$$
\sum_{i=0}^{\infty} T(p^i)\, z^i=
\frac{ 1-p^2 \, T_2(p^2) \, z^2}{z^4 F(1/z)}.
$$
For a slightly different approach we refer to a paper \cite{Krieg} by Krieg
and a preprint by Ryan with an algorithm to
calculate the images, cf., \cite{Ryan}.
\end{section}
\begin{section}{Satake Parameters}\label{satake-parameters}
The usual argument that uses the Petersson product shows that the
spaces $S_{\rho}$ possess a basis of common eigenforms for the action of
the Hecke algebra.

If $F$ is a Siegel modular form in $M_{\rho}(\Gamma_g)$
for an irreducible representation $\rho=(\lambda_1,\ldots,\lambda_g)$ 
of ${\rm GL}(g,\CC)$ which is an eigenform of
the Hecke algebra $H$ then we get for each Hecke operator $T$
an eigenvalue $\lambda_F(T) \in \CC$, a real algebraic number. 
Now the determination of 
the local Hecke algebra $H_p \otimes \CC \cong \CC[Y]^{W_G}$
says that 
$$
{\rm Hom}_{\CC}(H_p,\CC) \cong (\CC^*)^{g+1}/W_G.
$$
In particular, for a fixed eigenform $F$
the map $H_p \to \CC$ given by $T \mapsto \lambda_F(T)$
is determined by (the $W_G$-orbit of) a $(g+1)$-tuple 
$(\alpha_0,\alpha_1,\ldots,\alpha_g)$ of non-zero complex numbers,
the \Definition{$p$-Satake parameters} 
\index{Satake parameter} of $F$. So for $i=1,\ldots,g$
the parameter $\alpha_i$ is the
image of $u_i/v_i$ and $\alpha_0$ that of $v_1\cdots v_g$
and $\tau_i \in W_G$ acts by $\tau_i(\alpha_0)=\alpha_0\alpha_i$,
$\tau_i(\alpha_i)=1/\alpha_i$ and $\tau_i(\alpha_j)=\alpha_j$ if $j\neq 0, i$. These
Satake parameters satisfy the relation
$$
\alpha_0^2\alpha_1\cdots \alpha_g=p^{\sum_{i=1}^g \lambda_i-(g+1)g/2}.
$$
This follows from the fact that $T_g(p^2)$, which corresponds to the
double coset of $p \cdot 1_{2g}$, is mapped to
$p^{-g(g+1)/2} (v_1\cdots v_g)^2 \, \prod_{i=1}^g (u_i/v_i)$
as we saw above.

For example, if $f=\sum_n a(n) q^n \in S_k(\Gamma_1)$ is a normalized eigenform
and if we write $a(p)=\beta+\bar{\beta}$ with $\beta\bar{\beta}=p^{k-1}$
then $(\alpha_0,\alpha_1)=(\beta, \bar{\beta}/\beta)$ or
$(\alpha_0,\alpha_1)=(\bar{\beta},\beta/\bar{\beta})$. 
Or if $f\in M_k(\Gamma_g)$ 
is the Siegel Eisenstein series of weight $k$ then
the Satake parameters at $p$ are: $\alpha_0=1$,
$\alpha_i= p^{k-i}$ for $i=1,\ldots,g$.

The formulas from Proposition 
\ref{phi} give now formulas for the eigenvalues of the
Hecke operators $T(p)$ and $T_i(p^2)$ in terms of these Satake parameters:
$$
\lambda(p)= \alpha_0 (1+\sigma_1+\ldots + \sigma_g)
$$
and similarly
$$
\lambda_i(p^2)= 
 \sum_{j,k\geq 0, j+i\leq k}^g  m_{k-j}(i)
p^{-{k-j+1\choose 2}} \alpha_0^2 \sigma_i \sigma_j,$$
where $\sigma_j$ is the $j$th elementary symmetric function 
in the $\alpha_i$ with $i=1,\ldots,\alpha_g$ and the $m_h(i)$ are
defined as in \ref{phi}.
\end{section}
\begin{section}{L-Functions} \index{L-functions}
It is customary to associate to an eigenform $f=\sum a(n) q^n 
\in M_k(\Gamma_1)$ of the Hecke algebra a Dirichlet series
$\sum_{n \geq 1} a(n) n^{-s}$ with $s$ a complex parameter
whose real part is $>k/2+1$. It is well-known that for a cusp form
this L-function
admits a holomorphic continuation to the whole $s$-plane and
satisfies a functional equation. The multiplicativity properties
of the coefficients $a(n)$ ensure that we can write it formally
as an Euler product
$$
\sum_{n>0} a(n) n^{-s}=\prod_p (1-a(p)p^{-s}+p^{k-1-2s})^{-1}.
$$
In defining $L$-series for Siegel modular forms one uses Euler products.

Suppose now that $f\in M_{\rho}(\Gamma_g)$ is an eigenform of the Hecke
algebra with eigenvalues $\lambda_f(T)$ for $T \in H_{p}^0$. Then the
assignment $T \mapsto \lambda_f(T)$ defines an element of
${\rm Hom}_{\CC}(H_{p}^0,\CC)$. We called the corresponding $(g+1)$-tuple 
of $\alpha$'s  the $p$-Satake parameters of $f$. 
The fact that $\ZZ[Y]^{W_G}$ is also the representation ring
of the complex dual group $\hat{G}$ of $Gi={\rm GSP}(2g,\QQ)$ 
(determined by the dual `root datum') 
is responsible for a connection with $L$-functions.
In our case we can use the Satake parameters to
define the following formal $L$-functions. Firstly, there is the
\Definition{spinor zeta function} $Z_f(s)$ 
\index{spinor zeta function}
with as Euler factor at $p$ the expression
$Z_{f,p}(p^{-s})^{-1}$ defined by
$$
(1-\alpha_0t)\prod_{r=1}^g \prod_{1 \leq i_1 < \cdots < i_r\leq g}
(1-\alpha_0\alpha_{i_1} \ldots \alpha_{i_r}t)
=(1-\alpha_0t)\prod_I (1-\alpha_0\alpha_It),
$$
where the product has $2^g$ factors corresponding to the
$2^g$ subsets $I \subseteq \{ 1,\ldots,g\}$. Secondly, there is the 
\Definition{standard zeta function} with as Euler factor 
$D_{f,p}(p^{-s})^{-1}$
at $p$ the expression
$$
D_{f,p}(t)=(1-t)\prod_{i=1}^g (1-\alpha_{i}t)(1-\alpha_{i}^{-1}t).
$$
\index{standard zeta function}
For example, for $g=1$ the spinor zeta function is
 $Z_f(s)=\sum a(n)n^{-s}$, the usual $L$-series
and the standard zeta function 
$D_{f}(s-k+1)=\prod(1+p^{-s+k-1})^{-1} \, \sum a(n^2)n^{-s}$,
that is related to the Rankin zeta function.
For $g=2$ and eigenform $f \in M_{j,k}(\Gamma_2)$ with $T(m)f=\lambda_f(m) f$
we have $Z_f(s)=\zeta(2s-j-2k+4)\sum_{m \in \ZZ_{>0}}
 \lambda_f(m) m^{-s}$.

We set
$$\Delta(f,s)=(2\pi)^{-gs} \pi^{-s/2} \Gamma(\frac{s+\epsilon}{2})
\prod_{j=1}^g \Gamma(s+k-j) D(f,s),
$$
where $\epsilon=0$ for $g$ even and $\epsilon=1$ for $g$ odd. Then
the function $\Delta(f,s)$ can be extended meromorphically
to the whole $s$-plane and satisfies a functional equation
$\Delta(f,s)=\Delta(f,1-s)$, cf.\ papers by B\"ocherer \cite{Bo2}, 
Andrianov-Kalinin \cite{A-K},
Piatetski-Shapiro and Rallis \cite{P-R}. 
If $f\in S_k(\Gamma_g)$ is a cusp form and $k \geq g$
then $\Delta(f,s)$ is holomorphic except for simple poles at $s=0$ and $s=1$.
It is even holomorphic if the eigenform does not lie in the space generated by
theta series coming from unimodular lattices of rank $2g$.
Also for $k<g$ we have information about the poles, cf., \cite{Mi}.
Andrianov proved that for $g=2$ the function
$\Phi_f(s) = \Gamma(s)\Gamma(s-k+2) (2\pi)^{-2s} Z_f(s)$ is
meromorphic with only finitely many poles and satisfies a functional
equation $\Phi_f(2k-2-s)= (-1)^k \Phi_f(s)$.

One instance where spinor zeta functions associated to Siegel
classical modular forms of weight $2$ occur is as $L$-functions
associated to the $1$-dimensional cohomology of simple abelian
surfaces.

We end by giving two additional references: the lectures notes by 
Courtieu and Panchishkin
\cite{CourPan} and a paper \cite{Yoshida} by Yoshida on motives associated to
Siegel modular forms.
\end{section}
\begin{section}{Liftings}\label{lifting} \index{lifting}
It is well-known that for a normalized cusp form which is an eigenform
$f=\sum_{n\geq 1} a(n)q^n$ 
of weight $k$ on $\Gamma_1$ we have 
the inequality $|a(p)| \leq 2p^{(k-1)/2}$ for every prime $p$, or equivalently,
the roots of the Euler factor $1-a(p)X+p^{k-1}X^2$ at $p$ have absolute
value $p^{-(k-1)/2}$. This was shown by Eichler for cusp forms of weight $k=2$
on the congruence subgroups $\Gamma_0(N) \subset {\rm SL}(2,\ZZ)$ and by
Deligne for general $k$ in two steps, by first reducing it to the 
Weil conjectures in 1968 (\cite{De1}) and then by proving the Weil 
conjectures in 1974. 

For $g=2$ the analogous Euler factor at $p$ for an eigenform $F$ of
the Hecke algebra is the expression
$$
{\Fcal}_p=1-\lambda(p)X+(\lambda (p)^2-\lambda (p^2)-p^{2k-4})X^2
-\lambda (p)p^{2k-3}X^3+p^{4k-6}X^4,
$$
with $\lambda(p)$ the eigenvalue of the cusp form $F\in S_k(\Gamma_2)$;
cf., the polynomial at the end of Section \ref{relations}.
The tacit assumptions of many mathematicians in the 1970's 
was that the absolute values
of the roots of ${\Fcal}_p$ were equal to $p^{-(2k-3)/2}$. For example,
for $k=3$ a classical cusp form $F$ of weight $3$ on a congruence
subgroup $\Gamma_2(n)$ with $n\geq 3$ 
determines a holomorphic $3$-form 
$F(\tau)\prod_{i\leq j} d\tau_{ij}$ on the complex
$3$-dimensional manifold $\Gamma_2(n) \backslash \H2$ that can be extended to
a compactification and we thus find an element of 
the cohomology group $H^3$, so we expect
to find absolute value $p^{-3/2}$. 
But then in 1978 Kurokawa and independently H.\ Saito (\cite{Ku})
found examples
of Siegel modular forms of genus $2$ contradicting this expectation.
Their examples are the very first examples that one encounters, like
the cusp form $\chi_{10} \in S_{10}(\Gamma_2)$. On the basis of
explicit calculations Kurokawa guessed that
$$
L(\chi_{10},s)=\zeta(s-9)\zeta(s-8) L(f_{18},s),
$$
with $f_{18}=\Delta\, e_6 \in S_{18}(\Gamma_1)$ the normalized cusp form of
weight $18$ on ${\rm SL}(2,\ZZ)$. For example, he found for $p=2$
$$
{\Fcal}_2= (1-2^8X)(1-2^9X)(1+528\, X + 2^{17} X^2)
$$
giving the absolute values $p^8$, $p^9$ and $p^{17/2}$ for the inverse
roots.
The examples he worked out suggested that in these cases
$L(F_k,s)=\zeta(s-k+1)\zeta(s-k+2)
L(f_{2k-2},s)$ with $f_{2k-2}\in S_{2k-2}(\Gamma_1)$ a normalized cusp
form and $F_k$ a corresponding Siegel 
modular form of weight $k$ which is an eigenform of the Hecke
algebra. 
On the basis of this he conjectured the existence of a `lift'
\index{lifting}
$$
S_{2k-2}(\Gamma_1) \longrightarrow S_k(\Gamma_2), \quad f \mapsto F
$$
with $L(F,s)=\zeta(s-k+1)\zeta(s-k+2)L(f,s)$. A little later, Maass
identified in $M_k(\Gamma_2)$ a subspace (`Spezialschar', 
\index{Spezialschar} nowadays
called the Maass subspace, cf., \cite{Maa2} ) \index{Maass subspace}
consisting of modular forms 
$F$ with a Fourier
development $F=\sum_{N \geq 0} a(N) e^{2 \pi i {\rm Tr}N\tau}$ satisfying
the property that $a(N)$ depends only on the discriminant $d(N)$ and
the content $e(N)$, i.e., if we write
$$
N=\left( \begin{matrix} n & r/2 \\ r/2 & m \end{matrix} \right)
$$
then $N$ corresponds to the positive definite quadratic form
$[n,r,m]:= nx^2+rxy+my^2$ with discriminant $d=4mn-r^2$ and
content $e={\rm g.c.d.}(n,r,m)$. We shall write $a([n,r,m])$ for $a(N)$.
The condition that $F$ belongs to the Maass space can be formulated 
alternatively as 
$$
a([n,r,m])=\sum_{d>0,\, d|(n,r,m)} d^{k-1} a([1,r/d,mn/d^2])
$$
We shall write $M^*_k(\Gamma_2)$ or $S^*_k(\Gamma_2)$ for the Maass
subspace of  $M_k(\Gamma_2)$ or $S_k(\Gamma_2)$. It was then conjectured
(`Saito-Kurokawa Conjecture') 
\index{Saito-Kurokawa Conjecture} 
that there is a 1-1 correspondence between
eigenforms in $S_{2k-2}(\Gamma_1)$ and eigenforms in the Maass space
$S^*_k(\Gamma_2)$
given by an identity between their 
$L$-functions. More precisely, we now have the
following theorem.

\begin{theorem}\label{spezialschar} The Maass
 subspace $S_k^*(\Gamma_2)$ is invariant under
the action of the Hecke algebra and there is a 1-1 correspondence between
eigenspaces in $S_{2k-2}(\Gamma_1)$ and Hecke eigenspaces in
$S^*_k(\Gamma_2)$ given by 
$$
f \leftrightarrow F \quad \iff \quad L(F,s)=\zeta(s-k+1)\zeta(s-k+2)L(f,s)
$$
with $L(F,s)$ the spinor $L$-function of $F$.
\end{theorem}

The main part of the theorem is due to Maass, but it was completed by 
Andrianov and Zagier, see \cite{Maa2,An2,Zag1}.

We can make an extended picture as follows.
The map $ F \mapsto \phi_{k,1}$ that sends a Siegel modular form to its first
Fourier-Jacobi coefficient induces an isomorphism $M_{k}^*(\Gamma_2)
\cong J_{k,1}$, the space of Jacobi forms, 
and  the map $h=\sum c(n)q^n \mapsto 
\sum_{n\equiv -r^2 \, (\bmod \, 4)} c(n)q^{n+r^2)/4} \zeta^r$
gives an isomorphism of the Kohnen plus space 
$M_{k-1/2}^{+}$ with $ J_{k,1}$ 
fitting in a diagram
$$
\begin{matrix}
M_k^*(\Gamma_2) & {\buildrel \sim \over \longrightarrow}
 & J_{k,1} & {\buildrel \sim \over \longleftarrow} & 
M_{k-1/2}^{+} \\
&&&&\downarrow{\cong} \\
&&&& M_{2k-2}(\Gamma_1)\\
\end{matrix}
$$
where the vertical map is the Kohnen isomorphism
\index{Kohnen isomorphism}. 
Note that the vertical map
is quite different from the horizontal two maps. The vertical isomorphism is notcanonical at all, but depends on the choice of a discriminant $D$.

We now sketch a proof of theorem \ref{spezialschar}.
A classical Siegel modular
form  $F \in M_k(\Gamma_2)$ has a Fourier-Jacobi series
$F(\tau,z,\tau')=\sum \phi_m(\tau,z)e^{2 \pi i m \tau'}$ with
$\phi_m(\tau,z) \in J_{k,m}$, the space of Jacobi forms of weight $k$
and index $m$. The reader may check this by himself. We have on the
Jacobi forms a sort of Hecke operators $V_m \colon J_{k,m} \to J_{k,ml}$
with $\phi|_{k,m} V_l (\tau,z)$ given explicitly by
$$
l^{k-1}\sum_{\Gamma_1 \backslash O(l)}
(c\tau +d)^{-k} e^{2 \pi i m l (-cz^2/(c\tau +d))} \phi((a\tau+b)/(c\tau+d),
lz/(c\tau+d)).
$$
On coefficients, if $\phi=\sum_{n,r} c(n,r) q^n \zeta^r$ then
$$
\phi|_{k,m} V_l = \sum_{n,r} \sum_{a|(n,r,l)} a^{k-1}c(nl/a^2,r/a)q^n\zeta^r.
$$
One now checks using generators of $\Gamma_2$ that for $\phi \in J_{k,1}$
the expression
$$
v(\phi):=\sum_{m\geq 0} (\phi |V_m) (\tau,z) e^{2\pi i m \tau'}
$$
is a Siegel modular form in $M_k(\Gamma_2)$. 

We then have a map $M_k(\Gamma_2) \to \oplus_{m=0}^{\infty} J_{k,m}$
by associating to a modular form its Fourier-Jacobi coefficients;
we also have a map in the other direction $J_{k,1}\to M_k(\Gamma_2)$
given by $\phi \to v(\phi)$ and the composition
$$
J_{k,1} \to M_k(\Gamma_2) \to \oplus_m J_{k,m} {\buildrel {\rm pr}
\over \longrightarrow} J_{k,1}
$$
is the identity. So $v: J_{k,1} \to M_k(\Gamma_2)$ is injective
and the image consists of those modular forms $F$ with the property
that $\pi_m= \phi_1|V_m$. This implies the following relation for the
Fourier coefficients for $[n,r,m]\neq [0,0,0]$
$$
a([n,r,m])=\sum_{d|(n,r,m)} d^{k-1} c((4mn-r^2)/d^2),
$$
where $C(N)$ is given by
$$
c(N)=\begin{cases} a([n,0,1]) & N=4n \\ a([n,1,1]) & N=4n-1. \end{cases}
$$
In particular, we see that the image is the Maass subspace because
$$
a([n,r,m])=\sum_{d|(n,r,m)} d^{k-1} a([nm/d^2,r,1]).
$$
On the other hand, it is known that $J_{k,1} \cong M^{+}_{k-1/2}$.
Combination of the two isomorphisms yields what we want.

Duke and Imamo{\v g}lu conjectured in \cite{D-I1} 
a generalization of this and some evidence
was given by Breulmann and Kuss \cite{B-K}. Then Ikeda 
\index{Ikeda lift} generalized
the Saito-Kurokawa lift of modular forms from one variable 
to Siegel modular forms of
degree $2$ in \cite{Ik} in 1999 under the condition
that $g \equiv k \, (\bmod \, 2)$ to a lifting
from an eigenform $f \in S_{2k}(\Gamma_1)$ to an 
eigenform $F \in  S_{g+k}(\Gamma_{2g})$
such that the standard zeta function of $F$ is given in terms of 
the usual $L$-function of $f$ by
$$
\zeta(s) \prod_{j=1}^{2g} L(f,s+k+g-j).
$$
The Satake parameters of $F$ are $\beta_0,\beta_1,\ldots,\beta_{2g}$
with
$$
\beta_0=p^{gk-g(g+1)/2}, \, \beta_i=\alpha \, p^{i-1/2}, \, 
\beta_{g+i}=\alpha^{-1}p^{i-1/2} \quad {\rm for} \quad i=1,\ldots,g
$$
with $f=\sum a(n)q^n$ and 
$$
(1-\alpha p^{k-1/2}X)(1-\alpha^{-1}p^{k-1/2}X)=1-a(p)X+p^{2k-1}X^2,
$$
cf., \cite{Murokawa}.
(In particular, such lifts do not satisfy the Ramanujan inequality.)
Kohnen (\cite{Ko1}) has interpreted it as an explicit linear map
$S_{k+1/2}^{+} \longrightarrow S_{k+g}(\Gamma_{2g})$ given by 
$$
f=\sum_{ (-1)^kn
\equiv 0,1 (\bmod 4)} c(n)q^n F\mapsto
\sum_{N} a(N) e^{2 \pi {\rm Tr} i N\tau},
$$
with $a(N)$ given by an expression 
$\sum_{a|f_N} a^{k-1}\varphi(a,N)c(|D_N|/a^2)$ and $\phi(a,N)$ an
explicitly given integer-valued numbertheoretic function.

One defines also a Maass space with $M_k^{*}(\Gamma_g)$ consisting of $F$ such
that $a(N)=a(N^{\prime})$ if the discriminants of $N$ and $N^{\prime}$ are
the same and in addition $\phi(a,N)=\phi(a,N^{\prime})$ for all divisors 
$a$ of $f_N=f_{N^{\prime}}$. 
Under the additional assumption that $g \equiv 0, \, 1 \, (\bmod \, 4)$ 
Kohnen and Kojima prove in \cite{K-K}
that the image of the lifting is the Maass space.
\begin{example} Let $k=6$ and $g=2$. Then the Ikeda lift is a map from
$S_{12}(\Gamma_1) \to S_{8}(\Gamma_4)$ and the image of $\Delta$
is a cusp form that vanishes on the closure of the Jacobian locus
(i.e., the abelian $4$-folds that are Jacobians of curves of genus $4$),
\cite{B-K}.
Or take $k=g=6$ and get a lift $S_{12}(\Gamma_1) \to S_{12}(\Gamma_{12})$.
This lifted form occurs in the paper \cite{B-F-W}.
\end{example}

Miyawaki observed in \cite{Miyawaki} that the standard $L$-function
of a non-zero cusp form $F$ of weight $12$ on $\Gamma_3$ is a product
$D_{\Delta}(F,s) L(\phi_{20},s+10)L(\phi_{20},s+9)$, with $\Delta\in
S_{12}(\Gamma_1)$ and $\phi_{20} \in S_{20}(\Gamma_1)$ the normalized
Hecke eigenforms of weight $12$ and $20$. He conjectured a lifting and his
idea was refined by Ikeda to the following conjecture.

\begin{conjecture} {\rm (Miyawaki-Ikeda)} 
\index{Miyawaki-Ikeda conjecture}
Let $k$ and $n$ be natural numbers 
with $k-n$ even. Furthermore, let $f\in S_{2k}(\Gamma_1)$ be
a normalized Hecke eigenform and $F_{2n} \in S_{k+n}(\Gamma_{2n})$ 
the Ikeda lift of $f$. Then there exists for every eigenform $g \in
S_{k+n+r}(\Gamma_r)$ with $n,r \geq 1$ a Siegel modular eigenform
${\mathcal F}_{f,g} \in S_{k+n+r}(\Gamma_{2n+r})$ such that
$$
D_{{\mathcal F}_{f,g}}(s)=Z_g(s) \, \prod_{j=1}^{2n} L_f(s+k+n-j),
$$
with $L_f=Z_f$ the usual $L$-function.
\end{conjecture}
In \cite{Ikeda2} Ikeda constructs a lifting from Siegel modular cusp
forms of degree $r$ to Siegel cusp forms of degree $r+2n$. This is a partial
confirmation of this conjecture.

Finally, I would like to mention a conjectured lifting from
vector-valued Siegel modular forms of half-integral weight
to vector-valued Siegel modular forms of integral weight
due to Ibukiyama\index{Ibukiyama lift}. 
He predicts in the case of genus $g=2$
for even $j\geq 0$ and $k\geq 3$ an isomorphism
$$
S_{j,k-1/2}^{+}(\Gamma_0(4), \psi) {\buildrel \sim \over \longrightarrow}
S_{2k-6,j+3}(\Gamma_2)
$$
which should generalize the Shimura-Kohnen lifting 
\index{Shimura-Kohnen lifting} $S_{k-1/2}^{+}(\Gamma_0(4))
\cong S_{2k-2}(\Gamma_1)$, see \cite{Ibu4}. 
Here $\psi(\gamma)=\left( {-4 \over \det(d)} \right)$.
\end{section}
\begin{section}{The Moduli Space of Principally Polarized Abelian
Varieties}
It is a fundamental fact, due to Mumford, that the moduli space of 
principally polarized abelian varieties exists as an algebraic stack
${\mathcal A}_g$ over the integers. The orbifold $\Gamma_g \backslash
{\mathcal H}_g$ is the complex fibre ${\mathcal A}_g(\CC)$ of this
algebraic stack. This fact has very deep consequences for the arithmetic
theory of Siegel modular forms, but an exposition of this
exceeds the framework of these lectures.
Also the various compactifications, 
the Baily-Borel \index{Baily-Borel compactification}
or Satake compactification \index{Satake compactification}
and the toroidal compactifications constructed by Igusa and Mumford et.\ al.\
exist over $\ZZ$ as was shown by Faltings. We refer to an extensive, 
but very condensed survey of this theory in \cite{F-C}.
In particular, Faltings constructed the Satake compactification
over $\ZZ$ as the image of a toroidal compactification 
$\tilde{\mathcal A}_g$ by the sections of a sufficiently big
power of $\det(\EE)$, the determinant of the Hodge bundle.
A corollary of Faltings' results is that
the ring of classical Siegel modular forms with integral
Fourier coefficients is finitely generated over $\ZZ$.

In the following sections we shall sketch how one can use 
some of these facts to extract information on the Hecke eigenvalues 
of Siegel modular forms.

The action of the Galois group of $\QQ$ on the points of
${\mathcal A}_g(\bar{\QQ})$ that correspond to abelian
varieties with complex multiplication is described in
Shimura's theory of canonical models\index{canonical model}. 
This theory can also
explain the integrality of the eigenvalues of Hecke operators.
For this we refer to two papers by Deligne, see \cite{De2,De3}.
\end{section}
\begin{section}{Elliptic Curves over Finite Fields}
Suppose we did not have the elementary approach to
$g=1$ modular forms using holomorphic functions on the upper half plane
like the Eisenstein series and $\Delta$. How would we get the arithmetic
information hidden in the Fourier coefficients of Hecke eigenforms?
Would we encounter $\Delta$?

We claim that one would by playing with elliptic curves over 
finite fields.
Let $\FF_q$ with $q=p^m$ be a finite field of characteristic $p$ and
cardinality $q$. An elliptic curve $E$ defined over $\FF_q$ can be given as
an affine curve by an equation
$$
y^2+a_1xy+a_3y=x^3+a_2x^2+a_4x+a_6,
$$
with $a_i \in \FF_q$ and with non-zero discriminant (a polynomial in the coefficients). 
We can then
count the number $\# E(\FF_q)$ of $\FF_q$-rational points of $E$. A result
of Hasse tells us that $\# E(\FF_q)$ is of the form $q+1-\alpha-\bar{\alpha}$
for some algebraic integer $\alpha$ with $|\alpha|=\sqrt{q}$. We can do this
for all elliptic curves $E$ defined over $\FF_q$ up to $\FF_q$-isomorphism
and we could ask (as Birch did in \cite{Bi}) 
for the average of $\# E(\FF_q)$, or better for
$$
\sum_E \frac{q+1-\# E(\FF_q)}{\# {\rm Aut}_{\FF_q}(E)},
$$
where ${\rm Aut}_{\FF_q}(E)$ is the group of $\FF_q$-automorphisms of $E$,
or more generally we could ask for the average of the expression
$$
h(k,E):= \alpha^k+\alpha^{k-1}\bar{\alpha}+\ldots + \alpha \bar{\alpha}^{k-1}
+\bar{\alpha}^k,
$$
i.e.\ we sum
$$
\sigma_k(q)= - \sum_E \frac{h(k,E)}{\#{\rm Aut}_{\FF_q}(E)}
$$
where the sum is over all elliptic curves $E$ defined over $\FF_q$ up to
$\FF_q$-isomorphism. (As a rule of thumb, whenever one counts mathematical
objects one should count them with weight $1/\# {\rm Aut}$ with ${\rm Aut}$
the group of automorphisms of the object.) If we do this for $\FF_3$
we get the following table, where we also give the $j$-invariant
of the curve $y^2=f$

\smallskip
\vbox{
\bigskip\centerline{\def\quad{\hskip 0.6em\relax}
\def\quod{\hskip 0.5em\relax }
\vbox{\offinterlineskip
\hrule
\halign{&\vrule#&\strut\quod\hfil#\quad\cr
height2pt&\omit&&\omit&&\omit&&\omit&\cr
&$f$ && $\# E(k)$ && $1/\# {\rm Aut}_k(E)$ && j &\cr
\noalign{\hrule}
&$x^3+x^2+1$ && $6$ && $1/2$ && $-1$ & \cr
&$x^3+x^2-1$ && $3$ && $1/2$ && $1$  &\cr
&$x^3-x^2+1$ && $5$ && $1/2$ && $1$  & \cr
&$x^3-x^2-1$ && $2$ && $1/2$ && $-1$ & \cr
&$x^3+x    $ && $4$ && $1/2$ && $0$  & \cr
&$x^3-x    $ && $4$ && $1/6$ && $0$  & \cr
&$x^3-x+1  $ && $7$ && $1/6$ && $0$  & \cr
&$x^3-x-1  $ && $1$ && $1/6$ && $0$  & \cr
} \hrule}
}}
\bigskip

\noindent
and obtain 
the following frequencies for the number of $\FF_3$-rational
points:

\smallskip
\vbox{
\bigskip\centerline{\def\quad{\hskip 0.6em\relax}
\def\quod{\hskip 0.5em\relax }
\vbox{\offinterlineskip
\hrule
\halign{&\vrule#&\strut\quod\hfil#\quad\cr
height2pt&\omit&&\omit&&\omit&&\omit&&\omit&&\omit&&\omit&&\omit&\cr
\noalign{\hrule}
& n && 1 && 2 && 3 && 4&& 5&& 6 && 7 &\cr
\noalign{\hrule}
&freq && 1/6 && 1/2 && 1/2 && 2/3 && 1/2 && 1/2 && 1/6 &\cr
} \hrule}
}}
\bigskip

Note that $\sum 1/{\rm Aut}_{\FF_q}(E) =q$ and  
$\sum_{E \colon j(E)=j}1/{\rm Aut}_{\FF_q}(E) =1$ (see
\cite{vdG-vdV} for a proof); so a `physical point' of the moduli space
contributes $1$.

If we work this out not only for $p=3$, but for several primes
($p=2,3,5,7$ and $11$) we get the following values:

\smallskip
\vbox{
\bigskip\centerline{\def\quad{\hskip 0.6em\relax}
\def\quod{\hskip 0.5em\relax }
\vbox{\offinterlineskip
\hrule
\halign{&\vrule#&\strut\quod\hfil#\quad\cr
height2pt&\omit&&\omit&&\omit&&\omit&&\omit&&\omit&\cr
\noalign{\hrule}
& $p$ && $2$ && $3$ && $5$ && $7$ && $11$ &\cr
\noalign{\hrule}
& $\sigma_{10}$ && $-23$ && $253$ && $4831$ && $-16743$ && $534613$ &\cr
} \hrule}
}}
\bigskip

Anyone who remembers the cusp form $\Delta=\sum_{n>0} \tau(n) q^n=
q-24\, q^2 +252 \, q^3 -3520 \, q^4 + 4830 \, q^5 + \ldots $ 
will not fail to notice
that $\sigma_{10}(p)=\tau(p)+1$ for the primes listed in this example. 
And in fact, the relation $\sigma_{10}(p)=\tau(p)+1$ 
holds for all primes $p$.
The reason behind this is that the cohomology of the $n$th power of 
the universal elliptic curve ${\mathcal E} \to {\mathcal A}_1$
is expressed in terms of cusp forms on ${\rm SL}(2,\ZZ)$.
To describe this we recall the local system $\WW$ on ${\mathcal A}_1$
associated to $\eta^{-1}$ times the standard representation of
${\rm GSp}(2,\QQ)$ in Section \ref{roots}. The fibre of this local
system over a point $[E]$ given by the elliptic curve $E$
can be identified with the cohomology group $H^1(E,\QQ)$.
Or consider the universal elliptic curve 
(in the orbifold sense)
$\pi :{\Ecal} \to {\Acal}_1$ obtained as the quotient
${\rm SL}(2,\ZZ) \times \ZZ^2 \backslash \H1 \times \CC$, 
where the action of
$(a,b;c,d) \in {\rm SL}(2,\CC)$ on $(\tau,z) \in \H1 \times \CC$ 
is $((a\tau+b)/(c\tau+d), (c\tau+d)^{-1} z)$. Associating to an
elliptic curve its homology $H_1(E,\QQ)$ defines a local system
that can be obtained as a quotient 
${\rm SL}(2,\ZZ) \backslash {\mathcal H}_1 \times \QQ^2$.
Then the dual of this
local system is $\WW := R^1\pi_*\QQ$. We now put
$$
\WW^k := {\rm Sym}^k(\WW),
$$
a local system with a $k+1$-dimensional fibre for $k\geq 0$. 
We now have the following cohomological interpretation of cusp 
forms on
${\rm SL}(2,\ZZ)$, cf.\ \cite{De1}.

\begin{theorem} {\rm (Eichler-Shimura)} \index{Eichler-Shimura}
For even $k\in \ZZ_{\geq 2}$ we have an isomorphism of the
compactly supported cohomology of $\WW^k$
$$
H_c^1({\Acal}_1, \WW^k \otimes \CC) \cong 
S_{k+2}\oplus \bar{S}_{k+2} \oplus \CC
$$
with $S_{k+2}$ the space of cusp forms of weight $k+2$ on ${\rm SL}(2,\ZZ)$
and $\bar{S}_{k+2}$ the complex conjugate of this space.
\end{theorem}
Replacing $\WW$ by $\WW_{\RR}$ we have the exact sequence
$$
0 \to \EE \to \WW \otimes_{\RR} O \to \EE^{\vee} \to 0
$$
with $O$ the structure sheaf and an 
induced map $\EE^{\otimes k} \to \WW^k \otimes_{\RR} O$. Now the
de Rham resolution
$$
0 \to \WW^k\otimes_{\RR} \CC  \to \WW^k \otimes_{\RR} O {\buildrel d \over
\longrightarrow} \WW^k \otimes \Omega^1 \to 0
$$
defines a connecting homomorphism
$$
H^0({\Acal}_1, \Omega^1(\WW^k)) \to H^1({\Acal}_1, \WW^k \otimes \CC).
$$
The right hand space has a natural complex conjugation and we 
thus find also a complex conjugate map 
$$
\overline{H^0({\Acal}_1, \Omega^1(\WW^k))} 
\to H^1({\Acal}_1, \WW^k \otimes \CC).
$$
A cusp form $f \in S_{k+2}$ defines a section of 
$H^0({\Acal}_1, \Omega^1(\WW^k))$ by putting
$f(\tau) \mapsto f(\tau) d\tau dz^k$. We thus have a cohomological 
interpretation of the space of cusp forms.

As observed above
the moduli space ${\Acal}_1$ is defined over the integers $\ZZ$. This means that we also have the moduli space ${\Acal}_1 \otimes \FF_p$ of elliptic curves
in characteristic $p>0$. It is well-known that one can obtain a lot of
information about cohomology by counting points over finite fields. 
(Here we work with $\ell$-adic \'etale cohomology for $\ell \neq p$.)
And, indeed, there exists an analogue of the Eichler-Shimura isomorphism 
in characteristic $p$
and the relation $\sigma_{10}(p)=\tau(p)+1$ 
is a manifestation of this.
In fact a good notation for writing this relation is
$$
H_c^1({\Acal}_1,\WW^{10}) = S[12]+1,
$$
where the formula
$$
H_c^1({\Acal}_1,\WW^{2k})\cong S[2k+2] +1 \quad \hbox{\rm for $k\geq 1$}
$$
may be interpreted complex-analytically as the Eichler-Shimura isomorphism
and in characteristic $p$ as the relation
$$
\sigma_{2k}(p)=1+\hbox{\rm Trace of $T(p)$ on $S_{2k+2}$}.
$$
(A better interpretation is as a relation in a suitable $K$-group
and with $S[2k+2]$ as the motive \index{motive} associated to $S_{2k+2}$. 
This motive 
can be constructed in the $k$th power of ${\mathcal E}$ 
as done by Scholl \cite{Scholl} or using moduli space of $n$-pointed
elliptic curves as done by Consani and Faber, \cite{Consani-Faber}.)

This $1$ in the formula 
$H_c^1({\Acal}_1,\WW^{2k})\cong S[2k+2] +1$
is really a nuisance. To get rid of it we 
consider
the natural map
$$
H_c^1({\Acal}_1,\WW^k) \to H^1({\Acal}_1,\WW^k)
$$
the image of which is called the \Definition{interior cohomology} 
\index{interior cohomology} 
and denoted by $H_{!}^1({\Acal}_1,\WW^k)$. 
We thus have an elegant and sophisticated form of the 
Eichler-Shimura isomorphism
$$
H^1_c({\Acal}_1,\WW^k) = S[k+2]+1, \quad H_{!}^1({\Acal}_1,\WW^k)=S[k+2].
$$
The $1$ is the $1$ in $1+p^{k+1}$, the eigenvalue of the action of 
$T(p)$  on the
Eisenstein series $E_{k+2}$ of weight $k+2$ on ${\rm SL}(2,\ZZ)$.

The moral of this is that we can obtain information on the traces of
Hecke operators on the space $S_{k+2}$ by calculating $\sigma_{k}(p)$,
i.e., by counting points on elliptic curves over $\FF_p$. Even from a purely
computational point of view this is not a bad approach to calculating
the traces of Hecke operators.

\end{section}

\begin{section}{Counting Points on Curves of Genus $2$}\label{genus2}
With the example of $g=1$ in mind it is natural to ask whether also for
$g=2$ we could obtain information on modular forms using curves of genus $2$
over finite fields. In joint work with Carel Faber (\cite{F-G})
we showed that we can.

For $g=2$ the quotient space $\Gamma_2 \backslash \H2$ is the analytic space
of the moduli space ${\Acal}_2$ of principally polarized abelian surfaces.
A principally polarized abelian surface is the Jacobian of a smooth 
projective irreducible algebraic curve or it is a product of two
elliptic curves. If the characteristic is not $2$ a curve of genus $2$
can be given as an affine curve with equation $y^2=f(x)$ with $f$ a
polynomial of degree $5$ or $6$ without multiple zeros.

The moduli space ${\Acal}_2$ exists over $\ZZ$ and provides us with a
moduli space ${\Acal}_2 \otimes \FF_p$ for every characteristic $p>0$.
Also here we have a local system which is the analogue of the local system
$\WW$ that we saw for $g=1$:
$$
\VV:= {\rm GSp}(4,\ZZ) \backslash \H2 \times \QQ^4,
$$
where the action of $\gamma=(a,b;c,d) \in
{\rm GSp}(4,\ZZ)$ is given by $\eta^{-1}$ times the standard representation.
Or in more functorial terms, we consider the universal
family $\pi: {\Xcal}_2 \to {\Acal}_2$ and then $\VV$ is the direct image
$R^1\pi_*(\QQ)$. The fibre of this local system over the point $[X]$
corresponding to the polarized abelian surface $X$ is $H^1(X,\QQ)$.
The local system \index{local system}
$\VV$ comes equipped with a symplectic pairing
$\VV \times \VV \to \QQ(-1)$. Just as for $g=1$ we made the local systems 
$\WW^k$ out of the basic one $\WW$ we can construct more
local systems out of $\VV$ but now parametrized by two indices $l$ and $m$
with $l \geq m \geq 0$. Namely, the irreducible representations of
${\rm Sp}(4,\QQ)$ are parametrized by such pairs $(l,m)$ and we thus
have local systems $\VV_{l,m}$ with $l \geq m \geq 0$ such that
$\VV_{l,0}={\rm Sym}^{l}(\VV)$ and $\VV_{1,1}$ is the `primitive part'
of $\wedge^2 \VV$. A local system $\VV_{l,m}$ is called \Definition{regular}
if $l > m >0$. \index{regular local system}

Just as in the case $g=1$ we are now interested in the cohomology of the 
local systems $\VV_{l,m}$. We put
$$
e_c({\Acal}_2,\VV_{l,m})= \sum_i (-1)^i [ H_c^i({\Acal}_2,\VV_{l,m})].
$$
Here we consider the alternating sum of the cohomology groups with
compact support in the Grothendieck group of mixed Hodge structures.

We also have an $\ell$-adic analogue of this that can be used in 
positive characteristic. It is obtained from $R^1\pi_*(\QQ_{\ell})$
and lives over ${\Acal}_2 \otimes \ZZ[1/\ell]$; we consider the \'etale
cohomology of this sheaf. We simply use the same name $\VV_{l,m}$
and assume that $\ell$ is different from the characteristic $p$.

Using a theorem of Getzler (\cite{Getzler}
(on ${\Mcal}_2$)
tells us what the Euler characteristic 
$\sum_i (-1)^i \dim H_c^i({\Acal}_2,\VV_{l,m})$ over $\CC$ is.
This Euler characteristic equals the Euler characteristic of the
$\ell$-adic variant over a finite field.

The first observation is that because of the action of the 
hyperelliptic
involution these cohomology groups are zero for $l+m$ odd.

Our strategy is now to make a list of all $\FF_q$-isomorphism classes
of curves of genus $2$ over $\FF_q$ and to determine for each of them
$\# {\rm Aut}_{\FF_q}(C)$ and the characteristic polynomial of Frobenius.
So for each curve $C$ we determine algebraic integers $\alpha_1,
\bar{\alpha}_1, \alpha_2, \bar{\alpha}_2$ of absolute value $\sqrt{q}$
such that
$$
\# C(\FF_{q^i})= q^i+1-\alpha_1^i-\bar{\alpha}_1^i -\alpha_2^i-\bar{\alpha}_2^i 
$$
for all $i\geq 1$. These $\alpha$'s can be calculated using this identity for $i=1$ and $i=2$.
We also must calculate the contribution from the
degenerate curves of genus $2$, i.e., the contribution from the principally
polarized abelian surfaces that are products of elliptic curves.

Having done that we are able to calculate the trace of Frobenius
on the alternating sum of $H^i_c({\Acal}_2 \otimes \FF_q , \VV_{l,m})$ , 
where by $\VV_{l,m}$
we mean the $\ell$-adic variant, a smooth $\ell$-adic sheaf
on ${\Acal}_2 \otimes \FF_q$. 
In practice, it means that we sum a certain symmetric expression
in the $\alpha$'s divided by $\# {\rm Aut}_{\FF_q}(C)$, analoguous
to the $\sigma_k(q)$ for genus $1$.

What does this tell us about Siegel modular forms of degree $g=2$? To get
the connection with modular forms we have to replace the compactly
supported cohomology by the interior cohomology, i.e., by the image of
$H_c^i({\Acal}_2,\VV_{l,m}) \to H^i({\Acal}_2,\VV_{l,m})$ which is
denoted by $H_!^i({\Acal}_2,\VV_{l,m})$. So let us define
$$
e_{\rm Eis} ({\Acal}_2,\VV_{l,m})= e_c({\Acal}_2,\VV_{l,m})-e_!({\Acal}_2,
\VV_{l,m}).
$$
If we do the same thing for $g=1$ we find $e_{\rm Eis}({\Acal}_1,\WW^k)=
-1$ for even $k>0$.

Let $\LL$ be the $1$-dimensional Tate Hodge structure of weight $2$
\index{Tate Hodge structure}. It corresponds to the second cohomology of $\PP^1$. In terms of counting points one reads $q$ for $\LL$.
Our first result is (cf., \cite{F-G})

\begin{theorem} Let $(l,m)$ be regular. Then $e_{\rm Eis}({\Acal}_2,\VV_{l,m})$
is given by
$$
-S[l+3]-s_{l+m+4}\LL^{m+1}+S[m+2]+s_{l-m+2}\cdot 1 +
\begin{cases} 1 & \hbox{\rm $l$ even} \\
0 & \hbox{\rm $l$ odd,}\\ \end{cases}
$$
where $s_{n}=\dim S_n(\Gamma_1)$.
\end{theorem}

Faltings has shown (see \cite{F-C}) that $H_!^3({\Acal}_2,\VV_{l,m})$
possesses a Hodge filtration
$$
0 \subset F^{l+m+3} \subset F^{l+2} \subset F^{m+1} \subset F^0=
H_!^3({\Acal}_2,\VV_{l,m}).
$$
Moreover, if $(l,m)$ is regular then $H_{!}^i({\Acal}_2,\VV_{l,m})=(0)$
for $i\neq 3$. Furthermore, Faltings shows that
$$
F^{l+m+3}\cong S_{l-m,m+3}(\Gamma_2).
$$
Here $S_{j,k}(\Gamma_2)$ is the space of Siegel modular forms 
for the representation ${\rm Sym}^j \otimes \det^k$ of ${\rm GL}(2,\CC)$.
This is the sought-for connection with vector valued Siegel modular forms
and the analogue of $H^1_{!}({\Acal}_1,\WW^k)=F^0
\supset F^{k+1} \cong S_{k+2}(\Gamma_1)$ for $g=1$.
Faltings gives an interpretation of all the steps in the Hodge filtration
in terms of the cohomology of the bundles ${\mathcal W}(\lambda)$.

However, although for $g=1$ the Eichler-Shimura isomorphism tells us that we
know $H^1_{!}({\Acal}_1,\WW^k)$ once we know $S_{k+2}(\Gamma_1)$,
for $g=2$ 
there might be pieces of cohomology hiding in $F^{l+2} \subset F^{m+1}$
that are not detectable in $F^{l+m+3}$ or in $F^0/F^{m+1}$ and indeed there
is such cohomology. The contribution to this part of the cohomology is called
the contribution from \Definition{endoscopic lifting from $N={\rm GL}(2)
\times {\rm GL}(2) / \GG_m$}. \index{endoscopic contribution}

We conjecture on the basis of our numerical calculations that this endoscopic
contribution is as follows.

\begin{conjecture} Let $(l,m)$ be regular. Then the endoscopic contribution
is given by
$$
e_{\rm endo}({\Acal}_2,\VV_{l,m})= -s_{l+m+4} S[l-m+2]\, \LL^{m+1}.
$$
\end{conjecture}

There is a very extensive literature on endoscopic lifting (cf.\ \cite{KR}), 
but a precise
result on the image in our case seems to be absent. Experts on 
endoscopic lifting should be able to prove this conjecture.
Actually, since we know the Euler characteristics of the
interior cohomology and have Tsushima's dimension formula
it suffices to construct a subspace of dimension 
$2s_{l+m+4}s_{l-m+2}$ in the endoscopic part
via endoscopic lifting for regular $(l,m)$.

In terms of Galois representations a Siegel modular form 
(with rational Fourier coefficients) should correspond to a 
rank $4$ part of the cohomology or
a $4$-dimensional irreducible Galois representation. 
A modular form in the 
endoscopic part corresponds to a rank $2$ part and a $2$-dimensional
Galois representation. Modular forms coming from the Saito-Kurokawa
lift give $4$-dimensional representations that split off two $1$-dimensional pieces.

In analogy with the case of $g=1$ we now set
$$
S[l-m,m+3] := H^3_{!}({\Acal}_2,\VV_{l,m})-H^3_{\rm endo}({\Acal}_2,\VV_{l,m}).
$$
This should be a motive analogous to the motive $S[k]$ we encountered for $g=1$
and lives in a  power of the universal abelian surface 
over ${\mathcal A}_2$.
The trace of Frobenius on \'etale $\ell$-adic
$H^3_{!}({\Acal}_2,\VV_{l,m})-H^3_{\rm endo}({\Acal}_2,\VV_{l,m})$ 
should be the trace of the Hecke operator $T(p)$ on the 
space of modular forms $S_{l-m,m+3}$.

\end{section}

\begin{section}{The Ring of Vector-Valued Siegel Modular Forms 
for Genus $2$}\label{vvgenus2}
The quest for vector-valued Siegel modular forms starts with genus $2$.
We can consider the 
direct sum $M=\oplus_{\rho} M_{\rho}(\Gamma_2)$ (see Section 
\ref{modularforms}), \index{ring of modular forms}
where $\rho$ runs through
the set of irreducible polynomial representations of ${\rm GL}(2,\CC)$. 
Each such
$\rho$ is given by a pair $(j,k)$ such that $\rho = {\rm Sym}^j(W) \otimes
\det (W)^{k}$, with $W$ the standard representation of ${\rm GL}(2,\CC)$.
(Note that in the earlier notation we have $(\lambda_1-\lambda_2,\lambda_2)=(j,k)$.)
So we may write $M=\oplus_{j,k \geq 0} M_{j,k}(\Gamma_2)$ and we know that 
$M_{j,k}(\Gamma_2)=(0)$ if $j$ is odd. If $F$ and $F'$ are Siegel modular
forms of weights $(j,k)$ and $(j',k')$ then the product 
is a modular forms of weight of weight $(j+j',k+k')$. 
The multiplication is obtained from the canonical map
${\rm Sym}^{j_1}(W)\otimes \det(W)^{k_1} \otimes
{\rm Sym}^{j_2}(W)\otimes \det(W)^{k_2} \to 
{\rm Sym}^{j_1+j_2}(W)\otimes \det(W)^{k_1+k_2}$
obtained from multiplying polynomials in two variables.

There is the Siegel operator that goes
from $M_{j,k}(\Gamma_2)$ to $M_{j+k}(\Gamma_1)$. 
For $j>0$ the Siegel operator
gives a map to $S_{j+k}(\Gamma_1)$ and for $j>0, k>4$ the map
$\Phi: M_{j,k}(\Gamma_2) \to S_{j+k}(\Gamma_1)$ is surjective. 
For these facts on the Siegel operator
we refer to Arakawa's paper \cite{Ar}.
The Siegel operator is multiplicative: $\Phi(F\cdot F')=\Phi(F)\, \Phi(F')$.

There is a dimension formula
for $\dim M_{j,k}(\Gamma_2)$, due to Tsushima, \cite{Tsus2}. But apart
from this not much is known about  vector-valued Siegel modular forms.
The direct sum $\oplus_{k} M_{j,k}(\Gamma_2)$
for fixed $j$ is a module over the ring $M^{\rm cl}=\oplus M_{0,k}(\Gamma_2)$
of classical Siegel modular forms and we know generators of this module for
$j=2$ and $j=4$ and even $j=6$ due to Satoh and Ibukiyama, cf.\
\cite{Satoh,Ibu1,Ibu2}. 

One way to construct vector-valued Siegel modular forms
from classical Siegel modular forms is differentiation,
the simplest example being given by a pair $f\in M_a(\Gamma_2)$,
$g \in M_b(\Gamma_2)$ for which one sets
$$
[f,g]:=\frac{1}{b} \, f \nabla g -\frac{1}{a} g \nabla f
$$
with $\nabla f$ defined by
$$
2\pi i \nabla f = a (2i y)^{-1} f + \left(
\begin{matrix} \partial / \partial \tau_{11} & \partial / \partial \tau_{12} \\
\partial / \partial \tau_{12} & \partial / \partial \tau_{22} 
\end{matrix} \right) f.
$$
The point is that $[f,g]$ is then a modular form in $M_{2,a+b}(\Gamma_2)$.
Using this operation (an instance of Cohen-Rankin operators)
Satoh showed in \cite{Satoh} that $\oplus_{k\equiv 0(2)} 
M_{2,k}$ is generated over the ring $\oplus_{k} M_k(\Gamma_2)$
of classical Siegel modular forms by such $[f,g]$ with $f$ and $g$
classical Siegel modular forms.

We give a little table with dimensions for
$\dim S_{j,k}(\Gamma_2)$ for $4\leq k \leq 20$, $0 \leq j \leq 18$ with $j$
even:
$$
 \left[ \begin {array}{ccccccccccccccccccc}j\backslash k &&4&5&6&7&8&9&10&11&12&13&14&15&16&17&18&19&20\\
\\
0&&0&0&0&0&0&0&1&0&1&0&1&0&2&0&2&0&3\\
\noalign{\medskip}2&&0&0&0&0&0&0&0&0&0&0&1&0&2&0&2&0&3
\\\noalign{\medskip}4&&0&0&0&0&0&0&1&0&1&0&2&1&3&1&4&2&6
\\\noalign{\medskip}6&&0&0&0&0&1&0&1&1&2&1&3&2&5&3&7&4&9
\\\noalign{\medskip}8&&0&0&0&0&1&1&2&1&3&2&5&4&7&5&9&7&13
\\\noalign{\medskip}10&&0&0&0&0&0&1&2&1&3&2&5&5&8&6&11&9&15
\\\noalign{\medskip}12&&0&0&1&1&2&2&4&4&6&5&9&8&13&11&17&15&22
\\\noalign{\medskip}14&&0&0&0&1&2&2&4&4&6&6&10&10&15&13&19&18&26
\\\noalign{\medskip}16&&0&0&1&1&3&3&6&5&9&8&13&13&19&17&25&23&33
\\\noalign{\medskip}18&&0&1&1&2&4&5&7&8&11&11&17&17&23&23&31&30&40
\end {array} \right] 
$$

The ring $\oplus_{j,k} M_{j,k}(\Gamma_2)$ is not finitely generated
as was explained to me by Christian Grundh. Here is his argument.

\begin{lemma} The ring
$\oplus_{j,k} M_{j,k}(\Gamma_2)$ is not finitely generated.
\end{lemma}
\begin{proof}
Suppose that $g_n$ for $n=1,\ldots,r$ are the generators with 
weights $(j_n,k_n)$. If we have a modular form $g$ of weight $(j,k)$
with $j> \max(j_n, n=1,\ldots,r)$ then $g$ is a sum of products
of $g_n$, two of which at least 
have $j_n >0$, hence by the properties of $\Phi$ we see
that then $\Phi(g)$ is a sum of products of cusp forms, hence lies
in the ideal generated by $\Delta^2$ of the ring of elliptic modular
forms. But for $j>0, k>4$ the map $\Phi: M_{j,k}(\Gamma_2)
\to S_{j+k}(\Gamma_1)$ is surjective, so we have forms $g$ in 
$M_{j,k}(\Gamma_2)$ that land in the ideal generated by $\Delta$,
but not in the ideal generated by $\Delta^2$.
Thus the ring cannot be generated by $g_n$ for $n=1,\ldots,r$.
\end{proof}

Just as $\Delta$ is the first cusp form for $g=1$ that one encounters
the first vector-valued cusp form that one encounters for $g=2$
is the generator of $S_{6,8}(\Gamma_2)$.  The adjective `first' 
refers to the fact that the weight of the local
system $\VV_{j+k-3,k-3}$ is $j+2k-6$.
Our calculations (modulo the endoscopic conjecture)
allow the determination of the eigenvalues $\lambda(p)$
and $\lambda(p^2)$ for $p=2,3,5,7$. We then can calculate the
characteristic polynomial of Frobenius and the slopes of it
on $S_{6,8}(\Gamma_2)$. \index{slope}

\smallskip
\vbox{
\bigskip\centerline{\def\quad{\hskip 0.6em\relax}
\def\quod{\hskip 0.5em\relax }
\vbox{\offinterlineskip
\hrule
\halign{&\vrule#&\strut\quod\hfil#\quad\cr
height2pt&\omit&&\omit&&\omit&&\omit&\cr
&$p$&&$\lambda(p)$&&$\lambda(p^2)$&&{\rm slopes}&\cr
height2pt&\omit&&\omit&&\omit&&\omit&\cr
\noalign{\hrule}
height2pt&\omit&&\omit&&\omit&&\omit&\cr
&$2$&&$0$&&$-57344$&&$13/2,25/2$&\cr
&$3$&&$-27000$&&$143765361$&&$3,7,12,16$&\cr
&$5$&&$2843100$&&$-7734928874375$&&$2,7,12,17$&\cr
&$7$&&$-107822000$&&$4057621173384801$&&$0,6,13,19$&\cr
height2pt&\omit&&\omit&&\omit&&\omit&\cr
} \hrule}
}}
At our request Ibukiyama (\cite{Ibu1}) has constructed a
vector-valued Siegel modular form 
$0\neq F \in S_{6,8}\,$, using a theta series for the lattice
$\Gamma =\{ x \in \QQ^{16}\, : \, 2x_i\in \ZZ, x_i-x_j \in \ZZ,
\sum_{i=1}^{16} x_i \in 2\ZZ\} $.
One puts $a=(2,i,i,i,i,0,\ldots,0) \in \CC^{16}$ and one denotes by
$(\, , \, )$
the usual scalar product. If $F=(F_0,\ldots,F_6)$ is the vector of
functions on ${\Hcal}_2$ defined by
$$
F_{\nu}=\sum_{x,y \in \Gamma} (x,a)^{6-\nu} (y,a)^{\nu} e^{\pi i ((
x,x)\tau_{11}+2(x,y)\tau_{12}+(y,y)\tau_{22})} \qquad (\nu=0,\ldots,6)
$$
with $\tau= (\tau_{11}, \tau_{12}; \tau_{12}, \tau_{22}) \in {\Hcal}_2$,
then Ibukiyama's result is that $F\neq 0$ and $F \in S_{6,8}$.
The vanishing of $\lambda(2)$ agrees with this.

Here are two more examples of $1$-dimensional spaces, 
the space $S_{18,5}$ and the last one, $S_{28,4}$. 
In these examples and the other ones we assume the validity of
our conjecture on the endoscopic contribution.
The eigenvalues $\lambda(p)$
grow approximately like $p^{(j+2k-3)/2}$, 
i.e. $p^{25/2}$ and $p^{33/2}$.

\smallskip
\vbox{
\bigskip\centerline{\def\quad{\hskip 0.6em\relax}
\def\quod{\hskip 0.5em\relax }
\vbox{\offinterlineskip
\hrule
\halign{&\vrule#&\strut\quod\hfil#\quad\cr
height2pt&\omit&&\omit&&\omit&\cr
&$p$&& $\lambda(p)$  on $S_{18,5}$&&$\lambda(p)$ on $S_{28,4}$&\cr
height2pt&\omit&&\omit&&\omit&\cr
\noalign{\hrule}
height2pt&\omit&&\omit&&\omit&\cr
&$2$&&$-2880 $&&$35040$&\cr
&$3$&&$-538920$&&$30776760$&\cr
&$5$&&$118939500$&&$522308049900$&\cr
&$7$&&$1043249200$&&$18814963644400$&\cr
&$11$&&$-9077287359096$&&$132158356344353064$&\cr
&$13$&&$-133873858788740$&&$-1710588414695522180$&\cr
&$17$&&$667196591802660$&&$-17044541241181641180$&\cr
&$19$&&$2075242468196920$&&$888213094972004807320$&\cr
&$23$&&$-8558834216776560$&&$-43342643806617018857520$&\cr
&$29$&&$64653981488634780$&&$-172663192093972503614820$&\cr
&$31$&&$-5977672283905752896$&&$1826186223285615270299584$&\cr
&$37$&&$56922208975445092780$&&$-29747516862655204839491540$&\cr
height2pt&\omit&&\omit&&\omit&\cr
} \hrule}
}}
\bigskip
In principle our database allows for the calculation of the traces
of the Hecke operators $T(p)$ with $p\leq 37$ 
on the spaces $S_{j,k}$ {\sl for all values} 
$j,k$. In the cases at hand these numbers tend to be `smooth', i.e., they
are highly composite numbers as we illustrate with 
the two $1$-dimensional spaces $S_{j,k}$
for $(j,k)=(8,8)$ and $(12,6)$ (where the trace equals the eigenvalue
of $T(p)$).

\smallskip
\vbox{
\bigskip\centerline{\def\quad{\hskip 0.6em\relax}
\def\quod{\hskip 0.5em\relax }
\vbox{\offinterlineskip
\hrule
\halign{&\vrule#&\strut\quod\hfil#\quad\cr
height2pt&\omit&&\omit&&\omit&\cr
&$p$&& $\lambda(p)$  on $S_{8,8}$&&$\lambda(p)$ on $S_{12,6}$&\cr
height2pt&\omit&&\omit&&\omit&\cr
\noalign{\hrule}
height2pt&\omit&&\omit&&\omit&\cr
&$2$&&$2^6\cdot 3 \cdot 7$&&$-2^4\cdot 3 \cdot 5$&\cr
&$3$&&$-2^3\cdot 3^2 \cdot 89$&&$2^3 \cdot 3^5 \cdot 5 \cdot 7$&\cr
&$5$&&$-2^2\cdot 3 \cdot 5^2 \cdot 13^2 \cdot 607$&&$2^2\cdot 3 \cdot 5^2 \cdot 7 \cdot 79 \cdot 89$&\cr
&$7$&&$2^4\cdot 7 \cdot 109 \cdot 36973$&&
 $-2^4 \cdot 5^2 \cdot 7 \cdot 119633$&\cr
&$11$&&$2^3\cdot 3 \cdot 4759 \cdot 114089$&&
$ 2^3 \cdot 3 \cdot 23 \cdot 2267 \cdot 2861$&\cr
&$13$&&$-2^2\cdot 13 \cdot 17 \cdot 109 \cdot 3404113$&&
$2^2 \cdot 5 \cdot 7 \cdot 13 \cdot 50083049$&\cr
&$17$&&$2^2\cdot 3^2 \cdot 17 \cdot 41 \cdot 1307 \cdot 168331$&&$
-2^2\cdot 3^2 \cdot 5 \cdot 7 \cdot 13 \cdot 47 \cdot 14320807$&\cr
&$19$&&$-2^3\cdot 5 \cdot 74707 \cdot 9443867$&&
$-2^3 \cdot 5 \cdot 7^3 \cdot 19 \cdot 2377 \cdot 35603$&\cr
height2pt&\omit&&\omit&&\omit&\cr
} \hrule}
}}
\bigskip
Satoh had calculated a few eigenvalues of Hecke operators $T(m)$
acting on $S_{14,2}(\Gamma_2)$, (for $m=2,3,4,5,9$ and $25$)
cf.\ \cite{Satoh}, and our values agree with his.

\end{section}
\begin{section}{Harder's Conjecture}
In his study of the contribution of the boundary of the moduli space
to the cohomology of local systems on the symplectic group, more precisely
of the Eisenstein cohomology, Harder arrived at a conjectural congruence
between modular forms for $g=1$ and Siegel modular forms for $g=2$,
cf., \cite{Harder1,Harder2}. The second reference is his colloquium
talk in Bonn (February 2003) which can be found in this volume.
One can view his conjectured congruences as a generalization of the
famous congruence for the Fourier coefficients of the $g=1$ cusp
form $\Delta =\sum \tau(n)q^n$ of weight $12$
$$
\tau(p)\equiv p^{11} +1 \, (\bmod \, 691).
$$

To formulate it we start with a $g=1$ cusp form $f \in S_r(\Gamma_1)$ of
weight $r$ that is a normalized eigenform of the Hecke operators. We write
$f=\sum_{n\geq 1} a(n) q^n$ with $a(n)=1$. To $f$ we can associate the
$L$-series $L(f,s)$ defined by $L(f,s)= \sum_{n\geq 1} a(n)/n^s$ for
complex $s$ with real part $>k/2+1$. If we define $\Lambda(f,s)$ by
$$
\Lambda(f,s)= \frac{\Gamma(s)}{(2\pi)^s} L(f,s) = \int_{0}^{\infty} f(iy)y^{s-1} dy
$$
then $\Lambda(f,s)$ admits a holomorphic continuation to the whole $s$-plane
and satisfies a functional equation $\Lambda(f,s)=i^k \Lambda(f,k-s)$.
It is customary to call the values $\Lambda(f,t)$ for $t=k-1, k-2, 
\ldots, 0$ the \Definition{critical values}. \index{critical value}
In view of the functional equation we may
restrict to the values $t=k-1,\ldots,k/2$. 

A basic result due to Manin and Vishik is the following.
\index{Manin-Vishik theorem}

\begin{theorem} There exist two real numbers (`periods') $\omega_{+}, 
\omega_{-}$ such that the ratios
$$
\Lambda(f,k-1)/\omega_{-}, \, \Lambda(f,k-2)/\omega_{+},
\ldots, \Lambda(f,k/2)/\omega_{(-1)^{k/2}}
$$
are in the field of Fourier coeffients 
$\QQ_f=\QQ(a(n) : n\in \ZZ_{\geq 1})$.
\end{theorem}

If the Fourier coefficients are rational integers 
we may normalize these ratios 
so that we get integers in a minimal way.
In practice one observes that one usually finds many small primes dividing
these coordinates.
By small we mean here less than $k$ (or something close to this). 
Occasionally, there is a larger prime
 dividing these critical values
of $\Lambda(f,s)$. 

Instead of calculating the integrals one may use a
slightly different approach by employing the so-called period polynomials,
\cite{K-Z}, which are defined for $f \in S_k(\Gamma_1)$
by $r=i\, r^{+}+r^{-}$ with
$$
r^{+}(f)= \sum_{0 \leq n \leq k-2, n\, \hbox{\rm even}} 
(-1)^{n/2} \left( {k-2 \choose n}\right)
r_n(f) X^{k-2-n}
$$
and
$$
r^{-}(f)= \sum_{0 < n<k-2, n \, \hbox{\rm odd}} 
(-1)^{(n-1)/2} \left( {k-2 \choose n}\right)
r_n(f) X^{k-2-n}
$$
with $r_n(f)=\int_{0}^{\infty} f(it)t^n dt$ for $n=0,\ldots,k-2$.
Then the coefficients of these period polynomials give up to `small' primes
the critical $L$-values. These can be calculated purely algebraically
and these are the ones that I used. By slight abuse of notation
I denote these ratios again by the same symbols 
$(\Lambda(f,k-1): \Lambda(f,k-3): \ldots )$. 
See also \cite{Dum} for more on the critical values.

For example, if we do this for $f=\Delta \in S_{12}$ then we get
$$
(\Lambda(f,10): \Lambda(f,8):\Lambda(f,6))=(48:25:20).
$$
and see only `small' primes.
The first example where we see larger primes 
is the normalized eigenform
$f=\Delta e_4e_6 \in S_{22}$. We find
$$
(\Lambda(f,20):\Lambda(f,18): \ldots : \Lambda(f,12))=
(2^5\cdot 3^3 \cdot 5 \cdot 19:2^3\cdot 7 \cdot 13^2:3\cdot 5 \cdot 7 \cdot 13: 2\cdot 3 \cdot 41 : 2 \cdot 3 \cdot 7)
$$
where obviously $41$ is the exception. 
We shall say for short $41 | \Lambda(f,14)$. 
What is the meaning of these exceptional primes dividing the critical values?

Harder made the following conjecture.

\begin{conjecture}
{\rm (Harder's Conjecture)} 
\index{Harder's conjecture}
Let $f\in S_r(\Gamma_1)$ be a normalized eigenform
with field of Fourier coefficients $\QQ_f$.
If a `large' prime $\ell$ of $\QQ_f$ 
divides a critical value $\Lambda(f,t)$ then there
exists a Siegel modular form $F\in S_{j,k}(\Gamma_2)$ 
of genus $2$ and weight $(j,k)$ with $j=2t-r-2$ and
$k=r-t+2$ that is an eigenform for the Hecke algebra with 
eigenvalue $\lambda(p)$ for $T(p)$ with field $\QQ_F$ 
of eigenvalues $\lambda(p)$ 
and such that for a suitable prime $\ell'$ of the compositum of
$\QQ_f$ and $\QQ_F$ dividing $\ell$ one has
$$
\lambda(p)\equiv p^{k-2} + a(p) +p^{j+k-1} \, (\bmod \, \ell')
$$
for all primes $p$.
\end{conjecture}

(Here the $\lambda(p)$ are algebraic integers lying in a totally
real field $\QQ_F$. Harder formulated the conjecture for the case
$L=\QQ$.)

For example, 
if $f=\Delta e_4 e_6\in S_{22}(\Gamma_1)$ is the unique normalized
cusp form of weight $22$ then 
$41 | \Lambda(f,14)$, so Harder predicts that
the space $S_{4,10}$ should contain a non-zero 
eigenform $F$ with eigenvalues $\lambda(p)$ satisfying 
$\lambda(p)\equiv p^8+a(p)+p^{13} \, (\bmod \, 41)$
for all $p$.  A mimimum consistency is that at least 
$\dim S_{4,10}(\Gamma_2) \neq 0$,
actually as it turns out this dimension is~$1$.
\end{section}

\begin{section}{Evidence for Harder's Conjecture}
Since we can calculate the trace of the Hecke operators $T(p)$
on the spaces $S_{j,k}(\Gamma_2)$ for all primes $p\leq 37$
(modulo the conjecture on the endoscopic contribution)
we can try to check the conjecture by Harder (and gain evidence
for the conjecture on the endoscopic contribution at the same time). 
As we just saw, the first case where we
have a `large' prime dividing a critical $L$-value is
the eigenform $f=\Delta e_{10} \in S_{22}(\Gamma_1)$ of weight $22$.
Here the prime $41$ divides the critical values $L(f,14)/\Omega^+$.
The conjecture predicts a congruence between the Fourier coefficients
of $f=\sum_{n=1}^{\infty} a(n)q^n$ and the eigenvalues $\lambda(p)$
of a form $F$ in the $1$-dimensional space $ S_{4,10}(\Gamma_{2})$. 
We give the tables with the eigenvalues
$a(p)$ of $f$ and $\lambda(p)$ of $F \in S_{4,10}(\Gamma_2)$
for the primes $p\leq 37$.

\smallskip
\vbox{
\bigskip\centerline{\def\quad{\hskip 0.6em\relax}
\def\quod{\hskip 0.5em\relax }
\vbox{\offinterlineskip
\hrule
\halign{&\vrule#&\strut\quod\hfil#\quad\cr
height2pt&\omit&&\omit&&\omit&\cr
& $p$ && $a(p)$ && $\lambda(p)$ &\cr
\noalign{\hrule}
  & $2$ &&  $- 288$ && $-1680$  & \cr
  & $3$ &&  $-128844$ &&$55080$ & \cr
  & $5$ &&  $21640950$ &&$-7338900$ & \cr
  & $7$ &&  $-768078808$ &&$609422800$  & \cr
  &$11$ &&  $-94724929188$ &&$25358200824$ & \cr  
  &$13$ &&  $-80621789794$ && $-263384451140$ & \cr
  &$17$ &&  $3052282930002$ && $-2146704955740$ & \cr
  &$19$ &&  $-7920788351740$ &&$43021727413960$ & \cr
  &$23$ &&  $-73845437470344 $ && $-233610984201360$ &\cr
  &$29$ &&  $-4253031736469010 $ && $-545371828324260$ & \cr
  &$31$ &&  $1900541176310432 $ && $830680103136064$ & \cr
  &$37$ &&  $22191429912035222 $ && $11555498201265580$& \cr
} \hrule}
}}
\bigskip
\begin{proposition} 
The congruence $\lambda(p)\equiv p^8+a(p)+p^{13}
(\bmod \, 41 )$ for the eigenvalues $\lambda(p)$ and $a(p)$
on $S_{4,10}(\Gamma_2)$ and $S_{22}(\Gamma_1)$ holds
for all primes $p\leq 37$.
\end{proposition}

In this way we can check Harder's conjecture for many cases
given in the tables below in the following sense. 
If both $\dim S_r(\Gamma_1)=1$ and $\dim S_{j,k}(\Gamma_2)=1$ and
if $\ell$ is a prime $> r$ dividing the critical $L$-value then we
checked the congruence 
$\lambda(p) -a(p)-p^{j+k-1}-p^{k-2} \equiv 0 (\bmod \, \ell)$
for all primes $p\leq 37$. In case  $\dim S_r(\Gamma_1)=2$
and  $\dim S_{j,k}(\Gamma_2)=1$ I checked that in the quadratic
field $\QQ(a(p))$ the expression 
$\lambda(p) -a(p)-p^{j+k-1}-p^{k-2}$
has a norm divisible by $\ell$ for all primes $p\leq 37$.
With a bit of additional effort one can check the congruence in
the real quadratic field. For example, take $r=24$ and let 
$$
f=\sum a(n)q^n =q-(54-12\sqrt{144169})\, q^2 + \ldots
$$
be a normalized eigenform in $S_{24}(\Gamma_1)$. In the quadratic field
$\QQ(\sqrt{144169})$ the prime $73$ splits as $\pi \cdot \pi^{\prime}$
with $\pi=(73,53+36\sqrt{144169})$. Let $\lambda(p)$ be the eigenvalue
under $T(p)$ of the generator of $S_{12,7}(\Gamma_2)$. Then we can check the congruence
$$
\lambda(p)\equiv p^5+a(p)+p^{18} \, (\bmod \, \pi)
$$
for all $p \leq 37$.

In case $\dim S_{j,k}(\Gamma_2)=2$  I can calculate the characteristic
polynomial $g$ of $T(2)$. In general this is an irreducible polynomial
$g$ of degree $8$ over $\QQ$. The corresponding number field $L$ possesses
just one subfield $L$ of degree $2$ over $\QQ$
and $g$ decomposes in two polynomials
of degree $4$ that are irreducible over $K$. I then checked that
the expression $\lambda(p) -a(p)-p^{j+k-1}-p^{k-2}$ has a norm in the composite
field $(\QQ(a(2)),K)$ which is divisible by our congruence prime $\ell$.

For example, we treat the case of the local system $V_{18,6}$
with $(\ell,m)=(18,6)$. The characteristic
polynomial $g$ of Frobenius at the prime $2$ is:
$$
1+t_1\,X+t_2\,{X}^{2}+t_3\,{X}^{3}+t_4
\,{X}^{4}+\\
2^{27}\,t_3{X}^{5}+2^{54}t_2\,{X
}^{6}+2^{81}t_1{X}^{7}+
2^{108}\,{X}^{8}.
$$
with $t_1=12432$, $t_2=193574912$, $t_3=3043199287296$ and $t_4=31380514975776768$.
The corresponding degree $8$ field extension $K$ of $\QQ$ has one
quadratic subfield $\QQ(\sqrt{7\cdot 3607})$. Our polynomial $g$ splits into
the product of a quartic polynomial~$h$
$$
\begin{matrix}
18014398509481984\,{X}^{4}+(834297397248-9663676416\sqrt {25249})\, X^3+\\
(142913536-110592\sqrt {25249})\,X^2+
(6216-72\sqrt {25249})\, X+1 \\
\end{matrix}
$$
and its conjugate over this quadratic subfield $\QQ(\sqrt{25249})$
and we get
$ 
\lambda(2)= -6216 \pm 72 \sqrt{25249}
$.
The normalized eigenform in $S_{28}$ has Fourier coefficient
$
a(2)=- 4140 \pm 108\sqrt{18209})
$
and one checks that the norm of 
$$
6216 +72\sqrt{25249}+ 2^7+2^{20}-(4140+108\sqrt{18209})
$$ 
in the field $\QQ(\sqrt{25249},\sqrt{18209})$ is divisible by $4057$
as predicted by Harder.

But there are cases where the characteristic polynomial $g$ decomposes.
These are the cases $(j,k)=(18,7)$ where we have two factors of degree $4$
and $(j,k)=(8,13)$ where $g$ is a product of four quadratic factors.
In the cases $(j,k)=(18,7)$ there is a congruence modulo $3779$.
In fact, $g$ decomposes as the product of
$$
288230376151711744 X^4 - 4252017623040 X^3 + 45752320X^2 - 7920X + 1
$$
and
$$
288230376151711744X^4 + 17575006175232X^3 + 857571328X^2 + 32736X + 1
$$
and one calculates
$$
{\rm Norm}(4320+96\sqrt{51349}+2^{24}+2^5+32736)=282720345772032
$$
and this is divisible by $3779$.
In the cases $(j,k,r)=(32,4,38)$ there are two congruence primes
and one finds indeed a congruence for both of them.

The following tables list the congruence primes in question.
All of these are checked in the sense explained above.

\vfill
\eject

{\small
\smallskip
\vbox{
\bigskip\centerline{\def\quad{\hskip 0.6em\relax}
\def\quod{\hskip 0.5em\relax }
\vbox{\offinterlineskip
\hrule
\halign{&\vrule#&\strut\quod\hfil#\quad\cr
height2pt&\omit&&\omit&&\omit&&\omit&&\omit&&\omit&\cr
&$ r$&&$\dim(S_r)$&&$(j,k)$&&$\dim(S_{j,k})$&&L-value&&primes&\cr
height2pt&\omit&&\omit&&\omit&&\omit&&\omit&&\omit&\cr
\noalign{\hrule}
height2pt&\omit&&\omit&&\omit&&\omit&&\omit&&\omit&\cr
&$20$&&$1$&&$(6,8)$&&$1$&&$2^2\cdot 3 \cdot 11$&&&\cr
height2pt&\omit&&\omit&&\omit&&\omit&&\omit&&\omit&\cr
\noalign{\hrule}
height2pt&\omit&&\omit&&\omit&&\omit&&\omit&&\omit&\cr
&$22$&&$1$&&$(4,10)$&&$1$&&$-2\cdot 3\cdot 17 \cdot 41$&&41&\cr
height2pt&\omit&&\omit&&\omit&&\omit&&\omit&&\omit&\cr
\noalign{\hrule}
height2pt&\omit&&\omit&&\omit&&\omit&&\omit&&\omit&\cr
&$22$&&$1$&&$(8,8)$&&$1$&&$3 \cdot 7 \cdot 13 \cdot 17$&&&\cr
height2pt&\omit&&\omit&&\omit&&\omit&&\omit&&\omit&\cr
\noalign{\hrule}
height2pt&\omit&&\omit&&\omit&&\omit&&\omit&&\omit&\cr
&$22$&&$1$&&$(12,6)$&&$1$&&$-2\cdot 7 \cdot 13^2$&&&\cr
height2pt&\omit&&\omit&&\omit&&\omit&&\omit&&\omit&\cr
\noalign{\hrule}
height2pt&\omit&&\omit&&\omit&&\omit&&\omit&&\omit&\cr
&$24$&&$2$&&$(12,7)$&&$1$&&$2^4\cdot 5 \cdot 7^2 \cdot 11 \cdot73$&&$73$&\cr
height2pt&\omit&&\omit&&\omit&&\omit&&\omit&&\omit&\cr
\noalign{\hrule}
height2pt&\omit&&\omit&&\omit&&\omit&&\omit&&\omit&\cr
&$24$&&$2$&&$(6,10)$&&$1$&&$3\cdot 11^2 \cdot 13^2 \cdot 17 $&&&\cr
height2pt&\omit&&\omit&&\omit&&\omit&&\omit&&\omit&\cr
\noalign{\hrule}
height2pt&\omit&&\omit&&\omit&&\omit&&\omit&&\omit&\cr
&$24$&&$2$&&$(8,9)$&&$1$&&$-2\cdot 7 \cdot 11 \cdot 29$&&$29$&\cr
height2pt&\omit&&\omit&&\omit&&\omit&&\omit&&\omit&\cr
\noalign{\hrule}
height2pt&\omit&&\omit&&\omit&&\omit&&\omit&&\omit&\cr
&$26$&&$1$&&$(4,12)$&&$1$&&$2 \cdot 11 \cdot 17 \cdot 19$&&&\cr
\noalign{\hrule}
height2pt&\omit&&\omit&&\omit&&\omit&&\omit&&\omit&\cr
&$26$&&$1$&&$(6,11)$&&$1$&&$3\cdot 5 \cdot 11 \cdot 19$&&&\cr
\noalign{\hrule}
height2pt&\omit&&\omit&&\omit&&\omit&&\omit&&\omit&\cr
&$26$&&$1$&&$(10,9)$&&$1$&&$-2\cdot 7 \cdot 11\cdot 29$&&$29$&\cr
\noalign{\hrule}
height2pt&\omit&&\omit&&\omit&&\omit&&\omit&&\omit&\cr
&$26$&&$1$&&$(14,7)$&&$1$&&$5\cdot 7 \cdot 97$&&$97$&\cr
\noalign{\hrule}
height2pt&\omit&&\omit&&\omit&&\omit&&\omit&&\omit&\cr
&$26$&&$1$&&$(16,6)$&&$1$&&$-2\cdot 11 \cdot 17 \cdot 19$&&&\cr
\noalign{\hrule}
height2pt&\omit&&\omit&&\omit&&\omit&&\omit&&\omit&\cr
&$26$&&$1$&&$(18,5)$&&$1$&&$-2^3\cdot 3 \cdot 43$&&$43$&\cr
\noalign{\hrule}
height2pt&\omit&&\omit&&\omit&&\omit&&\omit&&\omit&\cr
&$26$&&$1$&&$(8,10)$&&$2$&&$-3^2\cdot 7 \cdot 11\cdot 19$&&&\cr
\noalign{\hrule}
height2pt&\omit&&\omit&&\omit&&\omit&&\omit&&\omit&\cr
&$26$&&$1$&&$(12,8)$&&$2$&&$3\cdot 5^2\cdot 11 \cdot 17$&&&\cr
\noalign{\hrule}
height2pt&\omit&&\omit&&\omit&&\omit&&\omit&&\omit&\cr
&$28$&&$2$&&$(2,14)$&&$1$&&$2^3\cdot 5^2 \cdot  13^2\cdot 17^2\cdot 19 \cdot 23$&&&\cr
height2pt&\omit&&\omit&&\omit&&\omit&&\omit&&\omit&\cr
\noalign{\hrule}
height2pt&\omit&&\omit&&\omit&&\omit&&\omit&&\omit&\cr
&$28$&&$2$&&$(16,7)$&&$1$&&$2^5\cdot 3^4 \cdot 5\cdot 7 \cdot 13 \cdot 367$&&$367$&\cr
\noalign{\hrule}
height2pt&\omit&&\omit&&\omit&&\omit&&\omit&&\omit&\cr
&$28$&&$1$&&$(14,8)$&&$2$&&$2^4\cdot 11 \cdot 13^2 17 \cdot 19\cdot 23 \cdot 647$&&$647$&\cr
\noalign{\hrule}
height2pt&\omit&&\omit&&\omit&&\omit&&\omit&&\omit&\cr
&$28$&&$2$&&$(12,9)$&&$2$&&$2^3\cdot 7 \cdot 11\cdot 13 \cdot 23\cdot 4057$&&$4057$&\cr
\noalign{\hrule}
height2pt&\omit&&\omit&&\omit&&\omit&&\omit&&\omit&\cr
&$28$&&$2$&&$(8,11)$&&$1$&&$5\cdot 11^2 \cdot 13 \cdot 23\cdot 2027$&&$2027$&\cr\noalign{\hrule}
height2pt&\omit&&\omit&&\omit&&\omit&&\omit&&\omit&\cr
&$28$&&$2$&&$(18,6)$&&$1$&&$2^4\cdot 3^2\cdot 5^2 \cdot11 \cdot 13^2 \cdot 17^2 \cdot 19$&&&\cr
\noalign{\hrule}
height2pt&\omit&&\omit&&\omit&&\omit&&\omit&&\omit&\cr
&$28$&&$1$&&$(10,10)$&&$2$&&$2^2\cdot 5^2\cdot 11^2 \cdot 13^2 \cdot 17\cdot 23 \cdot 157$&&$157$&\cr
\noalign{\hrule}
height2pt&\omit&&\omit&&\omit&&\omit&&\omit&&\omit&\cr
&$28$&&$2$&&$(6,12)$&&$2$&&$5\cdot 11^2 \cdot 13^2\cdot 19 \cdot 23\cdot 823$&&$823$&\cr
\noalign{\hrule}
height2pt&\omit&&\omit&&\omit&&\omit&&\omit&&\omit&\cr
&$28$&&$2$&&$(20,5)$&&$1$&&$2^9\cdot 3^4 \cdot 5 \cdot 193$&&$193$&\cr
%
\noalign{\hrule}
height2pt&\omit&&\omit&&\omit&&\omit&&\omit&&\omit&\cr
&$30$&&$2$&&$(14,9)$&&$2$&&$2^8 \cdot 3 \cdot 5 \cdot 13 \cdot 1039$ && $1039$& \cr
\noalign{\hrule}
height2pt&\omit&&\omit&&\omit&&\omit&&\omit&&\omit&\cr
&$30$&&$2$&&$(6,13)$&&$1$&&$2^4\cdot 5 \cdot 11 \cdot 13 \cdot 19 \cdot 23$ &&& \cr
\noalign{\hrule}
height2pt&\omit&&\omit&&\omit&&\omit&&\omit&&\omit&\cr
&$30$&&$2$&&$(10,11)$&&$1$&&$3^4 \cdot 11 \cdot 13 \cdot 23 \cdot 97$ &&$97$& \cr
\noalign{\hrule}
height2pt&\omit&&\omit&&\omit&&\omit&&\omit&&\omit&\cr
&$30$&&$2$&&$(24,4)$&&$1$&&$2^{10}\cdot 3^4 \cdot 5^5 \cdot 7 \cdot 97$ &&$97$& \cr
\noalign{\hrule}
height2pt&\omit&&\omit&&\omit&&\omit&&\omit&&\omit&\cr
&$30$&&$2$&&$(20,6)$&&$2$&&$2^6 \cdot 3^3 \cdot 7 \cdot 11 \cdot 13 \cdot 17 \cdot 19 \cdot 23 \cdot 593$ &&$593$& \cr
\noalign{\hrule}
height2pt&\omit&&\omit&&\omit&&\omit&&\omit&&\omit&\cr
&$30$&&$2$&&$(4,14)$&&$2$&&$3^2 \cdot 5 \cdot 7^2 \cdot 13 \cdot 19^2 \cdot 23 \cdot 4289$ &&$4289$& \cr
\noalign{\hrule}
height2pt&\omit&&\omit&&\omit&&\omit&&\omit&&\omit&\cr
&$30$&&$2$&&$(18,7)$&&$2$&&$2^4 \cdot 3^2 \cdot 5 \cdot 11 \cdot 3779$ &&$3779$& \cr
\noalign{\hrule}
height2pt&\omit&&\omit&&\omit&&\omit&&\omit&&\omit&\cr
&$32$&&$2$&&$(4,15)$&&$1$&&$2^2 \cdot 5 \cdot 7^2 \cdot 13 \cdot 19 \cdot 23 \cdot 61$ &&$61$& \cr
\noalign{\hrule}
height2pt&\omit&&\omit&&\omit&&\omit&&\omit&&\omit&\cr
&$32$&&$2$&&$(2,16)$&&$2$&&$3^3\cdot 5^2 \cdot 7^2 \cdot 19^2 \cdot 23 \cdot 211$ &&$211$& \cr
\noalign{\hrule}
height2pt&\omit&&\omit&&\omit&&\omit&&\omit&&\omit&\cr
&$32$&&$2$&&$(22,6)$&&$2$&&$2^3\cdot 3^3 \cdot 5 \cdot 7 \cdot 13 \cdot 17 \cdot 19 \cdot 23 \cdot 7687$ &&$7687$& \cr
\noalign{\hrule}
height2pt&\omit&&\omit&&\omit&&\omit&&\omit&&\omit&\cr
&$32$&&$2$&&$(24,5)$&&$2$&&$2^{9} \cdot 3^5 \cdot 3119$ &&$3119$& \cr
\noalign{\hrule}
height2pt&\omit&&\omit&&\omit&&\omit&&\omit&&\omit&\cr
&$32$&&$2$&&$(8,13)$&&$2$&&$2 \cdot 7^3 \cdot 11^3 \cdot 13^2 \cdot 23$ &&& \cr
\noalign{\hrule}
height2pt&\omit&&\omit&&\omit&&\omit&&\omit&&\omit&\cr
&$34$&&$2$&&$(10,13)$&&$2$&&$2^3 \cdot 3^2 \cdot 5 \cdot 7 \cdot 13^2 \cdot 23^2 \cdot 29^2$ &&& \cr
\noalign{\hrule}
height2pt&\omit&&\omit&&\omit&&\omit&&\omit&&\omit&\cr
&$34$&&$2$&&$(28,4)$&&$1$&&$2^{10} \cdot 3^8 \cdot 5^5 \cdot 7 \cdot 103$ &&$103$& \cr
\noalign{\hrule}
height2pt&\omit&&\omit&&\omit&&\omit&&\omit&&\omit&\cr
&$34$&&$2$&&$(26,5)$&&$2$&&$2^{11} \cdot 3^3 \cdot 5^3 \cdot 15511$ &&$15511$& \cr
\noalign{\hrule}
height2pt&\omit&&\omit&&\omit&&\omit&&\omit&&\omit&\cr
&$34$&&$2$&&$(6,15)$&&$2$&&$2 \cdot 5^2 \cdot 7 \cdot 13 \cdot 23^2 \cdot 29 \cdot 233$ &&$233$& \cr
\noalign{\hrule}
height2pt&\omit&&\omit&&\omit&&\omit&&\omit&&\omit&\cr
&$38$&&$2$&&$(32,4)$&&$2$&&$2^8 \cdot 3^8 \cdot 5^4 \cdot 7^2 \cdot 67 \cdot 83$ &&$67$, $83$& \cr
} \hrule}
}}
\bigskip
} 

\end{section}

\eject

\printindex

  \end{document}